\documentclass[12pt]{article}
\usepackage{amsfonts}
\usepackage{epsfig}
\title{The Octagonal PET I: Renormalization and
Hyperbolic Symmetry}
\author{Richard Evan Schwartz \thanks{\hskip 5 pt Supported by 
N.S.F. Research Grant DMS-0072607}}

\newtheorem{theorem}{Theorem}[section]

\newtheorem{lemma}[theorem]{Lemma}

\newtheorem{corollary}[theorem]{Corollary}

\def\startproof{{\bf {\medskip}{\noindent}Proof: }}

\def\endproof{$\spadesuit$  \newline}

\def\C{\mbox{\boldmath{$C$}}}%
\def\H{\mbox{\boldmath{$H$}}}%
\def\N{\mbox{\boldmath{$N$}}}%
\def\Q{\mbox{\boldmath{$Q$}}}%
\def\R{\mbox{\boldmath{$R$}}}%
\def\Z{\mbox{\boldmath{$Z$}}}%

\begin{document}
\maketitle

\begin{abstract}
We introduce
a family of polytope exchange transformations
acting on parallelotopes in $\R^{2n}$,
for $n=1,2,3,...$ These
PETs are constructed using a pair of lattices
in $\R^{2n}$. The moduli
space of these PETs is $GL_n(\R)$.  
We study the case $n=1$ in detail.
In this case, we show that the
$2$-dimensional family is completely
renormalizable and that the $(2,4,\infty)$ 
hyperbolic reflection triangle group acts 
(by linear fractional transformations)
as the renormalization group on the moduli space.
These results have a number of geometric
corollaries for the system.
\end{abstract}

\section{Introduction}

\subsection{Background}

A {\it polytope exchange transformation\/} (or PET) is
defined by a polytope $X$ which has been partitioned
in two ways into smaller polytopes:
$$X=\bigcup_{i=1}^m A_i=\bigcup_{i=1}^m B_i.$$
For each $i$ there is some vector $V_i$ such that
$B_i=A_i+V_i$.  That is, some translation carries
$A_i$ to $B_i$.  One then defines a map $f: X \to X$
by the formula $f(x)=x+V_i$ for all $x \in {\rm int\/}(A_i)$.
It is understood that $f$ is not defined for points
in the boundaries of the small polytopes.  The
inverse map is defined by $f^{-1}(y)=y-V_i$ for all
$y \in {\rm int\/}(B_i)$.

The simplest examples of PETs are $1$-dimensional
systems, known as {\it interval exchange transformations\/}
(or IETs).
These systems have been extensively studied
in the past $30$ years, and there are close connections
between IETs and other areas of mathematics such
as Teichmuller theory.
See, for instance, [{\bf Y\/}] and [{\bf Z\/}] and the many references
mentioned therein.

The {\it Rauzy renormalization\/} [{\bf R\/}]
gives a satisfying renormalization theory for
the family of IETs all having the same number of
intervals in the partition.
The idea is that one starts with an $n$-interval IET, and then
considers the first return map to a specially chosen
sub-interval.  This first return map turns out to give
another $n$-interval IET.  This mechanism sheds a lot of
light on IETs.
Some examples of polygon exchange maps have been
studied in [{\bf AG\/}], [{\bf AKT\/}], [{\bf H\/}],
[{\bf Hoo\/}],
[{\bf LKV\/}],
[{\bf Low\/}], [{\bf S2\/}],
[{\bf S3\/}], and [{\bf T\/}].
Some definitive theoretical work concerning the
(zero) entropy of such systems is done in
[{\bf GH1\/}], [{\bf GH2\/}], and [{\bf B\/}].

Often, renormalization phenomena
are found in these systems:  The first return to
some subset is conjugate to the original
system, and this allows for a detailed understanding
of the system.  
In [{\bf S1\/}] and [{\bf Hoo\/}], a 
renormalizable family of polygon exchange maps is
constructed. In this setting, like
in Rauzy renormalization, the first return map to
a subset of one system is conjugate to the first
return map of another system in the family.  
The family in [{\bf Hoo\/}] is $2$-dimensional
and the family in [{\bf S1\/}] is $1$-dimensional.

In [{\bf S1\/}] we introduced a class of
PETs which we
called {\it double lattice PETs\/}. See \S 2
for a definition. These
PETs are defined in terms of a pair of
Euclidean lattices, and in each dimension
there is a large space of them.
We originally found (some of) 
the double lattice PETs as
compactifications of polygonal outer billiards systems,
but it seems reasonable to study these objects
for their own sake. One motivation for studying
these objects is to search for good
examples of renormalization schemes for
families of higher dimensional PETs.

In this paper we introduce a family
of double lattice PETs in every even dimension.
In dimension $2n$, the objects
are indexed by $GL_n(\R)$.  In the
two-dimensional setting, there is a $1$-parameter
family.  In this paper we will study the
$1$-parameter family of $2$-dimensional
examples in detail, showing that it has
a complete renormalization scheme. The $(2,4,\infty)$
reflection triangle group acts on the parameter space,
and points in the same orbit have closely related
dynamics - e.g., their limit sets have the same
Hausdorff dimension.

\subsection{The Planar Construction}

We will describe the $2n$-dimensional examples
systematically in \S \ref{higher}.  Here we explain the
$2$-dimensional case in a visual way.  Our construction
depends on a parameter $s \in (0,\infty)$.  Usually, but
not always, we take $s \in (0,1)$.  Below we suppress $s$
from most of our notation.

\begin{center}
\resizebox{!}{2.4in}{\includegraphics{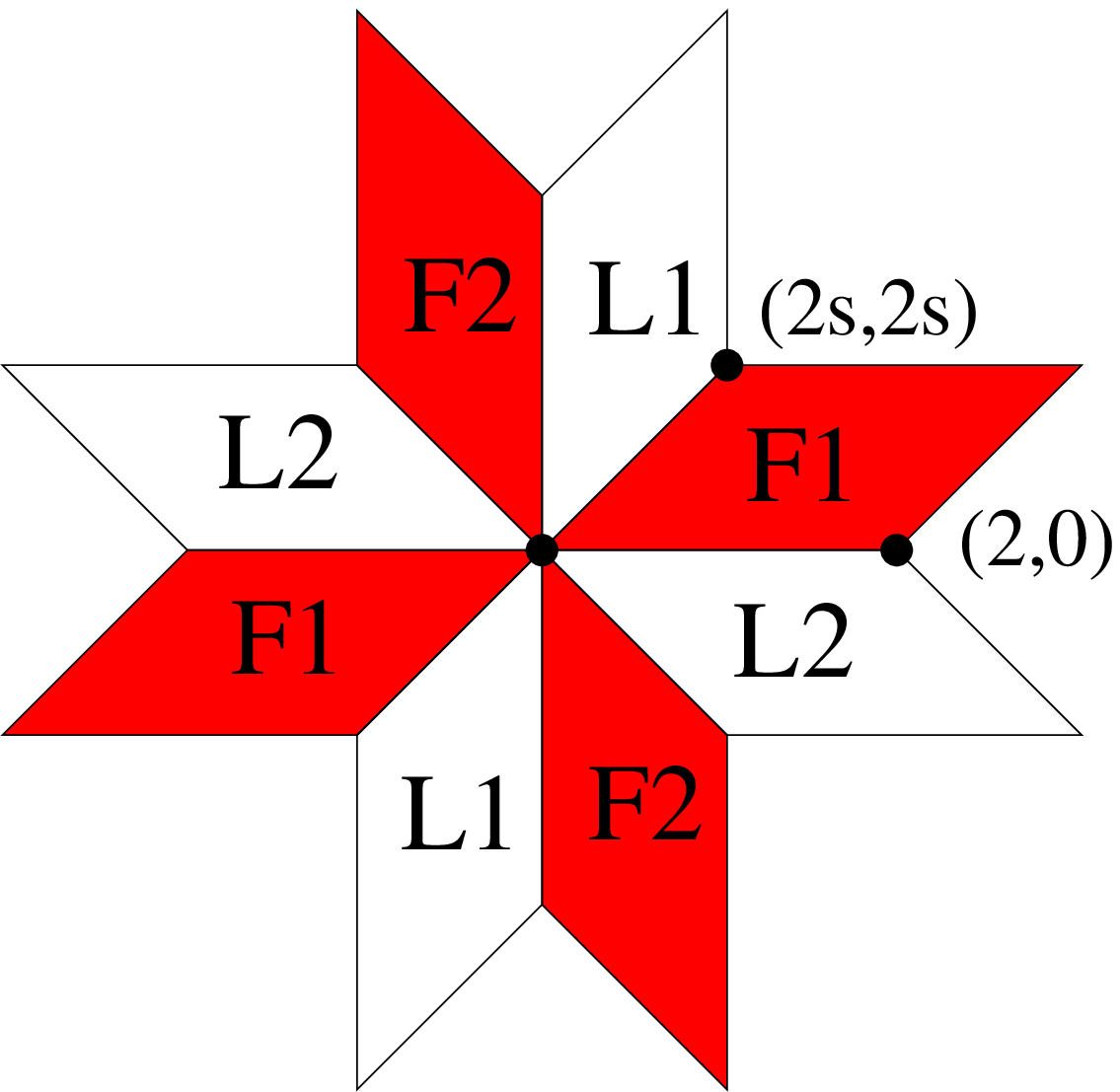}}
\newline
{\bf Figure 1.1:\/} The scheme for the PET.
\end{center}

The $8$ parallalograms in Figure 1.1 
are the orbit of a single parallelogram $P$ under a dihedral group
of order $8$.  Two of the sides of $P$ are
determined by the vectors
$(2,0)$ and $(2s,2s)$. For $j=1,2$, let $F_j$ denote the parallelogram
centered at the origin and translation equivalent to the ones
in the picture labeled $F_j$.  Let $L_j$ denote the lattice
generated by the sides of the parallelograms labeled $L_j$.
(Either one generates the same lattice.)

In \S 2 we will check the easy fact that
$F_i$ is a fundamental domain for
$L_j$, for all $i,j \in \{1,2\}$.
We define a system $(X',f')$, with $X'=F_1 \cup F_2$,
and $f': X' \to X'$, as follows.
Given $p \in F_j$ we let
\begin{equation}
f'(p)=p+V_p \in F_{3-j}, \hskip 30 pt
V_p \in L_{3-j}.
\end{equation}
The choice of $V_p$ is almost always
unique, on account of $F_{3-j}$ being
a fundamental domain for $L_{3-j}$.
When the choice is not unique,
we leave $f'$ undefined.
When $p \in F_1 \cap F_2$ we have
$V_p=0 \in L_1 \cap L_2$.
We will show in \S 2 that $(X',f')$ is a
PET.  

We prefer the map
$f=(f')^2$, which preserves both
$F_1$ and $F_2$. We set $X=F_1$.
Our system is
$f: X \to X$, which we denote by $(X,f)$.

\subsection{The Tiling and the Limit Set}

We define a {\it periodic tile\/} for $(X,f)$
as a maximal
convex polygon on which $f$ and its iterates
are completely defined and periodic.  So,
every point in a periodic tile has the same
period and every periodic point is contained
in a nontrivial periodic tile.
We call the union $\Delta$ of the periodic
tiles the {\it tiling\/}.

We define the
{\it aperiodic set\/} 
$\Lambda$ to be the set of aperiodic
points of $f$. The aperiodic set is
a subset of a somewhat more natural
set which we call the {\it limit set\/}
and denote by $\widehat \Lambda$.
The set $\widehat \Lambda$ is the
set of weakly aperiodic points.
  We call a point
$p \in X$ {\it weakly aperiodic\/} if
there is a sequence
$\{q_n\}$ converging to $p$ with the
following property. The first iterates
of $f$ are defined on $q_n$ and
the points $f^k(q_n)$ for $k=1,...,n$
are distinct.  
When $\Delta$ is dense (and it turns out
that this always happens for our system) we can say
alternately that  $\widehat \Lambda$
consists of those points $p$ such that
every neighborhood of $p$ intersects
infinitely many tiles of $\Delta$.  See
Lemma \ref{newdef}.  One advantage 
$\widehat \Lambda$ has over $\Lambda$ is
that $\widehat \Lambda$ is compact.

Figures 1.2 and 1.3 show (approximations of)
$\Delta$ and $\widehat \Lambda$ for two
quadratic irrational parameters.
In Figure 1.2, 
the picture is the same (locally) as what one sees
for outer billiards on the regular octagon.
In Figure 1.3, $\widehat \Lambda$ is the
union of two curves (though we do not give
a proof in this paper).
Isometric copies of these
curves arise in Pat Hooper's system
[{\bf Hoo\/}]. Hooper and I plan to explore this
``coincidence'' later.

\begin{center}
\resizebox{!}{1.3in}{\includegraphics{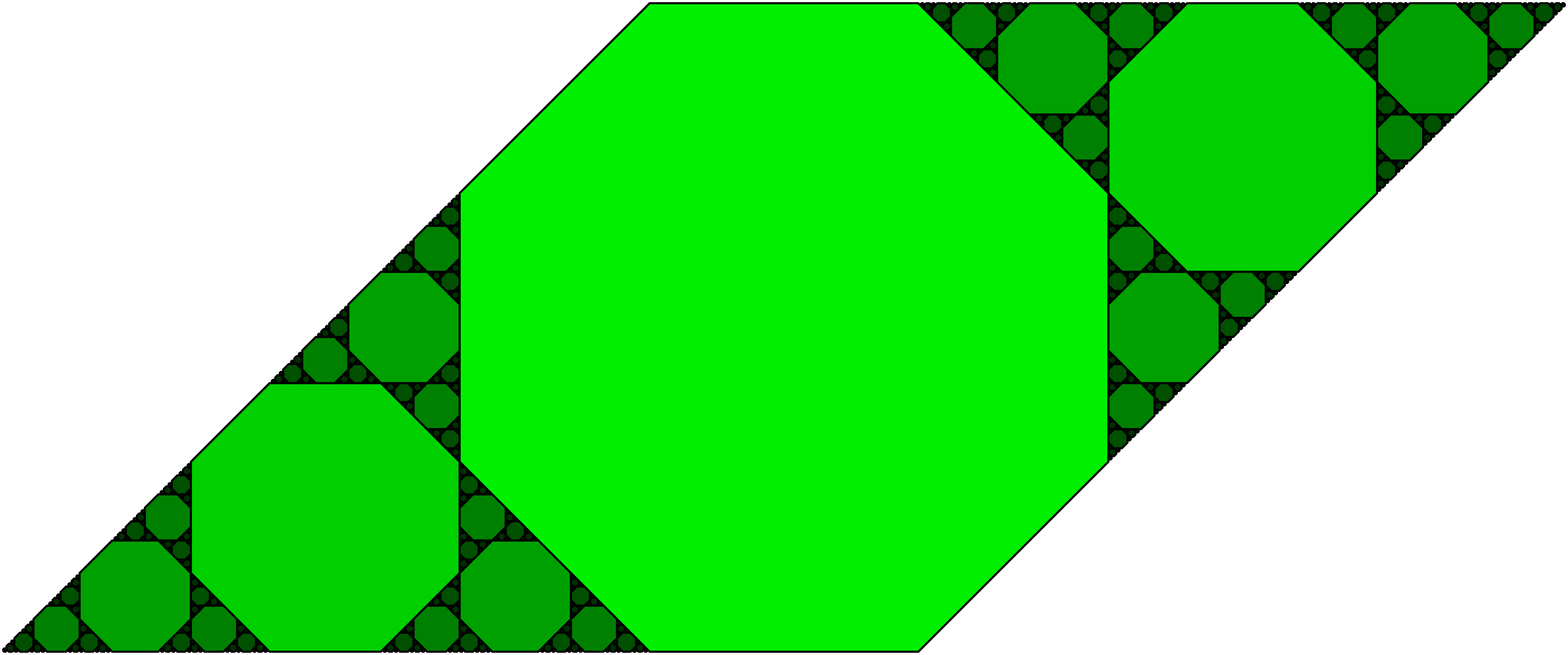}}
\newline
{\bf Figure 1.2:\/} The 
tiling associated to 
$s=\sqrt 2/2$
\end{center}

\begin{center}
\resizebox{!}{1.3in}{\includegraphics{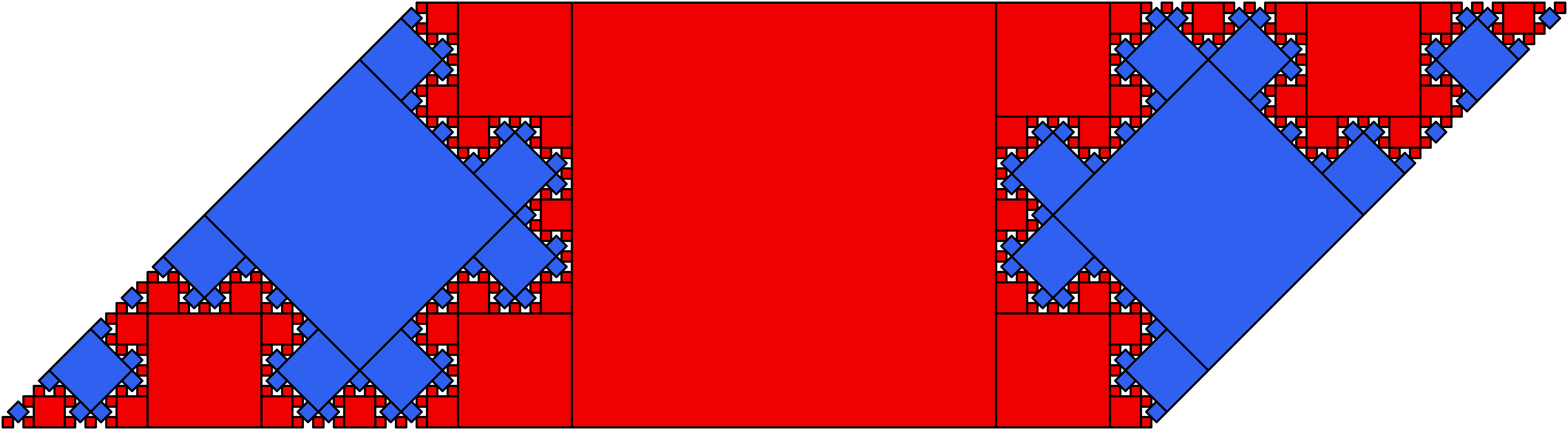}}
\newline
{\bf Figure 1.3:\/} The tiling associated to 
$s=\sqrt 3/2-1/2$.
\end{center}

Here are some results about the tiling.

\begin{theorem}
\label{cor}
When $s$ is rational, $\Delta_s$ is a finite
union of squares, semi-regular octagons, and
right-angled isosceles triangles.  When
$s$ is irrational, $\Delta_s$ is an infinite
union of squares and semi-regular octagons.
Moreover, the following is true.
\begin{enumerate}
\item $\Delta_s$ has at least one square
unless $s=\sqrt 2/2$.
\item $\Delta_s$ has only squares
if and only if the continued fraction
expansion of $s$ has the form
$[a_0,a_1,a_2,a_3,..]$ where $a_k$ is
even for all odd $k$.
\item $\Delta_s$ has infinitely many squares and
a dense set of shapes of semi-regular octagons,
for almost all $s$.
\end{enumerate}
\end{theorem}

A {\it semi-regular octagon\/} is an octagon with
$8$-fold dihedral symmetry.
Theorem \ref{cor} is a corollary of a more precise
statement about $\Delta_s$, Theorem \ref{one} below.
We defer the statement of Theorem \ref{one} because
it requires a build-up of terminology.

Here are some results about the limit set and the aperiodic set.
We only care about the irrational case.  These sets are
empty when $s$ is rational. See Lemma \ref{rat}.

\begin{theorem}
\label{two}
Suppose $s$ is irrational.
\begin{enumerate}
\item $\widehat \Lambda_s$ has zero area.
\item The projection of $\widehat \Lambda_s$
onto a line parallel to any $8$th root of unity
contains a line segment. Hence $\widehat \Lambda_s$
has Hausdorff dimension at least $1$.
\item $\widehat \Lambda_s$ is not contained in a
finite union of lines.
\item $\widehat \Lambda_s-\Lambda_s$ has zero
length for almost all $s$.  Hence
$\Lambda_s$ has Hausdorff dimension at least
$1$ for almost all $s$.
\end{enumerate}
\end{theorem}

\noindent
{\bf Remark:\/} 
In \S \ref{intersectX} we prove a more precise
result about when $\widehat \Lambda_s-\Lambda_s$
has zero length. \newline

The results above are consequences of 
the remormalization properties of
the family $\{(X_s,f_s)|\ s \in (0,1)$. We
explain this next.

\subsection{Renormalization}

We define the
{\it renormalizaton map\/}
$R: (0,1) \to [0,1)$ as follows.
\begin{itemize}
\item $R(x)=1-x$ if $x>1/2$
\item $R(x)=1/(2x)-{\rm floor\/}(1/(2x))$
if $x<1/2$.
\end{itemize}
$R$ relates to the $(2,4,\infty)$ reflection
triangle triangle much in the way that the classical
Gauss map relates to the modular group.

Define
\begin{equation}
Y=F_1-F_2 = X-F_2 \subset X.
\end{equation}

For any subset $S \subset X$, let
$f|S$ denote the first return map
to $S$, assuming that this map is
defined.  When we use this notation,
it means implicitly that the map is
actually defined, at least away
from a finite union of line segments.
We call $S$
{\it clean\/} if no point on $\partial S$ has
a well defined orbit.  This means, in particular,
that no tile of $\Delta$ crosses over
$\partial S$.

\begin{theorem}[Main]
Suppose $s \in (0,1)$ and $t=R(s) \in (0,1)$.
There is a clean set $Z_s \subset X_s$
such that $f_t|Y_t$ is conjugate to
$f_s^{-1}|Z_s$ by a map $\phi_s$.
\begin{enumerate}
\item $\phi_s$ commutes with reflection in the origin and
maps the acute vertices of $X_t$ to the
acute vertices of $X_s$.
\item  When $s<1/2$, the restriction
of $\phi_s$ to each component of $Y_t$
is an orientation reversing similarity,
with scale factor $s \sqrt 2$.
\item When $s<1/2$, either half of
$\phi_s$ extends to the trivial tile of
$\Delta_t$ and maps it to a tile in $\Delta_s$. 
\item When $s<1/2$, the only nontrivial
orbits which miss $Z_s$ are
contained in the 
$\phi_s$-images of the trivial tile of
$\Delta_t$.
These orbits have period $2$.
\item When $s>1/2$ the restriction
of $\phi_s$ to each component of $Y_s$
is a translation.
\item When $s>1/2$, all
nontrivial orbits intersect $Z_s$.
\end{enumerate}
\end{theorem}

The Main Theorem is an example of a result where
a picture says a thousand words.
Figures 1.4 and 1.5 show the Main Theorem in action for
$s<1/2$.  Figures 1.6 and 1.7 show the Main Theorem in
action for $s>1/2$.

\begin{center}
\resizebox{!}{1.35in}{\includegraphics{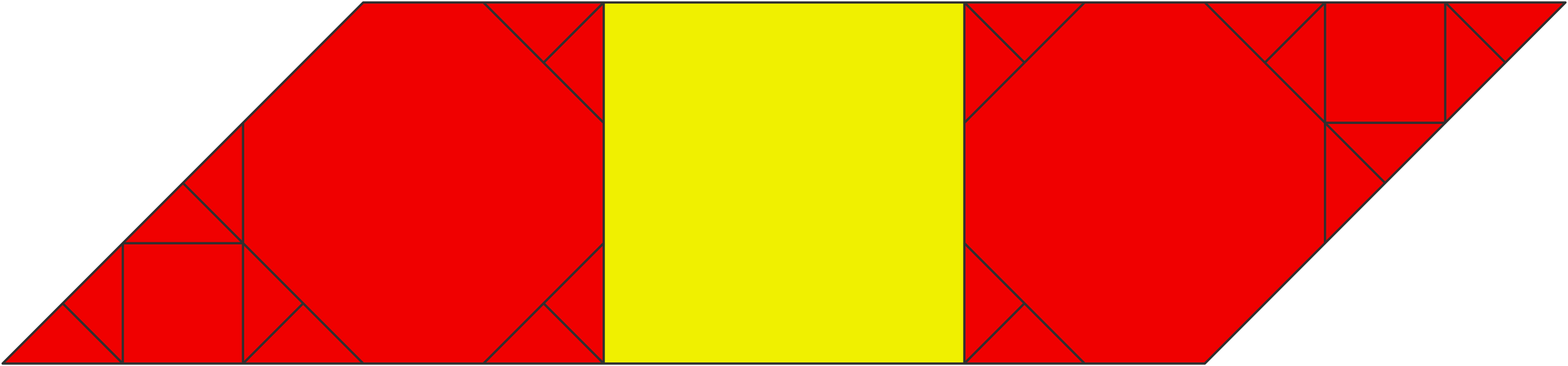}}
\newline
{\bf Figure 1.4:\/} $Y_t$ in red for $t=3/10=R(5/13)$.
\end{center}

\begin{center}
\resizebox{!}{1.4in}{\includegraphics{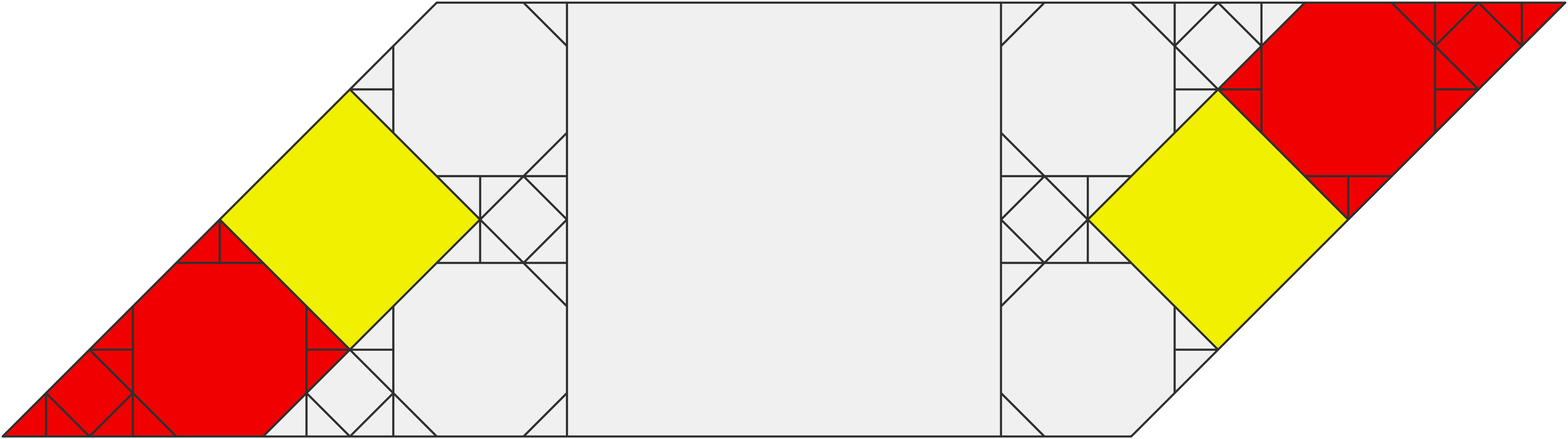}}
\newline
{\bf Figure 1.5:\/} $Z_s$ in red $s=5/13$.
\end{center}

\begin{center}
\resizebox{!}{1.35in}{\includegraphics{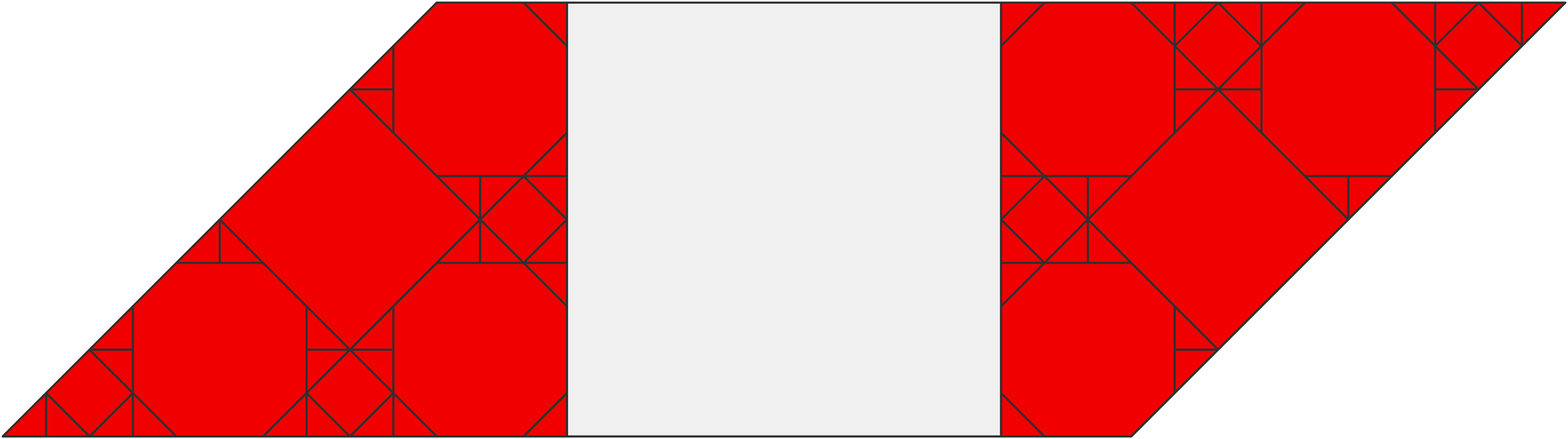}}
\newline
{\bf Figure 1.6:\/} $Y_t$ in red for $t=R(8/13)=5/13$.
\end{center}

\begin{center}
\resizebox{!}{2in}{\includegraphics{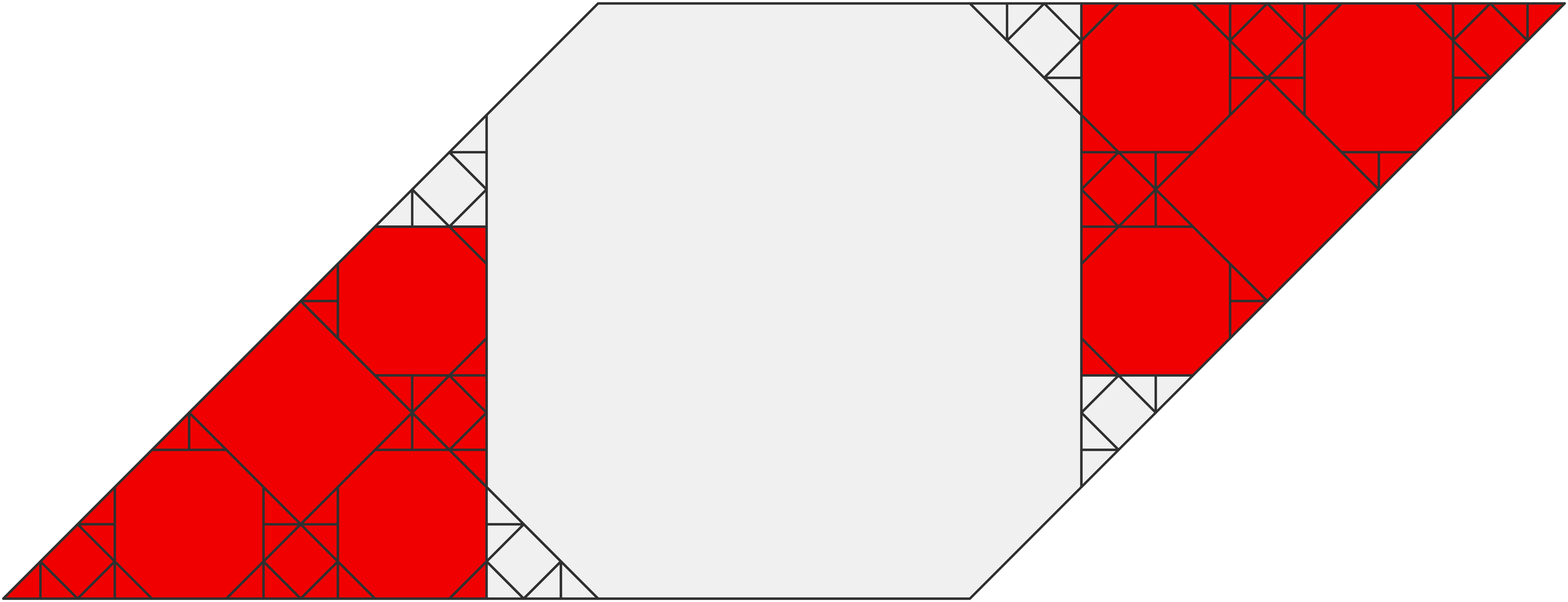}}
\newline
{\bf Figure 1.7:\/} $Z_s$ in red for $s=8/13$.
\end{center}

Now we explain the result behind Theorem \ref{cor}.
When $s>1/2$, the intersection
\begin{equation}
\label{fund}
O_s=(F_1)_s \cap (F_2)_s
\end{equation}
is the semi-regular octagon with vertices
\begin{equation}
\label{octagon}
(\pm s, \pm (1-s)), \hskip 30 pt
(\pm (1-s), \pm s).
\end{equation}
When $s<1/2$, the interection $O_s$ is
the square with vertices $(\pm s,\pm s)$.

Suppose we fix some irrational $s=s_0 \in (0,1)$.
Let $s_n=R^n(s)$.  Let $T_0$ be the identity map
and, referring to the Main Theorem, let
$T_n$ be the linear part of the composition
\begin{equation}
\phi_{s_0} \circ ... \circ \phi_{s_{n-1}}.
\end{equation}
The map $T_n$ is a similarity whose exact nature can
be computed using the information given in the
Main Theorem.

\begin{theorem}
\label{one}
When $s \in (0,1)$ is irrational,
a polygon arises in $\Delta_s$ if and only if it is
translation equivalent to
$T_n(O_{s_n})$ for some $n=0,1,2,...$. 
\end{theorem}

\noindent
{\bf Remark:\/}
When $s$ is rational, we get a very similar result, except
that $\Delta_s$ also contains some right-angled
isosceles triangles.  See Lemma \ref{st1}.

\subsection{Hyperbolic Symmetry}

Let $\H^2 \subset \C$ denote
the upper half plane model of the hyperbolic
plane. Let $\Gamma$ denote the
$(2,4,\infty)$ reflection triangle group,
generated by reflections in
the sides of the hyperbolic
triangle with vertices
\begin{equation}
\frac{i}{\sqrt 2}, \hskip 30 pt
\frac{1}{2}+\frac{i}{2}, \hskip 30 pt
\infty.
\end{equation}
We extend our parameter range so that
our system is defined for all $s \in \R$.
The systems at $s$ and $-s$ are identical.
$\Gamma$ acts on the parameter
set by linear fractional transformations.

We call the two systems $(X_s,f_s)$ and
$(X_t,f_t)$ {\it locally equivalent\/}
if the following is true.  There is
a finite union $L$ of lines such that,
for each point $p_s \in \widehat \Lambda_s-L$,
there is a point $p_t \in \widehat \Lambda_t$, together
with open neighborhoods $U_s$ and
$U_t$ of $p_s$ and $p_t$ respectively,
such that $\Delta_s \cap U_s$ is equivalent to
$\Delta_t \cap U_t$ by a similarity.
We also require the same statement to
be true with the roles of $s$ and $t$
reversed.  

Local equivalence is
strong:  For instance, the limit sets
of locally equivalent systems have the
same Hausdorff dimension. 

\begin{theorem}
\label{three}
Suppose $s$ and $t$ are in the
same orbit of $\Gamma$.  Then 
$(X_s,f_s)$ and $(X_t,f_t)$
are locally equivalent.  In particular,
the Hausdorff dimension of the limit
set, as a function of the parameter,
is a $\Gamma$-invariant function.
\end{theorem}

\noindent
{\bf Remarks:\/} \newline
(i) $\Gamma$ is contained with index $4$
in the group generated by reflections
in the ideal triangle with vertices
$0,1,\infty$. Using this fact, together with
a classic result about continued fractions, we
we will show that the forward orbit
$\{R^n(s)\}$ is dense in $(0,1)$ for almost
all $s \in (0,1)$.
\newline
(ii)
The need to exempt
a finite union of lines in the definition
of local equivalence seems partly
to be an artifact of our proof, but
in general one needs to disregard
some points to make everything work.
\newline
(iii) 
Given the ergodic nature of the action
of $\Gamma$, we can say that there is
some number $\delta_0$ such that
$\dim(\widehat \Lambda_s)=\delta_0$
for almost all $s$.  However, we don't
know the value of $\delta_0$.

\subsection{Further Results and Claims}

This paper now has a sequel [{\bf S0\/}].  Here
is the main result of that paper.

\begin{theorem}
Let $s \in (0,1)$ be irrational.
\begin{enumerate}
\item $\widehat \Lambda_s$ is a disjoint union of two
arcs if and only if $\Delta_s$ contains only squares.
(This happens if and only if $R^n(s)<1/2$ for all $n$.)
\item $\widehat \Lambda_s$ is a finite forest if and only
if $\Delta_s$ contains finitely many octagons.
(This happens if and only if $R^n(s)>1/2$ for finitely many $n$.)
\item $\widehat \Lambda_s$ is a Cantor set if and only if
$\Delta_s$ contains infinitely many octagons.
(This happens if and only if $R^n(s)>1/2$ for infinitely many $n$.)
\end{enumerate}
\end{theorem}

It turns out that our system here is
intimately related to outer billiards on semi-regular
octagons.  Recall that $O_s$ is the octagon from
Equation \ref{octagon}.
\newline
\newline
\noindent
{\bf Claim 2:\/}
\label{outer}
{\it
For any $s \in (1/2,1)$, the
limit set and perodic tiling produced by
outer billiards on $O_s$ are
locally isometric to
$\widehat \Lambda_s$ and $\Delta_s$,
respectively, except at finitely many points.
\/}
\newline
\newline
This claim is amply supported
by computer evidence, and I basically
know how to prove it.  I hope to
prove this claim in a sequel paper.

\subsection{Organization}

This paper is organized as follows.
In \S 2 we will define our PETs in
every even dimension and give some
basic information about them.  

In \S 3 we prove some basic results
about our PETs, most of which have
to do with the stability of orbits
under the perturbation of the parameter.

In \S 4 we prove some
some symmetry results about our PETs,
modulo $2$ computer calculations.

In \S 5 we prove the Main Theorem, modulo
$6$ more computer calculations.

In \S 6 we prove 
Theorem \ref{one} and Theorem \ref{cor}.

In \S 7 we deduce some length and area
estimates for the limit set.

In \S 8 we prove Theorem \ref{three} and
Theorem \ref{two}.

In \S 9 we will do the $8$ calculations
left over from \S 4-5.  These are exact
integer arithmetic calculations.

In \S 10 we list the raw data needed
for the calculations done in \S 9.

\subsection{OctaPET}

This paper has a companion java program,
called {\it OctaPET\/}. 
I discovered practically everything in
the paper while developing and using this program.
Also, the calculations
mentioned in \S 9 are all done using OctaPET.
I encourage you to use OctaPET while reading
the paper.  OctaPET relates to this paper much like a
song relates to musical notes written on a page.
You can download OctaPET from
\newline
{\bf http://www.math.brown.edu/$\sim$res/Java/OCTAPET.tar\/}
\newline
This is a tarred directory, which untars to a
directory called OctaPET.  This directory
contains a java program which you can 
compile and then run.

\subsection{Acknowledgements}

I would like to thank Nicolas Bedaride,
Pat Hooper, Injee Jeong,
John Smillie, and
Sergei Tabachnikov for interesting
conversations about topics related to
this work.  Some of this work
was carried out at ICERM in Summer 2012,
and some was
carried out during my sabbatical in
2012-13.
This sabbatical was funded from many sources.
I would like to thank the National
Science Foundation, All Souls College, Oxford,
the Oxford Maths Institute,
the Simons Foundation, the Leverhulme Trust, and
Brown University for their support during this time
period.

\newpage
\section{The Examples}

\subsection{Double Lattice PETs}

Generalizing the construction made in the introduction,
we will define what we mean by a double lattice PET.
The input to an $n$-dimensional
double lattice PET is a quadruple
$(L_1,L_2,F_1,F_2)$, where
\begin{itemize}
\item $L_1$ and $L_2$ are lattices in $\R^n$.
\item $F_1$ and $F_2$ are parallelotopes.
\item $F_i$ is a fundamental domain for
$L_j$ for all $i,j \in \{1,2\}$. 
\end{itemize}

Given $(F_1,F_2,L_1,L_2)$ we have a PET
defined on $X'=F_1 \cup F_2$ as follows.  For each
$x \in F_i$ we define (if possible)
$f'(x)=x+V_x$, where $V_x \in F_{3-i}$ is the unique vector
such that $x+V_x \in F_{3-i}$.  It may happen that
$V_x$ is not uniquely defined. In these cases,
we leave $f$ undefined. Note that $V_x=0 \in L_1 \cap L_2$ when
$x \in F_1 \cap F_2$, so that our definition is not ambigious
for such points. 

\begin{lemma}
$(X',f')$ is a PET.
\end{lemma}

\startproof
First, we check that $f$ is invertible.
The inverse map $(f')^{-1}$ is defined 
just as $f'$ is defined, but 
with the roles of $L_1$ and $L_2$
reversed: For each
$x \in F_i$ we define (if possible)
$f'(x)=x+W_x$, where $W_x \in F_{i}$ is the unique vector
such that $x+W_x \in F_{3-i}$.

The set $X'$ is partitioned as follows. Let
$U \subset X'$ denote the set of points for
which the assignment $p \to V_p$ is defined.
The set $X'-U$ is contained in a finite
union of hyperplanes.  Let $U' \subset X'$
denote the complement of these hyperplanes.
Then $U'$ is a finite union of convex
polytopes which partitions $X'$.

The assignment $p \to V_p$ is constant
on each component of $U'$.
Hence $V'=f'(U')$ is also a finite union
of polytopes.  But construction, the
polytopes in $V'$ have disjoint interior,
and their union has full measure.  Hence
$V'$ is a second partition of $X'$.
The two partitions $U'$ and $V'$ determine
$f$ in the manner of a PET.
\endproof

As in the introduction, it is convenient to set
$f=(f')^2$ and $X=F_1$.  Then
$(X,f)$ is also a PET, and the domain of
$f$ is the parallelotope $X$.

\subsection{Basic Construction}
\label{higher}

In this section we construct examples of
double lattice PETs in every even dimension.
To describe the examples in a natural way,
it is useful to work in $\C^n$,
complex $n$-space. Note, however, that when
it comes time to do calculations in the
case of interest, $\C$, we will revert to
working in $\R^2$.

Let $J$ denote multiplication by $i$.
Consider the following totally real subspaces
\begin{equation}
H=\{(z_1,...,z_n)|\ \Im(z_j)=0\ \forall j\} \hskip 12 pt
D=\{(z_1,...,z_n)|\ \Im(z_j)=\Re(z_j)\ \forall j\}.
\end{equation}
$H$ and $D$ respectively are the fixed point
sets of the reflections
\begin{equation}
R_H(z_1,...,z_n)=(\overline z_1,...,\overline z_n)
\hskip 12 pt
R_D(z_1,...,z_n)=(i\overline z_1,...,i\overline z_n).
\end{equation}
Clearly $J=R_D \circ R_H$.

Let $\{H_1,...,H_n\}$ and
$\{D_1,...,D_n\}$ be any $\R$-bases for $H$ and $D$
respectively.  Let $F_1$ be the parallelogram
centered at the origin, whose sides are
spanned by the vectors $\{H_1,...,H_n,D_1,...,D_n\}$.
We let $F_2=J(F_1)$.  Since
$J^2=-I$, and $F_j$ is centrally symmetric,
we have $J(F_2)=F_1$.

We let $L_1$ be the lattice spanned by
the vectors
$$H_1,...,H_n,R_H(D_1),...,R_H(D_n).$$
We let $L_2$ be the lattice spanned by
the vectors
$$D_1,...,D_n,R_D(H_1),...,R_D(H_n).$$

\begin{lemma}
$J(L_1)=L_2$ and $J(L_2)=L_1)$.
\end{lemma}

\startproof
Since $J^2=-{\rm Identity\/}$, we
have $J^2(L_j)=L_j$ for $j=1,2$.
So, it suffices to prove that
$J(L_1)=L_2$.
We have
$$J(H_j)=R_DR_H(H_j)=R_D(H_j) \in L_2,$$
$$J(R_H(D_j))=R_DR_HR_H(D_j)=R_D(D_j)=D_j \in L_2.$$
This does it.
\endproof

In the next result we will use the notation of
$\C^n$ but we point out in advance that our
argument only uses the $\R$-structure of $\C^n$.

\begin{lemma}
$F_1$ is a fundamental domain for both $L_1$ and $L_2$.
\end{lemma}

\startproof
The proof works the same way for both $L_1$ and $L_2$.
So, we will just prove that $F_1$ is a fundamental domain
for $L_1$.    Let $L_1'$ denote the lattice
generated by the sides of $F_1$.  A $\Z$-basis
for $L_1'$ is $\{H_1,...,H_n,D_1,....,D_n\}$.
Let $L_1'' \subset L_1 \cap L_1'$ be the $\Z$-span of
$\{H_1,...,H_n\}$.

Note that
$L_1$ and $L_1'$ have the same co-volume,
by symmetry.
But the volume of $F_1$ coincides with the
co-volume of $L_1'$.  Hence,
the volume of $F_1$ coincides with the
co-volume of $L_1$.  To finish the proof
we just have to show the following:
For any $p \in \C^n$ there is some
vector $V \in L_1$ such that
$p+V \in F_1$.

Let $\pi: \C^n \to H^{\perp}$ be orthogonal
projection.  The kernel of $\pi$ is exactly $H$.
Moreover, $\pi(F_1)$ is a fundamental domain for
$\pi(L_1')$ in $H^{\perp}$.   However $\pi(L_1)=\pi(L_1')$. So,
$\pi(F_1)$ is a fundamental domain for $\pi(L_1)$
as well.  Since $\pi(L_1)$ is
a fundamental domain for $\pi(F_1)$, there is
some $V_1 \in L_1$ such that the translate $H'$ 
of $H$, through $p+V_1$, intersects $F_1$.

Since $F_1$ is a parallelotope which intersects
$H$ in a parallelotope, the
intersection $H' \cap F_1$ is isometric
to $H \cap F_1$, and hence is a fundamental
domain for $L_1''$.  Hence, there is some
$V_2 \in L_1''$ such that $p+V_1+V_2 \in F_1$.
But $V=V_1+V_2\in L_1$.   So, $p+V \in F_1$.
This proves what we want.
\endproof

Since $J$ swaps $F_1$ with $F_2$ and
also swaps $L_1$ with $L_2$, the preceding
lemma shows that $F_i$ is a fundamental
domain for $L_j$ for all $i,j \in \{1,2\}$.
Now we know that the quaduple
$(F_1,F_2,L_1,L_2)$ defines a double
lattice PET.  We imagine that these
higher dimensional examples are interesting,
but so far the $2$-dimensional case is
hard enough for us.
\newline
\newline
\noindent
{\bf Remark:\/}
In our examples, both $F_1$ and $F_2$ are
centered at the origin.  One can define
a PET without this property, but some
computer experimentation suggests that the
character of the PET is much different when
$F_1$ and $F_2$ are not centered at the
origin.  In the $2$-dimensional example,
which is the only one we've looked at,
these ``exotic'' examples all have
limit sets which contain open sets.
That is, there are aperiodic orbits which
are dense in open sets.   Indeed, in
the $2$ dimensional example, this seems
to happen for every placement of
$F_1$ and $F_2$ except for the case when
they have a common center (which we might
as well take as the origin.)

\subsection{The Moduli Space}

We keep the notation from the previous section.
Let $M \in GL_n(\R)$ denote any invertible
real $n \times n$ matrix.  Note that
$M$ commutes with both $R_H$ and $R_D$.
In particular $M(H)=H$ and $M(D)=D$.
The system
$(F_1,F_2,L_1,L_2)$ is conjugate to the
system
$(M(F_1),M(F_2),M(L_1),M(L_2))$.
So, up to conjugacy, we might as well
consider systems in which
$\{H_1,...,H_n\}$ is the standard
basis for $\R^n$.

We can specify one of our PETs by giving
a lattice ${\bf s\/} \subset D$.  In
case $n=1$, the lattice ${\bf s\/}$ is $1$ dimensional
and generated by its shortest
vector $s+is$. Identifying
$\C$ with $\R^2$, this shortest vector
becomes $(s,s)$.  For later convenience,
we scale everything by a factor of $2$,
and this gives us the same collection
of objects as discussed in connection with
Figure 1.1.

We can always take
$s \in (0,1)$ because the case $s>1$ can be
reduced to the case $s \in (0,1)$ by
interchanging the roles of $D$ and $H$.
More precisely, the two parameters
$s$ and $s'=1/(2s)$ give rise to 
conjugate systems.  We will formalize
this idea in \S 4, in the Inversion Lemma.
Technically, we could take $s \in (0,\sqrt 2/2)$,
but this further restriction is not convenient to us.

Fixing $n$, let $\cal F$ denote the set
of PETs which arise from our construction
in $\C^n$.  It seems worth pointing out
that $\cal F$ contains a natural
collection of rational points.  These
correspond to taking ${\bf s\/}$ as
a sub-lattice of $(\Q[i])^n$.
We call such systems rational.
In case $n=1$, this just amounts to
taking $s \in \Q$.  Here is an
easy observation.

\begin{lemma}
\label{rat}
For any rational system, all the
orbits are periodic.
\end{lemma}

\startproof
Let $G=(\Z[i])^n$ be the usual
lattice of Gaussian integers.
We can conjugate a rational system
by a dilation so that
$L_1, L_2 \subset G$ and also
the vertices of $F_1$ and $F_2$
belong to $G$.

Letting $(X,f)$ denote the associated PET,
we observe that all the points in an orbit
of $f$ differ from each other by vectors
of $G$.  Moreover, such orbits are bounded.
Hence, they are finite.
\endproof

Consider the $2$-dimensional case, setting
$s=p/q$.  In this case, we can dilate so
that $H_1=q$ and $H_2=p+ip$.  Then
$F_1$ and $F_2$ have volume $(pq)^2$.
In this case, the maximum period of any
point in the system is $(pq)^2$.

\newpage
\section{Stability and Limiting Considerations}

\subsection{Intersection of the Lattices}

For the rest of the paper, we will consider the
case $n=1$ discussed in the previous chapter.
We will work in $\R^2$.  Our main goals in this
chapter are the Stability Lemma and the
Convergence Lemma, stated below.

\begin{lemma}
\label{intersect}
If $s$ is irrational, then $L_1 \cap L_2=0$.
\end{lemma}

\startproof
$L_1$ is the $\Z$-span of $(2,0)$ and $(2s,-2s)$.
and $L_2$ is the $\Z$-span of $(0,2)$ and $(2s,2s)$.
If this lemma is false, then we can find
an equation of the form
\begin{equation}
\label{irr}
(A+Bs,Bs) =(D,C+D).
\end{equation}
for integers $A,B,C,D$.
This forces
\begin{equation}
\label{irr2}
s=\frac{A}{D-B}=\frac{-C}{B+D},
\end{equation}
which is only possible if $s$ is rational.
\endproof

We say that
$(X,f)$ is {\it sharp\/} if, for any vector $V$,
the sets 
\begin{equation}
\{p|\ f(p)=p+V\}
\end{equation}
 are open and convex.
The issue is that such a set, if nonempty,
might be a finite union
of open convex polygons.
The systems $(X_s,f_s)$ are not sharp for
$s=1/n$ and $n=2,3,4...$.

\begin{lemma}
\label{provisional1}
Suppose that $s \in [1/4,1]$. Then
$(X_s,f_s)$ is sharp unless $s=1/n$
for $n=1,2,3,4$.
\end{lemma}

\startproof
we first observe that $(X_s,f_s)$ has the
following property for all $s$.  Suppose
we specify a pair of vectors $(V_1,V_2) \in L_1 \times L_2$.
Then the set of points
$$S(V_1,V_2)=\{p \in F_1|\ p+V_2 \in F_2, p+V_1+V_2 \in F_1\}$$
is always convex:  It is the intersection of $3$ 
parallelograms.  So, if $(X_s,f_s)$ is not sharp,
then we can find vectors $(V_1,V_2) \in L_1 \times L_2$
and $(V_1',V_2') \in L_1 \times L_2$ such that
\begin{itemize}
\item $S(V_1,V_2)$ and $S(V_1',V_2')$ are both nonempty,
\item $V_1+V_2=V_1'+V_2'$ 
\item $(V_1,V_2) \not = (V_1',V_2')$.
\end{itemize}
In this situation, we get a relation like
the one in Equation \ref{irr}.  

Both $L_1$ and $L_2$ have a standard basis.  Let
$\{v_{j1},v_{j2}\}$ be the standard basis for
$L_j$.  For intstance $v_{11}=(2,0)$ and $v_{12}=(2s,-2s)$.
Since $s \geq 1/4$, the vectors $V_j$ and $V_j'$ must
be fairly short.  Specifically, we have
\begin{equation}
V_1=a v_{11} + b v_{12}, \hskip 30 pt
|a| \leq 1 \hskip 20 pt |b| \leq 2.
\end{equation}
Similar bounds hold for the other vectors.
This means that the relation in Equation \ref{irr}
satisies 
$$\max(|A|,|C|) \leq 2, \hskip 30 pt
\min(|B+D|,|B-D|) \leq 4.$$
But then Equation \ref{irr2} forces $s=p/q$ where
$q \leq 4$.  This leaves only the cases
$s=1/n$ for $n=1,2,3,4$.
\endproof

Let $p \in X$ be a periodic point of period $n$.
We define the {\it displacement list\/} of $p$
to be the list of vectors $V_1,...,V_n$ so that
\begin{equation}
f^k(p)=p+\sum_{i=1}^k V_i, \hskip 30 pt k=1,...,n.
\end{equation}

\begin{corollary}
\label{provisional2}
Suppose $s \in [1/4,1]$ and
$s$ does not have the form $1/n$ for $n=1,2,3,4$.
Then the set of periodic points in $X_s$ having the
same period and displacement list is a single
tile of $\Delta_s$.
\end{corollary}.

\startproof
We check the result by hand for $s=1/n$, $n=1,2,3,4$.
For the remaining parameters, the system is sharp,
and the set in question is the intersection of
$n$ convex polygons.
\endproof

\noindent
{\bf Remark:\/}
In Corollary \ref{graph}, we promote Lemma
\ref{provisional1} and Corollary \ref{provisional2}
to results about
all $s \in (0,1)$.

\subsection{The Arithmetic Graph}
\label{ag}

Suppose that $p_0 \in X$ is some point
on which the orbit of $f$ is well defined.
Call this orbit $\{p_i\}$.
There are unique vectors $V_i \in L_1$
and $W_i \in L_2$ such that
$$f'(p_i)=p_i+V_i, \hskip 40 pt
f'(p_i+V_i)=p_i+V_i+W_i=p_{i+1}.$$
We call the sequence
$\{(V_i,W_i)\}$ the
{\it symbolic encoding\/} of the orbit.
We define the {\it arithmetic graph\/} to be the
polygon whose $i$th vertex
is $V_0+...+V_i$. We define the
{\it conjugate arithmetic graph\/} to be
the polygon whose $i$th vertex is
$W_0+...+W_i$.
\newline
\newline
{\bf Remark:\/}
The displacement list is determined from
the arithmetic graph, but the reverse is not
{\it a priori\/} true.  One way
to interpret Corollary \ref{provisional2} (and
the more general Corollary \ref{graph}) is
that, unless $s=1/n$, we can determine the
graph from the displacement list.
\newline

We hope to explain the structure of the
arithmetic graphs associated to systems
such as these in a later paper.  In this paper,
we are mainly interested in using the graph
to detect when an orbit is stable under
perturbation of the parameter.

\begin{lemma}[Stability]
If $s$ is irrational and $p$ is a periodic point
of $f$, then both the arithmetic graph of $p$ is
a closed polygon.
\end{lemma}

When the orbit is periodic, of period $n$,
we have
\begin{equation}
\sum_{i=1}^n(V_i+W_i)=0.
\end{equation}
We can re-write this as
$$
\sum_{i=1}^n V_i= - \sum_{i=1}^n W_i.
$$
But then the common sum belongs to
$L_1 \cap L_2$.  When $s$ is
irrational, Lemma \ref{intersect} then
tells us that
\begin{equation}
\label{stable}
\sum_{i=1}^n V_i=0, \hskip 30 pt
\sum_{i=1}^n W_i=0.
\end{equation}
These two equations are equivalent to the lemma.
\endproof

\subsection{Convergence Properties}

We say that a sequence of (solid) polygons $P_n$
converges to a (solid) polygon $P_{\infty}$ if
the sequence converges in the {\it Hausdorff metric\/}
on the set of compact sets. Concretely, for every
$\epsilon>0$ there should be an $N$ such that
$n>N$ implies that every point of $P_n$ is
within $\epsilon$ of $P_{\infty}$ and
{\it vice versa\/}.
In the next lemma, we write
$\Delta_{\infty}=\Delta_{s_{\infty}}$ and
$\Delta_{n}=\Delta_{s_{n}}$ for ease of
notation.

\begin{lemma}[Approximation]
\label{converge}
Let $s_{\infty} \in (0,1)$ be irrational, and
let $\{s_n\}$ be a sequence of rationals
converging to $s_{\infty}$.
Let $P_{\infty}$ be a tile of $\Delta_{\infty}$.
Then for all large $n$ there is a tile
$P_n$ of $\Delta_n$ such that $\{P_n\}$ converges
to $P_{\infty}$.
\end{lemma}

\startproof
For any object $A$ that depends on a
parameter, we let $A(n)$ denote the
object corresponding to the parameter $s_n$.

Let $p \in P_{\infty}$ be a periodic
point of period $N$.
Let $\{(V_i(\infty),W_i(\infty)\}$
be the symbolic encoding of the orbit.
The lattice $L_j(n)$ converges to
the lattice $L_j(\infty)$.  So, 
we can uniquely choose vectors
$V_i(n) \in L_1(n)$ and
$W_i(n) \in L_2(n)$ which converge
respectively to $V_i(\infty)$ and
$W_i(\infty)$.  

By the Stability Lemma,
$\sum_{i=1}^N V_i(\infty)=0.$  Hence
$\sum_{i=1}^N V_i(n) \to 0$ as
$n \to \infty$.
But, independent of $n$, there is some
$\epsilon$ so that the $\epsilon$ balls centered
on lattice points of $L_1(n)$ are disjoint.
This shows that
$\sum_{i=1}^N V_i(n)=0$
for large $n$.  The same argument works with
$W$ in place of $V$.  For $n$ large, the
symbolic encoding of $p$ starts out
$$V_1(n),W_1(n),...,V_N(n),W_N(n).$$
But then $p$ is a periodic point of
period $N$ for $n$ large, and
the above list is the whole symbolic
encoding.  In short, $p$ is a periodic
point of the same period relative
to $f_n$ as it is relative to $f_{\infty}$.
Let $P_n$ denote the periodic tile containing $p$.

The size of $n$ required for this continuity
argument depends only on the distance
from $p$ to $\partial P_{\infty}$.
For this reason, $P_{\infty}$ is
contained in the $\epsilon$ neighborhood of
$P_n$ for $n$ sufficiently large.
On the other hand,
if $q$ lies just outside $P_{\infty}$, the
symbolic encoding of $q$ relative to
$f_{\infty}$ differs from the
symbolic encoding of $p$ relative to
$f_{\infty}$.  But then, by continuity,
the same goes for large $n$ in place
of $\infty$.  This shows that
$P_n$ is contained in an arbitrarily
small neighborhood of $P_{\infty}$,
once $n$ is sufficiently large.

Putting everything together, we
see that $\{P_n\}$ converges to
$P_{\infty}$, as desired.
\endproof

\newpage

\section{Symmetry}

\subsection{Rotational Symmetry}

In this chapter we discuss the symmetry of the system $(X_s,f_s)$.
We begin with an obvious result.
Define
\begin{equation}
\iota(x,y)=(-x,-y)
\end{equation}
Note that $\iota(X)=X$ for every parameter.

\begin{lemma}[Rotation]
$\iota$ and $f_s$ commute for all $s \in (0,1)$.
\end{lemma}

\startproof
$\iota$
preserves $F_1$, $F_2$, $L_1$, and $L_2$.  For this reason $\iota$
commutes with $f$. 
\endproof

 The subsets of $\Delta_s$ and $\Lambda_s$ lying to
the right of the central tiles are reflected images of the subsets of
$\Delta_s$ and $\Lambda_s$ lying to the left of the central tiles. For
this reason, we will usually consider the picture just on the left
hand side. 

\subsection{Inversion Symmetry}

As we have already remarked,
we usually take the parameter $s$ to lie in $(0,1)$
but we can define all the objects for any $s \in (0,\infty)$.

\begin{lemma}[Inversion]
\label{invert}
$(X_t,f_t)$ and $(X_s,f_s)$ are conjugate if $t=1/2s$.
\end{lemma}

\startproof
An easy calculation shows that there
is a similarity $\phi: R_t \to R_s$ which maps the
horizontal (respectively diagonal) side of $R_t$ to
the diagonal (respectively horizontal) side of $R_s$.

We denote $F_1$, at the parameter $s$, by
$F_1^s$.  We make similar notations for the other
parameters.  We already know that
$\phi(F_1^t)=F_1^s$.
Let $J$ be rotation by $\pi/2$ clockwise.  The map
$\phi$ conjugates $J$ to $J^{-1}$.  For this
reason, we have
$\phi(F_j^t)=F_j^s$ and 
$\phi(L_j^t)=L_j^s$ for $j=1,2$.
This does it.
\endproof

Lemma \ref{invert} is really the source of the renormalization
map $R$.  However, Lemma \ref{invert} does not directly
apply to our situation.  Let's consider the situation 
in some detail.  Let $\rho(s)=1/2s$.  The map
$\rho$ is an involution of the parameter interval
$[1/2,1]$.  However, the Main Theorem requires us to
use the map $R(s)=1-s$ on this interval.  When 
$s \in 1/2$ we have $R(s)=t_1$ and
$\rho(s)=t_2$, where $t_1=t_2-k$, where $k$ is the
integer such that $t_2-k \in (0,1]$.  However, it
is not yet clear that the two systems at 
$t_1$ and $t_2$ are related.

\subsection{Insertion Symmetry}
\label{insertion}

When $s \in (1/2,1)$, the intersection
$F_1 \cap F_2$ is an octagon, which we
call the {\it central tile\/}.
When $s \leq 1/2$ or $s \geq 1$, 
the intersection 
$F_1 \cap F_2$ is a square. This square
generates a grid in the plane, and
finitely many squares in this grid lie in
$X=F_1$.  We call these squares the
{\it central tiles\/}.  See Figure 4.1.

\begin{center}
\resizebox{!}{2.3in}{\includegraphics{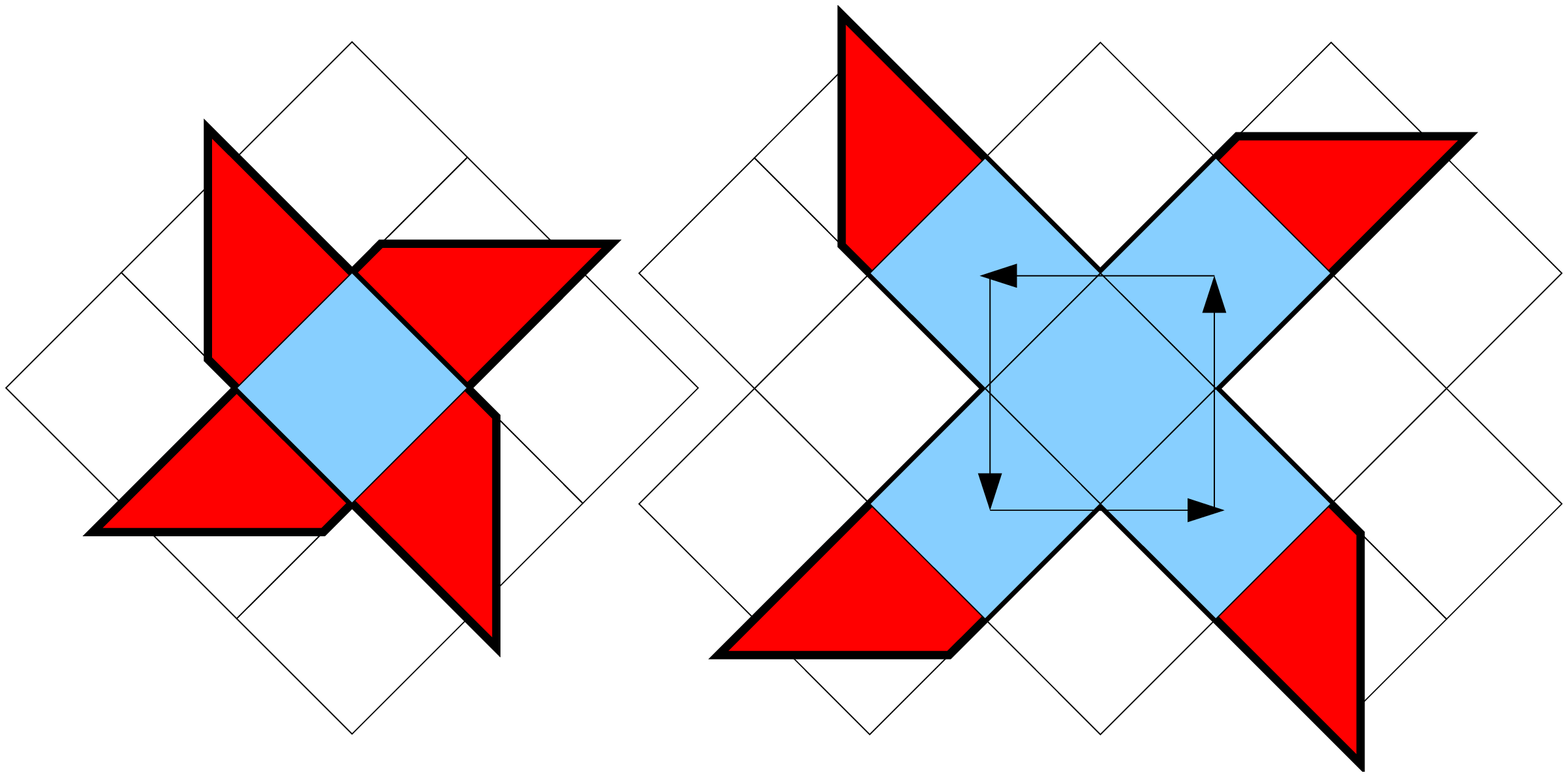}}
\newline
{\bf Figure 4.1:\/} 
The central tiles (blue) for $s=5/4$ and $t=9/4$.
\end{center}

Let $X^0$ denote the portion of $X$ which lies
to the left of the central tiles.  The set
$\iota(X^0)$ is the portion of $X$ which lies to
the right of the central tiles.  The union
$X^0 \cup \iota(X)$ is an $f$-invariant set.

\begin{lemma}[Insertion]
\label{insert}
Suppose $s \geq 1$ and $t=s+1$, or
suppose $s \leq 1/2$ and
$t=s/(2s+1)$.
The restriction of $f_s$ to
$X_s^0 \cup \iota(X_s^0)$ is
conjugate to the restriction of $f_t$ to
$X_t^0 \cup \iota(X_t^0)$.  The conjugacy
is a piecewise similarity.
\end{lemma}

\startproof
The case when $s<1/2$ is equivalent to the
case $s>1$ by the Inversion Lemma.  So,
we will take $s>1$ and $t=s+1$.

We consider how $X_s$ and $X_t$ sit
relative to the grid of diamonds mentioned
above.
When $s,t>1$, the squares in the grid are
diagonals of length $2$.
Significantly, the diagonals of the
diamonds are parallel to the vectors
$(\pm 2,0) \in L_1$ and $(0,\pm 2) \in L_2$.
This is true independent of the parameters.
There are two more diamonds contained
in $X_t$ than there are in $X_s$.
On the central tile, the map $f$
has the obvious
action shown in Figure 4.1.

The sets $X_s^0$ and $X_s^t$ have the
same relative position relative to
the diamond grid, and there is an
obvious translation carrying the one
set to the other.  This translation
extends to give piecewise translation
from the complement of the diamonds
in $X_s$ to the complement of the
diamonds in $X_t$.  Call points related
by this piecewise translation {\it partners\/}

Let $p_s \in X_s^0$ and $p_t \in X_s^t$ be partners.
Let $\lambda_s$ and $\lambda_t$
respectively be the vectors in
$(L_2)_s$ and $(L_2)_t$ such that
$p_s+\lambda_s \in (F_2)_s$ and
$p_t+\lambda_t \in (F_2)_t$.
We have either
$\lambda_s=\lambda_t+(2,0)$ or
$\lambda_s=\lambda_t+(0,2)$,
depending on whether or not
$p_s+\lambda_s$ and $p_t + \lambda_t$ lie
in the top or bottom of $(F_2)_s$ and
$(F_2)_t$ respectively.  The answer
(top/bottom) is the same for $s$ as it is for $t$.
In short, the two new points
$p_s+\lambda_s$ and $p_t+\lambda_t$ are
again partner points. Repeating this construction
again, we see that $f_s(p_s)$ and $f_t(p_t)$
are partner points. This is what we wanted to prove.
\endproof

\noindent
{\bf Remark:\/}
Informally, what Insertion Lemma says
is that, for $s \geq 1$,
the tiling $\Delta_{s+1}$ is
obtained from the tiling
$\Delta_s$ by inserting two new
large diamonds.

\begin{lemma}
\label{basic}
$\Delta_s$ consists entirely of squares and
right-angled isosceles triangles when
$s=1,2,3...$ and when 
$s=1/2, 1/4, 1/6,...$.
\end{lemma}

\startproof
We check this for $n=1$ just by a direct calculation.
See Figure 4.2.

\begin{center}
\resizebox{!}{1in}{\includegraphics{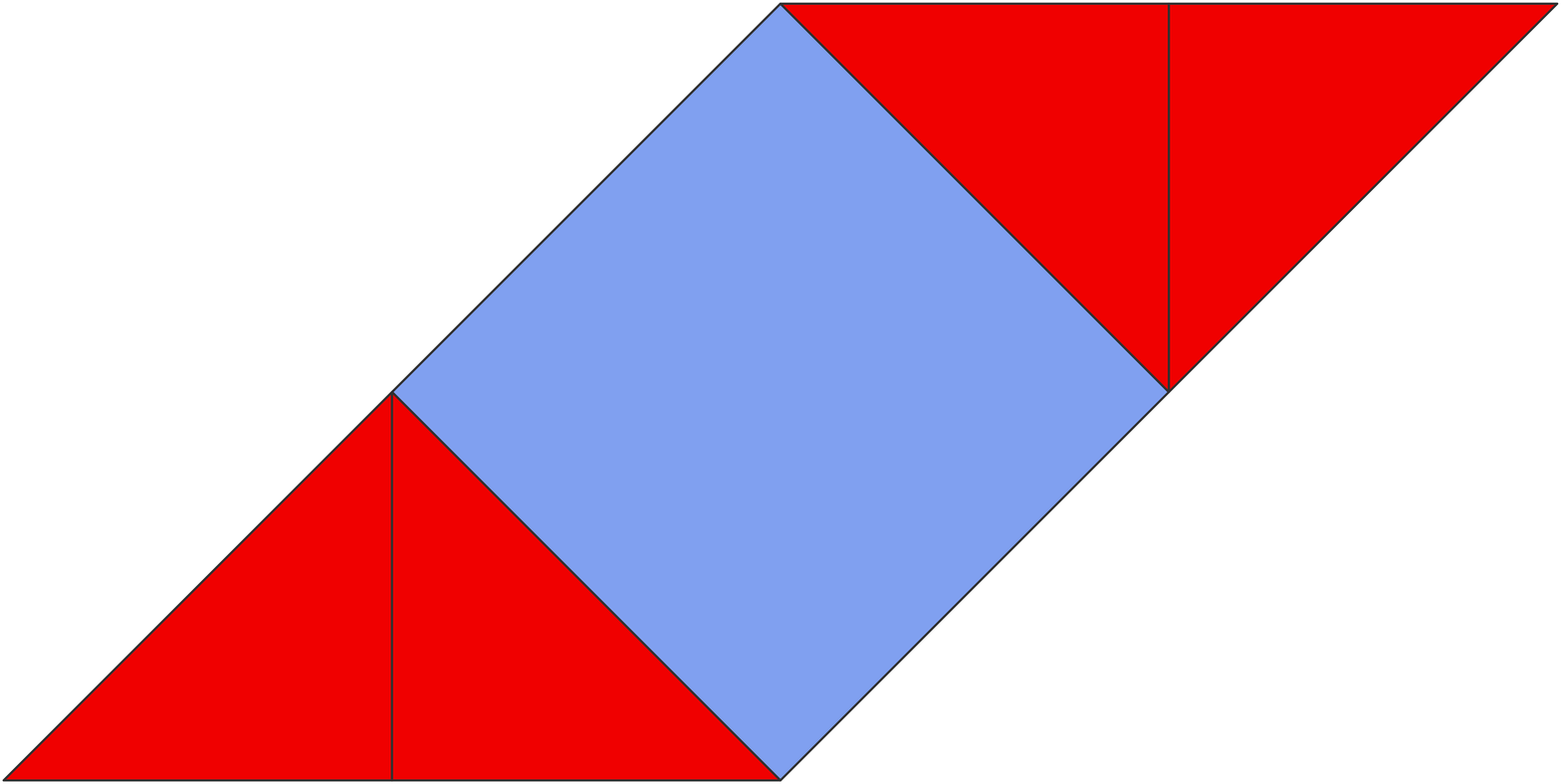}}
\newline
{\bf Figure 4.2:\/} 
The tiling $\Delta_s$ for $s=1$.
\end{center}

The cases $n=2,3,4...$ now follow from the Insertion
Lemma. The cases $n=1/2,1/4,1/6....$ follow from
the cases $n=1,2,3...$ and the Inversion Lemma.
\endproof

Combining the Insertion Lemma with Lemmas \ref{provisional1}
and  Corollary \ref{provisional2}, we have

\begin{corollary}
\label{graph} 
Let $s \in (0,1)$.  Suppose $s$ does not
have the form $1/n$
for $n=1,2,3,...$.  
Then $(X_s,f_s)$ is clean, and
the set of periodic points in $X_s$ having the
same period and displacement list is a single
tile of $\Delta_s$.
\end{corollary}

\subsection{Bilateral Symmetry}

In this section we always take $s \in (0,1)$.
We say that a line $L$ is a {\it line of symmetry\/}
for $\Delta_s$ if 
\begin{equation}
\Delta_s \cap \bigg(X_s \cap \rho(X_s)\bigg)
\end{equation}
is invariant under the reflection $\rho$ in $L$.
Note that $X_s$ itself need not be invariant
under $L$.

Figure 4.3 shows $3$ lines. $H$ is the line $y=0$ and
$V$ is the line $x=-1$ and
$D_s$ is the line of slope $-1$ through
the bottom vertex of the leftmost central square
of $\Delta_s$.

\begin{center}
\resizebox{!}{1in}{\includegraphics{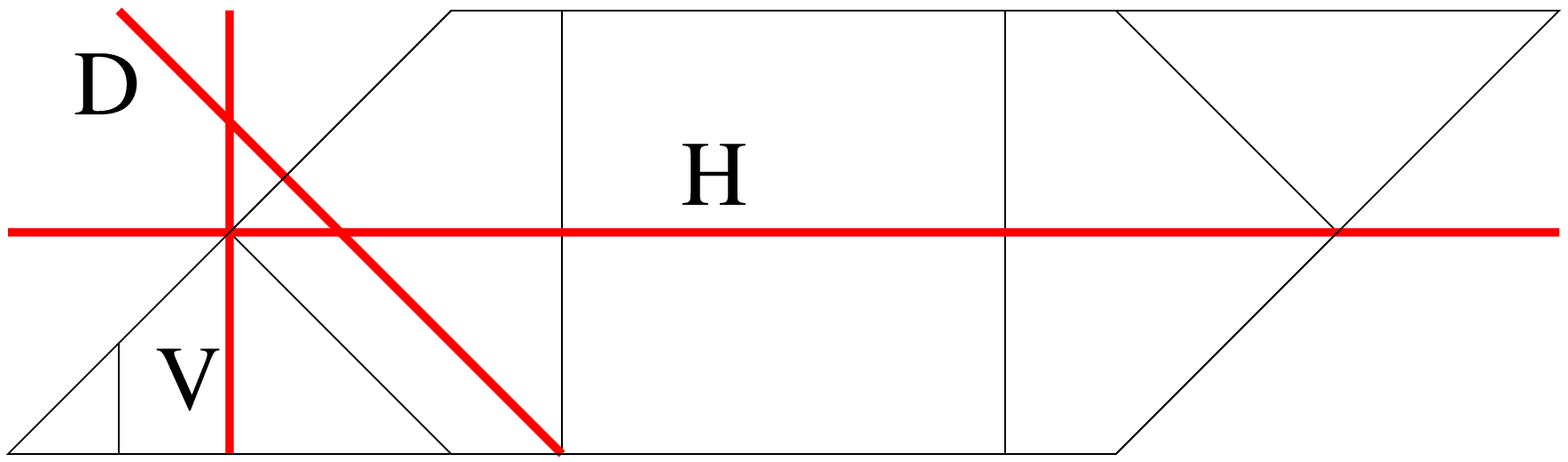}}
\newline
{\bf Figure 4.3\/} $H$ and $V$ and $D_s$ for the parameter $s=2/5$.
\end{center}

In this section we will prove (modulo $2$ finite calculations)
the following result.

\begin{lemma}[Bilateral]
For all $s \in (0,1)$, the lines
$H$, $V$, and $D_s$ are lines of
symmetry of $\Delta_s$.
\end{lemma}

We will prove this result through a
series of smaller lemmas. 

Define
\begin{equation}
A=X \cap \rho_H(X), \hskip 30 pt
B=X \cap \rho_V(X).
\end{equation}
Here $\rho_H$ is the reflection in $H$ and
$\rho_V$ is the reflection in $V$.
$A$ is the green hexagon shown in Figure 4.4.
The complement $X-A$ consists of two
triangles, $B$ and $\iota(B)$ as shown in Figure 4.4.

\begin{center}
\resizebox{!}{1in}{\includegraphics{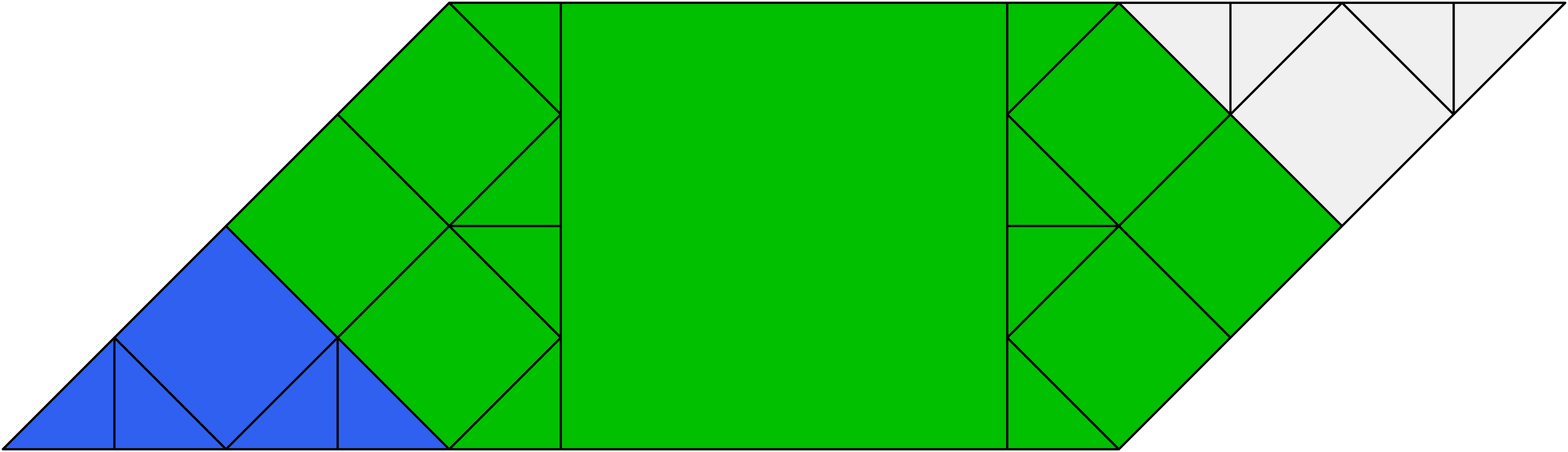}}
\newline
{\bf Figure 4.4\/} $A_s$ (green) and
$B_s$ (blue) and $\iota(B_s)$ (white) for $s=2/5$.
\end{center}

Each of the pieces $A$, $B$ and $\iota(B)$ has a
vertical line of bilateral symmetry.  The reflections
across these vertical lines gives rise to a piecewise
isometry of $X$.  We call this map $\mu$.
If $p \in A$ we define $\mu(p)$ to be reflection
in the vertical line of symmetry for $A$, etc.

The Insertion Lemma allows us to turn many
infinite-appearing calculations into finite
calculations.  The problem with trying to
compute something for every parameter 
in $(0,1)$ is that the number of domains
of continuity for the map $f_s$ tends to
$\infty$ as $s \to 0$.  However, if we
restrict our attention to $s \in [1/4,1)$,
then there is a uniform bound on the number
of regions of continuity.  In this range,
we can establish identities using a finite
calculation.  The picture for any parameter
$s<1/4$ is the same as some picture for
$s'>1/4$, up to the insertion of finitely
many central tiles.
In \S \ref{calc} we prove the following result by
direct calculation.

\begin{lemma}[Calculation 1]
If $s \in [1/4,1]$ then
$\mu_s \circ f_s \circ \mu_s=f_s^{-1}$
wherever both maps are defined.
\end{lemma}

Combining this calculation with the Insertion Lemma, we have
\begin{corollary}
\label{symm1}
Suppose $s \in (0,1)$ then
$\mu_s \circ f_s \circ \mu_s=f_s^{-1}$
wherever both maps are defined.
\end{corollary}

\noindent
{\bf Proof of Statements 1 and 2:\/}
Now we prove Statements 1 and 2 of the Bilateral Lemma.
Using the rotational
symmetry, it suffices to prove that $\Delta \cap A$
and $\Delta \cap B$ are invariant under the
action of $\mu$.  We will consider the sitution in the Hexagon
$A$.  The situation in the other regions has a similar
treatment.

Let $\tau$ be a tile of $\Delta$ that is
contained in $A$.
 All iterates of $f$ are
defined on the interior of $\tau$.  Let
$n$ be the order of $f$ on the interior of
$\tau$.  The first $n$ iterates of
$\mu f \mu$ are defined on an open
set $\tau'=\tau-L$. Here $L$ is a finite union
of line segments.  For each $p \in \tau'$,
the period of $\mu f \mu$ on $p$ is $n$.
Since $\mu$ is everywhere defined in the
interior of $A$, we see that $f^{-1}$ is defined
on all points of $\mu(\tau')$.  But all
the points in $\mu(\tau')$ have the same
displacement list.  Hence $\mu(\tau')$ is
convex. This is only possible if
$\tau'=\tau$.
\endproof

Define 
\begin{equation}
P=X \cap \rho_D(X), \hskip 30pt
Q=X^0-P.
\end{equation}
When $s<1/2$, the set $P_s$ is a pentagon
and $Q_s$ is an isosceles triangle.
When $s>1/2$, the set $P_s$ is a triangle
and $Q_s$ is empty.
Figures 4.5 shows the case $s<1/2$.

\begin{center}
\resizebox{!}{1.7 in}{\includegraphics{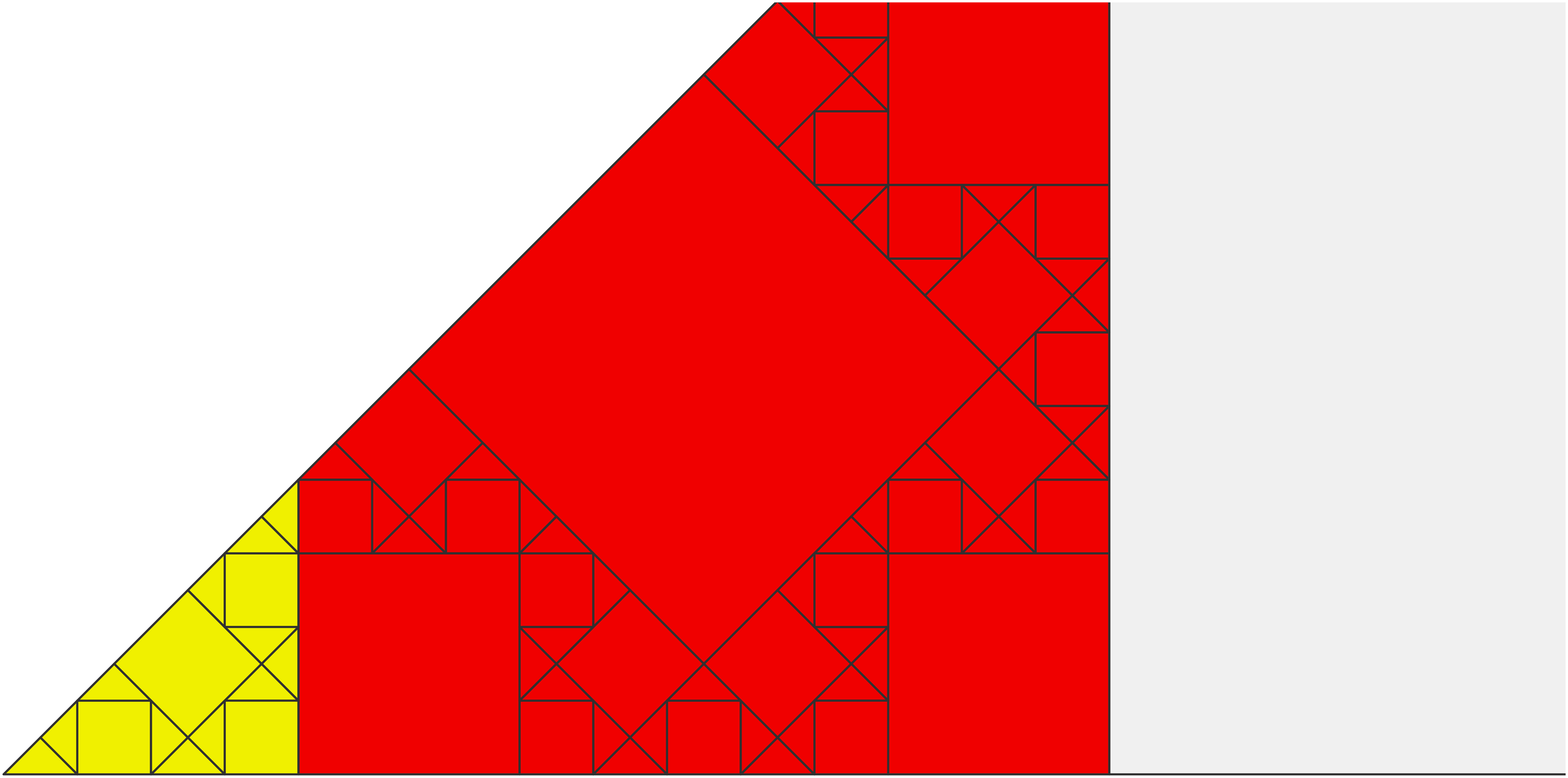}}
\newline
{\bf Figure 4.5:\/}
$P_s$ (red) and $Q_s$ (yellow) for $s=11/30$.
\end{center}

$X_s$ is partitioned into square central
tiles and the additional
tiles $P_s$, $Q_s$, $\iota(P_s)$ and
$\iota(Q_s)$. (When $s>1/2$, the tiles
$Q_s$ and $\iota(Q_s)$ do not exist.)
Each tile in this partition
has reflection symmetry, in a line of slope
$-1$.  Let $\nu_s: X_s \to X_s$ be the
piecewise isometry which does this reflection
on each piece.

In \S \ref{calc} we prove the following result by
direct calculation.

\begin{lemma}[Calculation 2]
If $s \in [1/4,1]$, then
$\nu_s \circ f_s \circ \nu_s=f_s^{-1}$
wherever both maps are defined.
\end{lemma}

Combining this calculation with the Insertion Lemma, we have

\begin{corollary}
\label{symm2}
If $s \in (0,1)$ then
$\nu_s \circ f_s \circ \nu_s=f_s^{-1}$
wherever both maps are defined.
\end{corollary}

\noindent
{\bf Proof of Statement 3:\/}
Statement 3 of the Bilateral Lemma is deduced
from Corollary \ref{symm2} in the same
way that Statements 1 and 2 are deduced
from Corollary \ref{symm1}.
\endproof

\newpage

\section{Proof of the Main Theorem}

\subsection{Discussion and Overview}

The Inversion Lemma and the Insertion Lemma go part of
the way towards proving the Main Theorem.  These two
results say that the systems $(X_s,f_s)$ and
$(X_t,f_t)$ are related, in the appropriate sense,
for pairs $(s,1/2s)$ and, assuming $s>1$, for
pairs $(s,s-1)$.
These results are not strong enough to
establish the Main Theorem. For instance,
when $s=2/5$ we have $R(s)=1/4$.  We can say
that the paramters $2/5$ and $5/4$ are related
by Condition 1 above.  However, Condition 2
does not apply to $(s,t)=(1/4,5/4)$ 
because $s<1$.  Similarly, if $s=3/4$ we have
$R(s)=1/4$.  Here, neither condition applies.

The reader might wonder why we care about $R$ in the
first place.  Perhaps we can prove all the corollaries
to the Main Theorem
just with the limited symmetries we have already
established.  The virtue of $R$ is that, for
every rational parameter $p/q$, one of the two
iterates $R(p/q)$ or $R^2(p/q)$ has denominator
smaller than $q$.  Thus, the map $R$ gives us a
an inductive mechanism for understanding our system
at all rational values.  Once we have a good understanding
of what happens at rational values, we can take limits.
We cannot do this much with just the two conditions
listed above.

Referring to Theorem \ref{three},
the existence of $\Gamma$ sheds light
on what we have said above.
The way we prove the Main Theorem,
roughly speaking, is to verify certain
facts on the generators of $\Gamma$,
by computation or symmetry, and then use the group
structure to extract global statements
about the renormalization map $R$.
What we are saying, in a sense, is that we haven't
checked {\it enough\/} of $\Gamma$ yet.
What is missing is a statement about what happens for
pairs $(s,s-1)$ with $s \in (1,2)$ and for pairs
$(s,1-s)$ with $s \in (1/2,1)$.

Here are the two results we prove in this
chapter.  The first is
equivalent to the half of the Main
Theorem corresponding to $s \in (0,1/2)$.

\begin{lemma}
\label{main1}
Suppose $s \in (1,2)$ and $t=s-1$.
Let $\phi_s: Y_t \to X_s$ be the map which is
a translation on each half of $Y_t$ and
maps the acute vertices of $Y_t$ to the
acute vertices of $X_s$.  Let $Z_s=\phi_s(Y_t)$.
Then $\phi_s$ conjugates $f_t|Y_t$
to $f_s|Z_s$, and $Z_s$ is a clean set.
Either half of $\phi_s$ extends to the trivial tile of
$\Delta_t$ and maps it to tiles
$\tau_1$ and $\tau_2$.
The only nontrivial
$f_s$-orbits which miss $Z_s$ are
contained in $\tau_1 \cup \tau_2$ and have
period $2$. 
\end{lemma}

Our other result is just a restatement of the half
of the Main Theorem corresponding to $s \in (1/2,1)$.

\begin{lemma}
\label{main2}
Suppose $s \in (1/2,1)$ and $t=1-s$.
Let $\phi_s: Y_t \to X_s$ be the map which is
a translation on each half of $Y_t$ and
maps the acute vertices of $Y_t$ to the
acute vertices of $X_s$.  Let $Z_s=\phi_s(Y_t)$.
Then $\phi_s$ conjugates $f_t|Y_t$
to $f_s^{-1}|Z_s$, and $Z_s$ is a
clean set.  All nontrivial $f_s$-orbits
intersect $Z_s$. 
\end{lemma}

To be sure, let's deduce the Main Theorem
from these results.
\newline
\newline
\noindent
{\bf Proof of the Main theorem:\/}
Lemma \ref{main2} is just a restatement of
the Main Theorem for $s \in (1/2,1)$.
Suppose that $s<1/2$.
By the Insertion Lemma, it suffices
to consider the case when $s \in (1/4,1/2)$.
By the Inversion Lemma,
the system $(X_s,f_s)$ is conjugate
to the system $(X_t,f_t)$, where
$t=1/2s$.   Here $t \in (1,2)$.
But now Lemma \ref{main1} applies to
the pair $(t,t-1)$ and $t-1=R(s)$.
When we combine the conjugacy given by the
Inversion Lemma with the one given by Lemma
\ref{main1}, we get the statement
of the Main Theorem.
\endproof

We would like to have conceptual proofs of
Lemmas \ref{main1} and \ref{main2}, but we do not.
Instead, we will give computational proofs.  The
difficulty in giving a computational proof is that
it seems to involve an infinite amount of 
calculation.  Consider, for instance, what
happens in Lemma \ref{main1} as $s \to 1$.
In this case, the area of $Y_t$ tends to $0$.
But then, the proportion of $X_s$ taken up
by $Z_s$ tends to $0$.  But then the amount
of time it takes for some orbits to return to
$X_s$ probably tends (and, in fact, does tend) to
$\infty$.  This makes a direct computer
verification difficult.  A similar
problem happens for Lemma \ref{main2} as
$t \to 1$.

In \S \ref{calc} we will prove the following
results by a direct and finite calculation.

\begin{lemma}[Calculation 3]
\label{aux1}
Lemma \ref{main1} holds for all $t \in [5/4,2]$.
\end{lemma}

\begin{lemma}[Calculation 4]
\label{aux2}
Lemma \ref{main2} holds for all $t \in [1/2,3/4]$.
\end{lemma}

The trick is to relate the system on the intervals
$[1/2,3/4]$ and $[5/4,2]$ to the larger intervals
$[1/2,1)$ and $(1,2]$.  We will do this by
establishing some auxilliary symmetry results.
One can view these auxilliary results as statements
about some of the other elements in the group
$\Gamma$.

\subsection{First Modular Symmetry}

Define the maps
\begin{equation}
\label{modular}
T(s)=\frac{s-2}{2s-3}, \hskip 30pt
\omega_s(x,y)=(3-2s)(x,y) \pm (2-2s,0).
\end{equation}
We take $s \in (1,4/3]$, so that
$u=T(s) \in (1,2]$.  The domain
of $\omega_s$ is the set $Y_u$
from the Main Theorem.  
The $(+)$ option for $\omega_s$ 
is taken when $x<0$ and
the $(-)$ option is taken then $x>0$.
We set $W_s=\omega_s(Y_u)$.  
A picture says a thousand words.
In the picture $W_s^0$ is the left half
of $W_s$ and $Y_u^0$ is the left half of $Y_u$.

\begin{center}
\resizebox{!}{1.3in}{\includegraphics{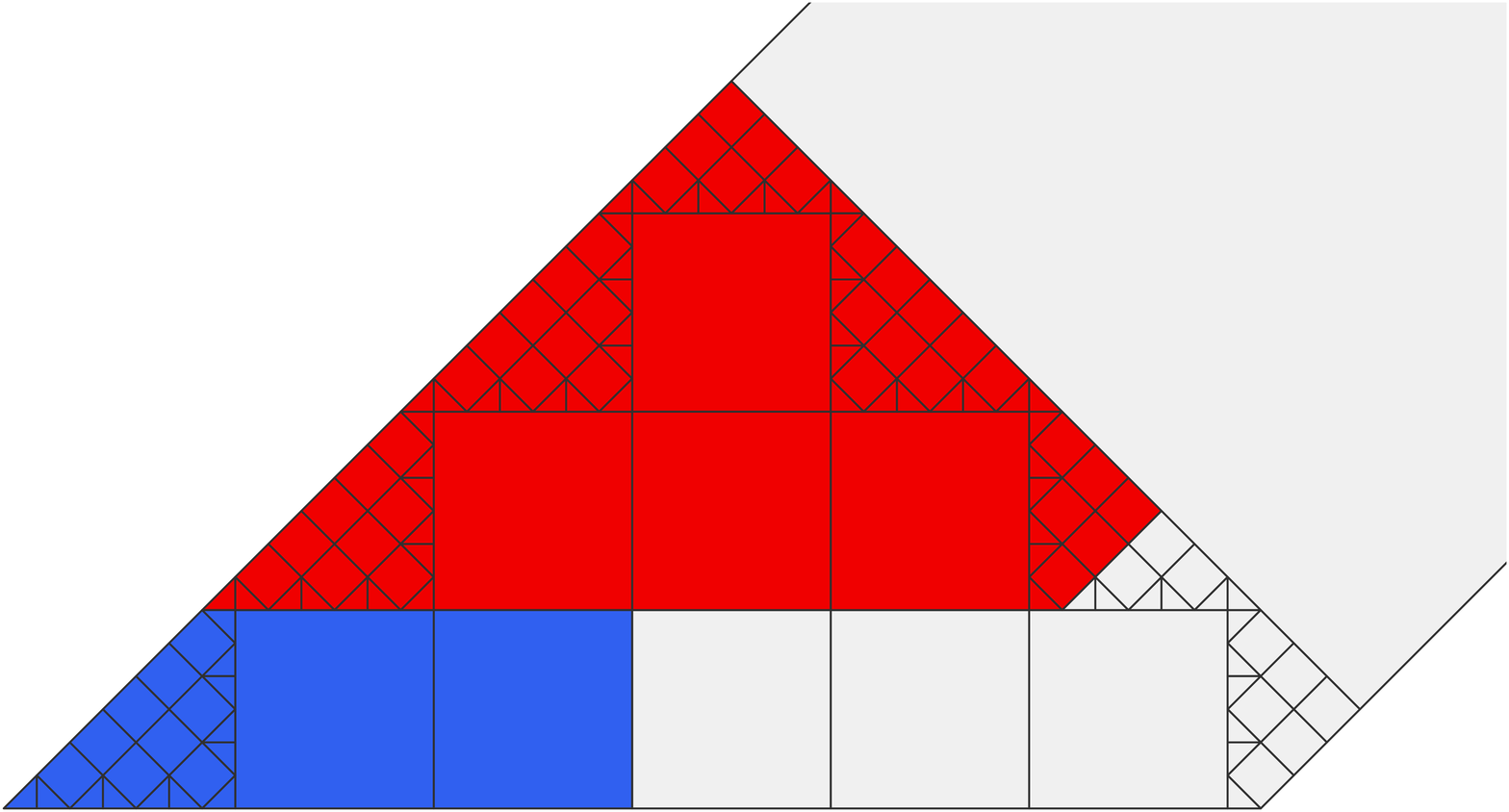}}
\newline
{\bf Figure 5.1:\/} 
$\Delta_s \cap W_s^0$ (red) for $s=22/19$.
\end{center}

\begin{center}
\resizebox{!}{1.3in}{\includegraphics{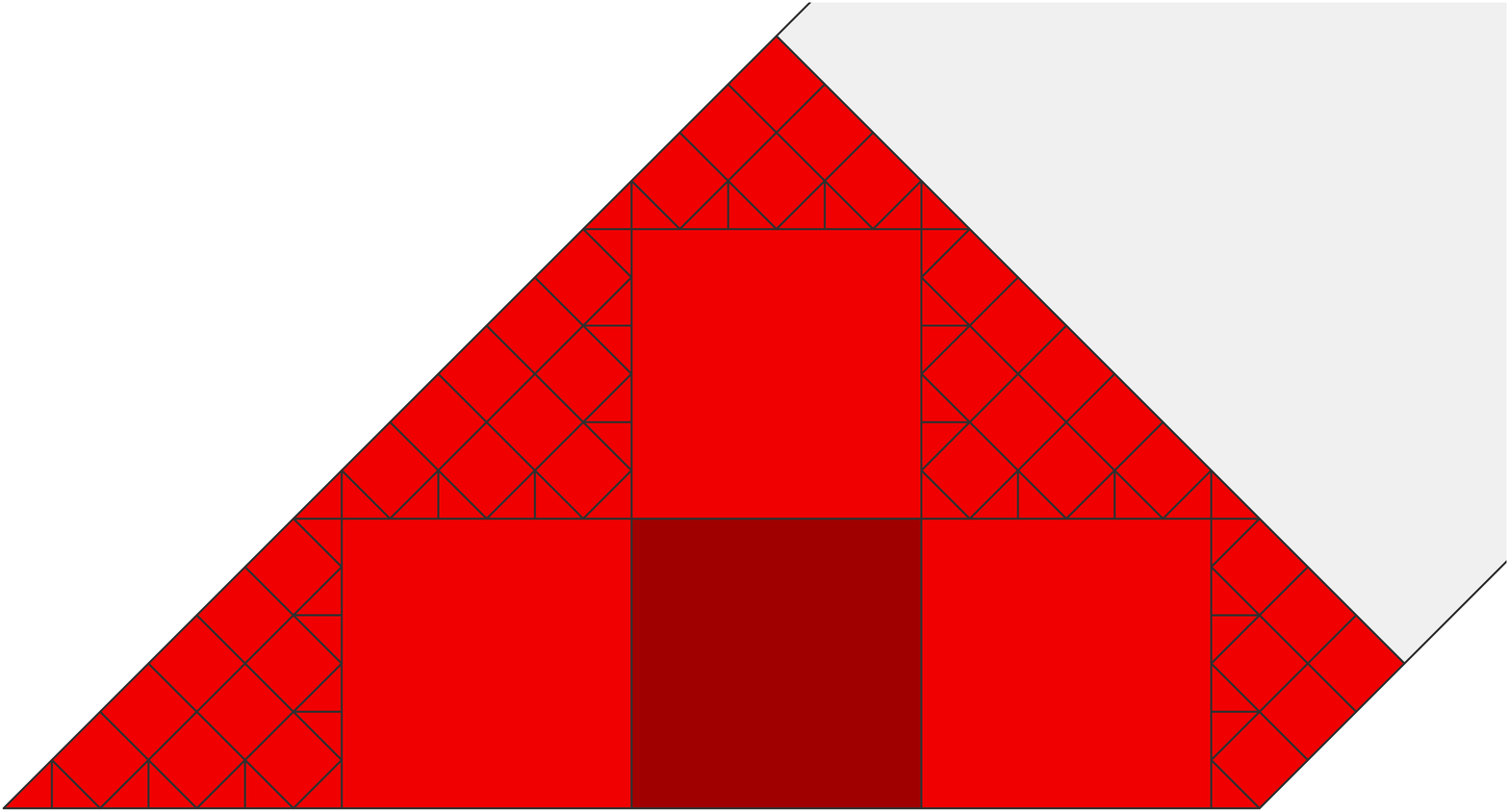}}
\newline
{\bf Figure 5.2:\/} 
 $\Delta_{u} \cap Y_{u}^0$ (red) for $u=T(22/19)=16/13$.
\end{center}

\begin{lemma}[Calculation 5]
\label{calc5}
Let $s \in (1,4/3]$ and $u=T(s)$.
Then $\omega_s$ conjugates $f_{u}|Y_{u}$ to
$f_s|W_s$. Moreover, $Z_s$ is a clean set.
\end{lemma}

\noindent
{\bf Remarks:\/} \newline
(i)
As a corollary, we see that
$\omega_s$ maps $\Delta_{u} \cap Y_{u}^0$ to
$\Delta_s \cap W_s^0$.  This explains why the
tilings in the red regions in 
Figures 5.1 and 5.2 look the same.
\newline
(ii) We only need Calculation 5 for
$s \in (1,5/4]$, but it is more convenient to
make the calculation on the larger interval.
\newline

The blue set in Figure 5.1 is isometric the left
half of $Y_{s-1}$.  We denote this set by
$Z_s^0$ and we set $Z_s=Z_s^0 \cup \iota(Z_s)$.
Here $\iota$ is reflection in the origin.
We let $\tau$ denote the square whose left side
coincides with the right side of $Z_{u}^0$. 
The square $\tau_{u}$ is a darker red
than the others in Figure 5.1. We define
$\tau_s$ just as we defined $\tau_u$,
with $s$ in place of $u$.
Let $\delta_s$ be the vector which spans the diagonal
of $\tau_s$, pointing from the bottom left vertex
to the top right vertex.

\begin{lemma}[Calculation 6]
\label{calc6}
Let $s \in (1,5/4]$, so that $u=T(s) \in (1,3/2]$.
Then
\begin{enumerate}
\item $\tau_u$ is a tile of $\Delta_u$, having period $2$.
\item $f_s^{-1}(p)=p+\delta_s$ for all $p \in Z_s^0$.
\item $f_s^{-1}(X_s - Z_s - W_s) \subset Z_s \cup \tau_s \cup \iota(\tau_s)$.
\end{enumerate}
\end{lemma}

We will establish these results in \S \ref{calc}.
Calculation 6 really just amounts to inspecting
the partitions for $f_{u}$ and
$f_s$.

\subsection{Proof of Lemma \ref{main1}}

Let $T$ be the map from Equation \ref{modular}.

\begin{lemma}
\label{orb}
Let $s \in (1,5/4)$ be any point.  Then
there is some positive $k$ such that
$T^k(s) \in (5/4,3/2)$.
\end{lemma}

\startproof
$T$ is a parabolic linear fractional
transformation fixing $1$ and having
the property that $T(5/4)=3/2$.  So,
the iterates $T^j(s), j=1,2,3...$
are increasing, but then cannot avoid
the interval $(5/4,3/2)$.
\endproof

Lemma \ref{main1} follows immediately from
Lemma \ref{orb} and from the following result.

\begin{lemma}
\label{induct}
If Lemma \ref{main1}
is true for some $u \in (1,3/2)$, then Lemma \ref{main1}
is also true for $s=T^{-1}(u)$. 
\end{lemma}

The rest of this section is devoted to proving
Lemma \ref{induct}.
Let $s'=s-1$ and
$u'=u-1$.   A calculation shows that
\begin{equation}
s'=\frac{u'}{2u'+1}.
\end{equation}
In other words, $s'$ and $u'$ are related
exactly as in the Insertion Lemma.
So, the dynamics relative to
$s'$ is the same as the dynamics relative to
$u'$, except that two more central
squares are inserted for $s'$.  These
central squares have period $2$.
We need to establish the same relation
between $f_{u}|Z_{u}$
and $f_s|Z_s$.

Define
\begin{equation}
Z_{u}^*=Z_{u} \cup \tau_{u} \cup
\iota(\tau_{u}).
\end{equation}
We are augmenting $Z_{u}$ by inserting
two period-$2$ squares at the two ends
of $Z_u$. Statement 1 of Calculation 6 guarantees
that these extra squares really are tiles of
$\Delta_{u}$.

\begin{lemma}
\label{proof1A}
There is a piecewise homothety $h_s$ which carries
$Z_{u}^*$ to $Z_s$ and respects the tilings
$\Delta_{u}$ and $\Delta_s$.
\end{lemma}

\startproof
By Statement 2 of Calculation 6,
and rotational symmetry,
the piecewise similarity
\begin{equation}
h_s=f_s \circ \omega_s
\end{equation}
maps $Z^*_{u}$ to
$Z_s$.  Thanks to Statement 1 of Calculation 6,
the tile $\tau_{u}$ and its rotated image
are really tiles of $\Delta_{u}$.
Thanks to Calculation 5,
the map $h_s$ maps the tiling
$Z^*_{u} \cap \Delta_{u}$ to
the tiling $Z_s \cap \Delta_s$.
\endproof

\begin{lemma}
\label{proof1B}
$h_s$ conjugates
$f_{u}|Z_{u}^*$ to
$f_s|Z_s$. 
\end{lemma}

\startproof
 Choose some point
$r_1 \in Z_{u}^*$.  Let
$p_1=h_s(r_1) \in Z_s$.
Let $p_n$ be the first return of
the forward $f_s$-orbit of $p_1$ to $Z_s$.
So, $p_2,...,p_{n-1}$ do not belong
to $Z_s$.
Define
\begin{equation}
q_j=f_s^{-1}(p_j), \hskip 30 pt j=1,n.
\end{equation}

By Statement 2 of Calculation 6, we have
$q_n \in \omega_s(Z_{u}^*)$. Define
\begin{equation}
r_n=\omega_s^{-1}(q_n).
\end{equation}
By Calculation 5, the
point $r_n$ lies in the forward
$f_{u}$-orbit of $r_1$.  To finish
our proof, we just have to show that
$r_n$ is the first return
of this orbit to $Z^*_{u}$.
If this is false, then there is some earlier
point $r_k \in Z^*_{u}$. But
then $h_s(r_k)=q_m \in Z_s$ for
some $m=2,...,(n-1)$. This is a contradiction.
\endproof

\begin{lemma}
\label{proof1C}
Any nontrivial $f_s$-orbit,
except those contained in
$\tau_s \cup \iota(\tau_s)$,
intersects $Z_s$.
\end{lemma}

\startproof
Consider first the orbit of a point
$p \in W_s$.  Let
$q=\omega_s^{-1}(p)$.   Since
the Main Theorem is true for
the parameter $u$, the orbit of
$q$ intersects $Z_u^*$.  But then,
by Calculation 5, the orbit of $q$
intersects $\omega_s(Z_u^*)=f_s^{-1}(Z_s).$
But then $f_s(q) \in Z_s$.

It remains to consider the orbit of
a point $p \in X_s-Z_s-W_s$.
If $p \in \tau_s \cup \iota(\tau_s)$,
there is nothing to prove.
Otherwise, Statement 3 of Calculation 6
finishes the proof.
\endproof

\subsection{The Second Modular Symmetry}

Now we turn our attention to the proof of
Lemma \ref{main2}.  We re-use some of
the notation from the other case.  Define
\begin{equation}
\label{modular2}
T=\frac{3x-2}{2x-1}, \hskip 30 pt
\omega_s(x,y)=(2s-1)(x,y) \pm (2s-2).
\end{equation}
$T$ here is the inverse of the one
in Equation \ref{modular}.
We take $s \in [3/4,1)$, so that
$u=T(s) \in [1/2,1)$.  The domain
of $\omega_s$ is the set $Y_u$
from the Main Theorem.  
The $(+)$ option for $\omega_s$ 
is taken when $x<0$ and
the $(-)$ option is taken then $x>0$.
We set $W_s=\omega_s(Y_u)$.  

\begin{center}
\resizebox{!}{1.2in}{\includegraphics{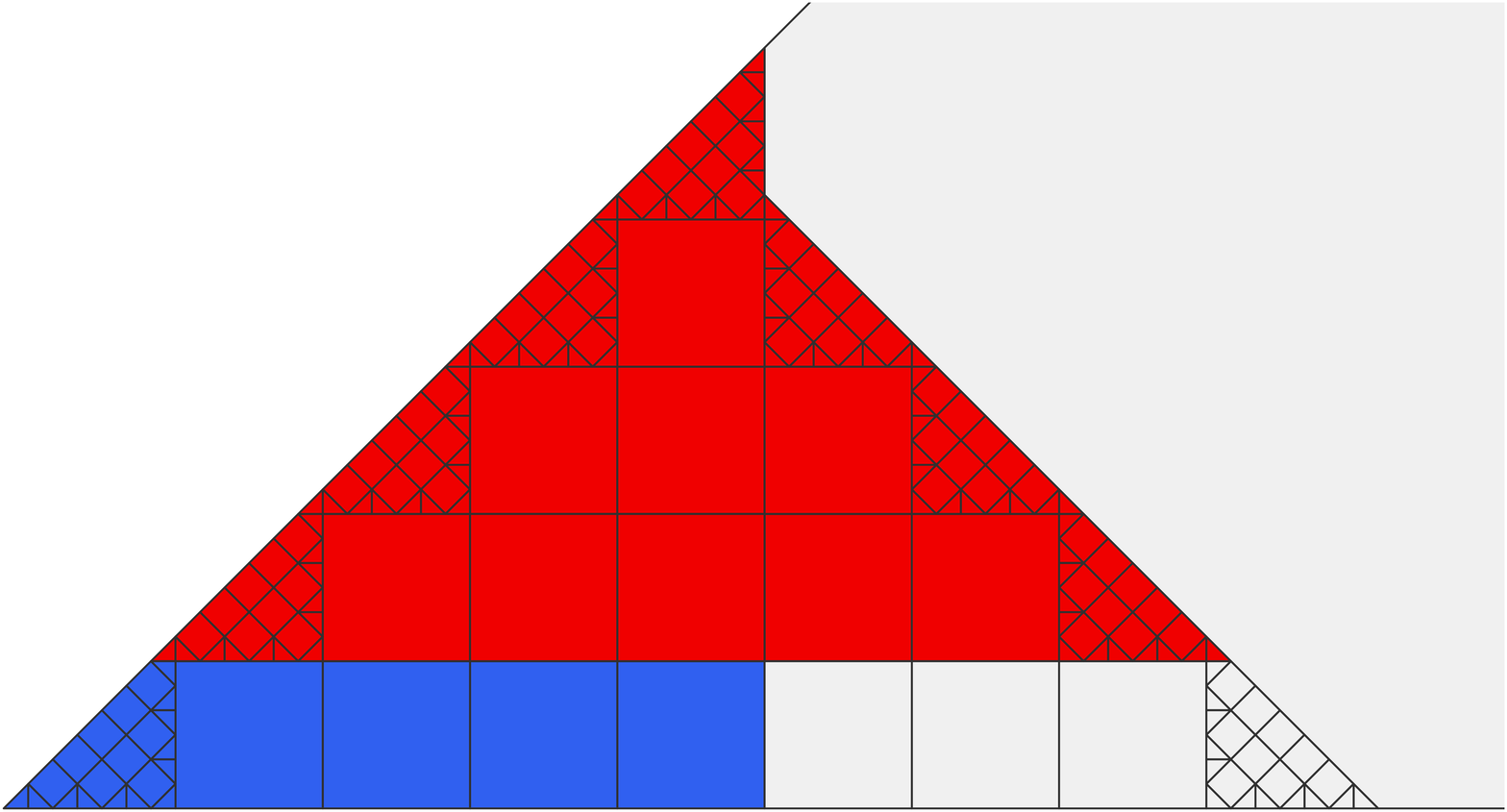}}
\newline
{\bf Figure 5.3:\/} 
Half of $W_s^0$ for $s=28/31$.
\end{center}

\begin{center}
\resizebox{!}{1.2in}{\includegraphics{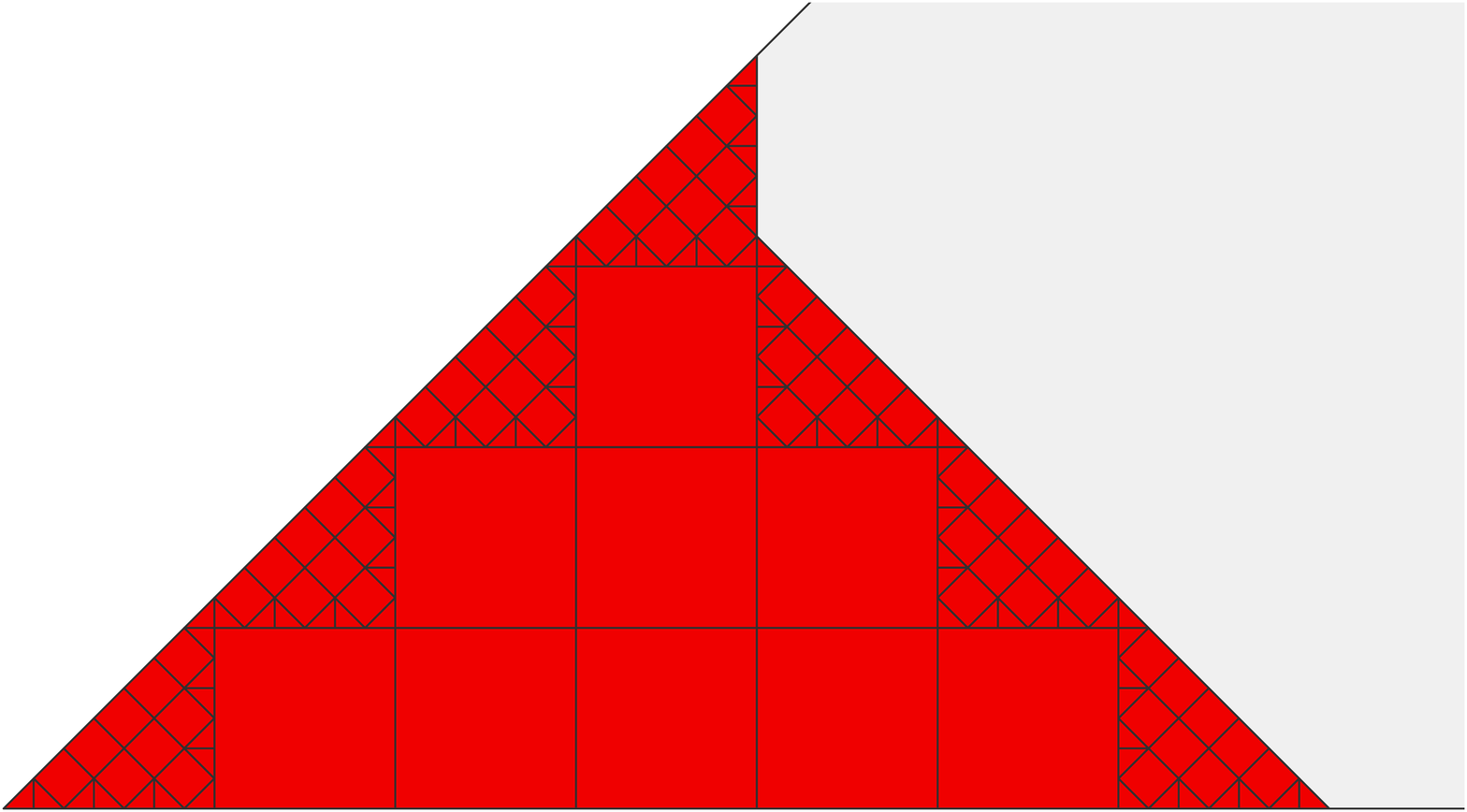}}
\newline
{\bf Figure 5.4:\/} 
 $\Delta_{u} \cap X_{u}^0$ for $u=T(28/31)=22/25$.
\end{center}

The rest of the definitions are done exactly
as in the previous section. The main
difference here is that the tile
$\tau_s$ and $\tau_u$ belong to
$Z_s$ and $Z_u$ respectively.
As above, $\delta_s$ is the vector which spans the diagonal
of $\tau_s$, pointing from the bottom left vertex
to the top right vertex.
Calculation 7 is the calculation parallel to Calculation 5.

\begin{lemma}[Calculation 7]
\label{calc7}
Let $s \in [3/4,1)$ and $u=T(s)$.
Then $\omega_s$ conjugates $f_{u}|Y_{u}$ to
$f_s|W_s$.  Moreover, $Z_s$ is a clean set.
\end{lemma}

Calculation 8 is the calculation parallel to Calculation 6.
Notice that there are some differences.
Item 1 of Calculation 6 refers to $\tau_u$.
We will discuss the reason for this difference below.
Item 2 of Calculation 8 refers to the modified set
\begin{equation}
(Z_s^0)^*=Z_s^0-\tau_s.
\end{equation}
We have to chop off the tile $\tau_s$ to make the statement true.
Also, $f_s$ appears in Item 2 of Calculation 8 whereas
$f_s^{-1}$ appears in Item 2 of Calculation 6.
Item 3 is a slightly different statement, but the
new statement works the same way in the proof of
Lemma \ref{proof1C}.

\begin{lemma}[Calculation 8]
\label{calc8}
Let $s \in [3/4,1)$.
Then
\begin{enumerate}
\item $\tau_s$ is a tile of $\Delta_s$, having period $2$.
\item $f_s(p)=p+\delta_s$ for all $p \in (Z_s^0)^*$.
\item $f_s(X_s - Z_s - W_s) \subset Z_s^*$.
\end{enumerate}
\end{lemma}

\subsection{Proof of Lemma \ref{main2}}

Let $T$ be the map from Equation \ref{modular2}.

\begin{lemma}
\label{orb2}
Let $s \in (3/4,1)$ be any point.  Then
there is some positive $k$ such that
$T^k(s) \in (1/2,3/4)$.
\end{lemma}

\startproof
Same proof as Lemma \ref{orb}.
\endproof

The rest of the proof of Lemma
\ref{main2} is like what we did for
Lemma \ref{main1}, but there are some
small differences. First of all, this
time we know that
$\tau_s$ is a tile of $Z_s$, so 
all the nontrivial orbits intersect
$Z_s$.

\begin{center}
\resizebox{!}{1.5in}{\includegraphics{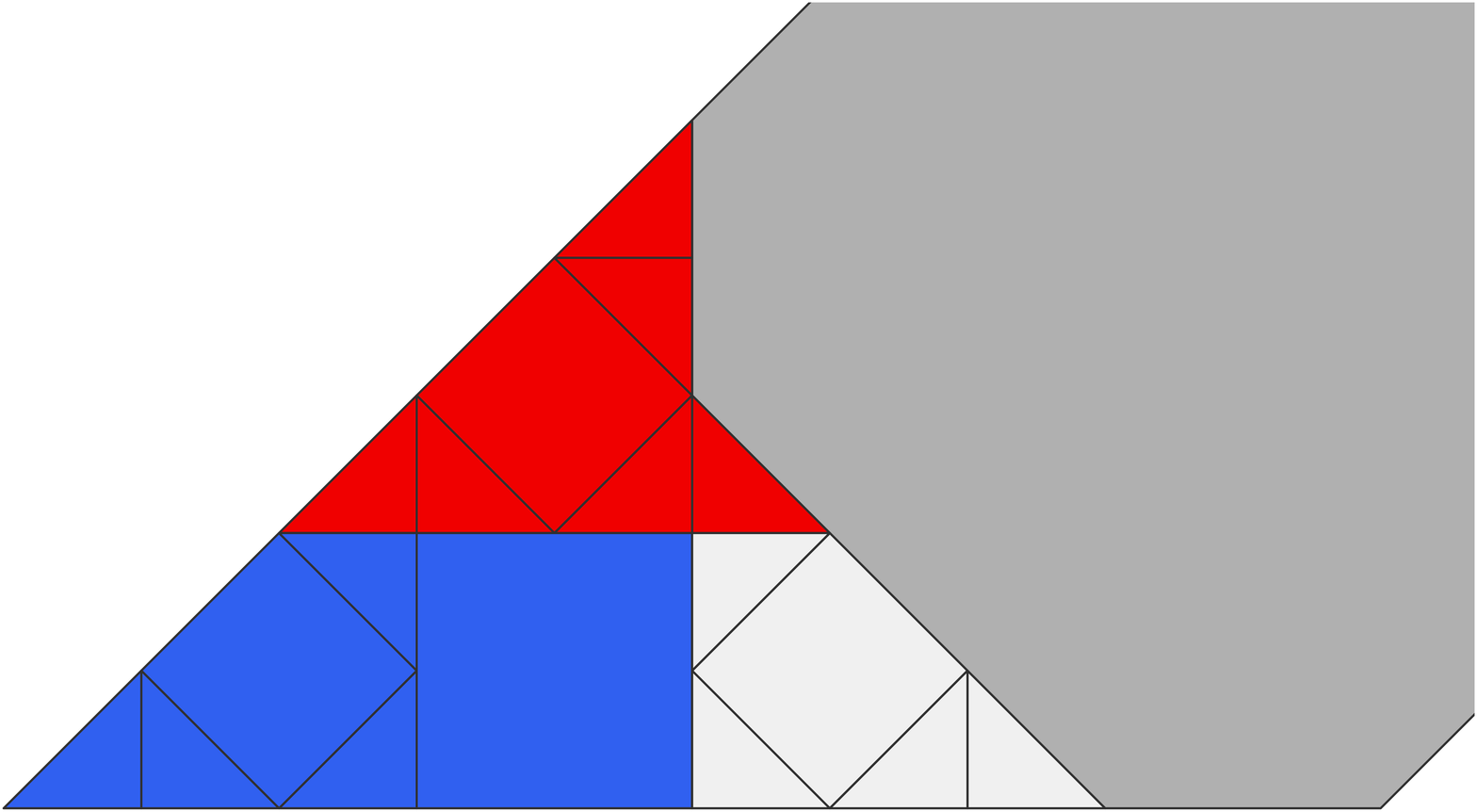}}
\newline
{\bf Figure 5.5:\/} 
Half of $W_s$ for $s=4/5$.
\end{center}

\begin{center}
\resizebox{!}{1.5in}{\includegraphics{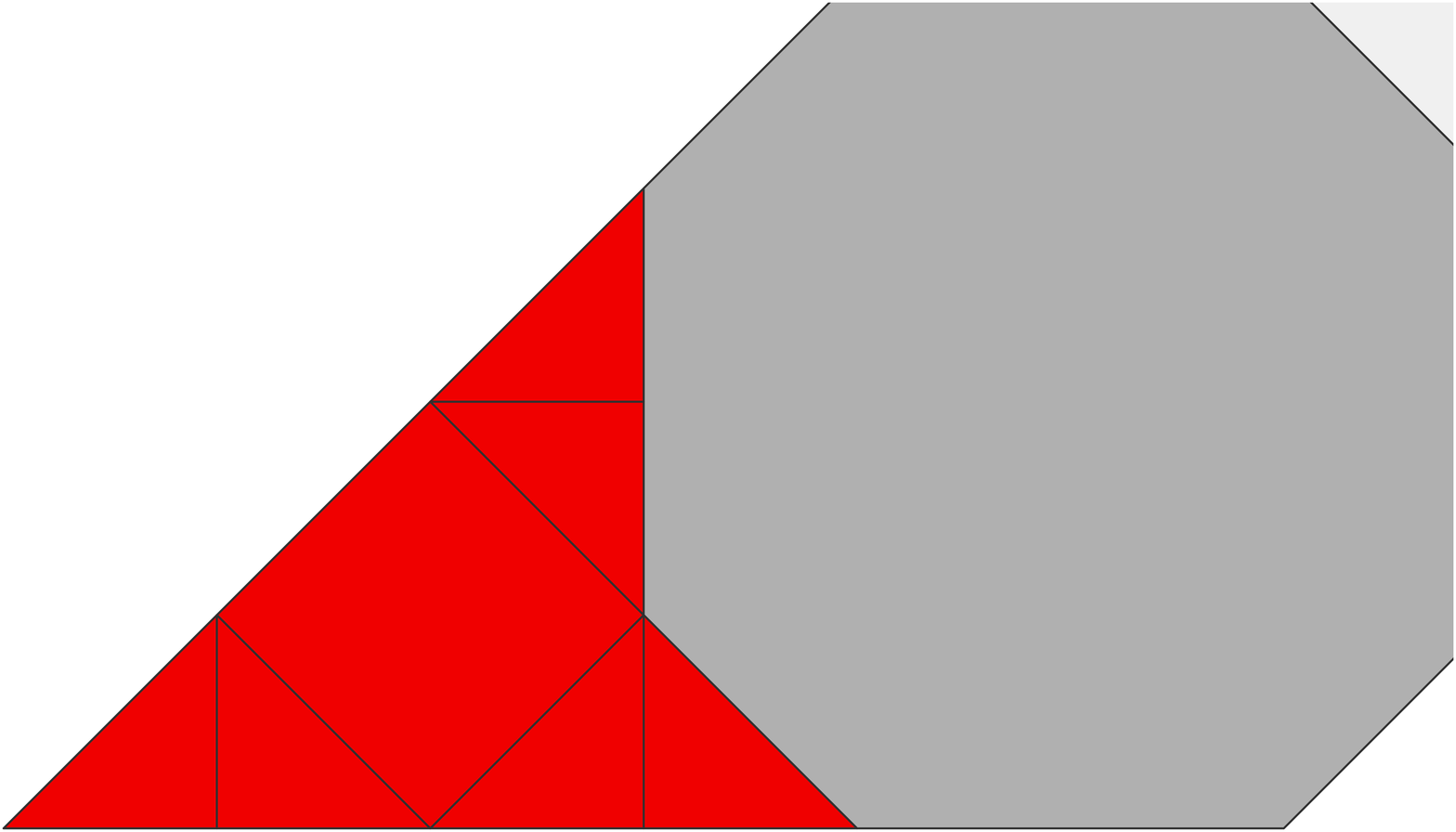}}
\newline
{\bf Figure 5.6:\/} 
 $\Delta_{u} \cap X_{u}^0$ for $u=2/3$.
\end{center}

When $s \in [3/4,5/6)$ and
$u \in [1/2,3/4)$ the tile $\tau_u$ does not exist.
Figures 5.5 and 5.6 show the example of $s=4/5$
and $u=2/3$.  However, we define
\begin{equation}
Z_s^*=(Z_s^0)^* \cup \iota\big((Z_s^0)^*\big) =
Z_s -\tau_s - \iota(\tau_s).
\end{equation}
and we use the
pair $(Z_s^*,Z_u)$ in place of
the pair $(Z_s,Z_u^*)$.  Making these changes,
Lemmas \ref{proof1A}, \ref{proof1B},
and \ref{proof1C} go through without a problem.

\newpage

\section{Properties of the Tiling}

\subsection{The Rational Case}
\label{rat1}

Even though we mainly care about the irrational
case of Theorem \ref{one}, 
we find it convenient to prove the analogous
result for the rational case.   When $s=1/2n$ we
define 
\begin{equation}
\phi_s(x,y)=\bigg(\frac{x}{2n},\frac{y}{2n}\bigg).
\end{equation}
This is a new definition, because the Main Theorem
does not apply for $s=1/2n$.
When $s=0$, we slightly abuse notation and
define $O_s$ to be (simultneously) each of the
$4$ isosceles triangles having vertices
$(0,0)$ and $(\pm 1,\pm 1)$.

\begin{lemma}
\label{st1}
When $s \in (0,1)$ is rational,
a polygon arises in $\Delta_s$ if and only if it is
translation equivalent to
$T_n(O_{s_n})$ for some $n$.
\end{lemma}

\startproof
Our proof goes by induction on the value of $n$ such
that $T^n(s)=0$.  

Consider the base case.  This corresponds to
$s=1/2n$, $n=1,2,3...$.
For $s=1/2$, we check that $\Delta_s$ consists
of $O_{1/2}$ and each of the $4$ triangles $O_0$
scaled down by a factor of $2$.  When $a=1/2n$,
we apply the Insertion Lemma and verify that the
$4$ triangles in $\Delta_{1/2n}$ have sizes consistent
with the statement that they are translates of
$T_0(O_0)$.

The general case follows from induction.
Let $s=s_0$.  By the Main Theorem, a tile
appears in $\Delta_{s_0}$ if and only if it is 
either $O_{s_0}=T_0(O_{s_0})$ or if it is
translation equivalent to a tile of the form
$$\phi_{s_0}(\sigma),$$
where $\sigma$ is a tile of
$\Delta_{s_1}$.   By induction, $\sigma$ is translation
equivalent to the tile
$$\phi_{s_1} \circ  \ldots \circ \phi_k(O_k),$$
for some $k$.  Hence, $\tau$ is translation
equivalent to the tile
$$\phi_{s_0} \circ  \ldots \circ \phi_k(O_k).$$
Conversely, all such tiles appear in $\Delta_{s_0}$,
by the same induction argument.
\endproof
\subsection{The Irrational Case}

Now we prove Theorem \ref{one}.

Let $s$ be some irrational parameter.
First of all, the argument in Lemma
\ref{st1} applies {\it verbatim\/} to
show that the a translate of the tile
$$T_k(O_u), \hskip 30 pt u=R^k(s)$$
does indeed appear in $\Delta_s$.
The argument is simply induction on $k$,
combined with the Main Theorem.

We need to prove the converse, showing that
these are the only tiles that appear in
$\Delta_s$ when $s$ is irrational.
Let $\{_ns\}$ be a sequence of rationals
converging to $s$.

$P$ be a tile of $\Delta_s$.
By Lemma \ref{converge}, we can find a tile
$_nP$ of $\Delta_{_ns}$ such that
$_nP \to P$ as $n \to \infty$.  By Lemma
\ref{st1}, there is an integer
$k_n$ such that
$_nP$ is translation equivalent to some
tile of the form
$$_nT_{k_n}(O_{u_n}), \hskip 30 pt
u_n=R^{k_n}(_ns).$$
Here $_nT_0, _nT_1, _nT_2,...$ are the maps
which arise in Lemma \ref{st1} relative to
the parameter $_ns$.

Now, the scale factor of $_nT_k$ tends to
$0$ with $k$, and $_nP$ has uniformly large
diameter.  Therefore, the integer $k_n$ is
uniformly bounded from above.  So, passing
to a subsequence, we can assume that $_nk$
is independent of $n$.  That is,
$P_n$ is a translation of 
\begin{equation}
\label{appx}
_nT_{k}(O_{u_n}), \hskip 30 pt
u_n=R^{k}(_ns).
\end{equation}
But \begin{itemize}
\item $_nT_k \to T_k$,
\item $u_n \to u=R^k(s)$,
\item  $O_{u_n} \to O_u$.
\end{itemize}
Hence $P$ is a translate of
$T_k(O_u)$.

The rest of this chapter is
devoted to proving Theorem \ref{cor}.

\subsection{Existence of Square Tiles}

Here we prove Statement 1 of Theorem \ref{cor}.
We will suppose that $s \not = \sqrt 2/2$ and
show that $\Delta_s$ has a square tile.
We will first suppose that $s$ is irrational.

\begin{lemma}
\label{copy}
Suppose that $s$ and $t=R(s)$ lie in
$(0,1)$.  If $\Delta_t$ has a
non-central square tile, then so
does $\Delta_s$.
\end{lemma}

\startproof
By the Main Theorem,
$\Delta_s$ has a similar copy
of every non-central tile of
$\Delta_t$, and these similar
copies are themselves non-central.
\endproof

\begin{lemma}
Suppose that $s$ is irrational.
If $R^n(s)<1/2$ for some even $s$, then
$\Delta_s$ has a square tile.
\end{lemma}

\startproof
We choose $n$ to be as small as possible.
If $n=0$ then $s<1/2$ and the central tile
of $\Delta_s$ is square.
Suppose $n \geq 2$.  Let
$t=R^{n-2}(s)$ and
$u=R^{n-1}(s)$ and 
$v=R^n(s)$.
We know that $v<1/2$.
If $u>1/2$ then $t<1/2$,
contradicting the minimality of $n$.
So, $u<1/2$.  But then we
apply the Main Theorem to the
pair $(u,R(u))$ and conclude that
some non-central tile of $\Delta_u$
is a square. Now we apply Lemma \ref{copy}
repeatedly to conclude that $\Delta_s$
has a non-central square tile.
\endproof

\begin{lemma}
\label{quarter}
Suppose that $s$ is irrational.
If $R^n(s)<1/4$ for some odd $n$,
then $\Delta_s$ has a square tile.
\end{lemma}

\startproof
If $n=0$ then $\Delta_s$ has a central
square tile.  So, assume that $n>0$.
Let $t=R^{n-1}(s)$ and $u=R^n(s)$.
Since $u<1/4$, the tiling
$\Delta_u$ has at least $3$ central
tiles, one of which is contained
in the set $Y_t \subset \Delta_u$.  But then
$Z_t \subset \Delta_t$ contains a non-central
square tile.  Repeated applications
of Lemma \ref{copy} finish the proof.
\endproof

\begin{lemma}
\label{half}
Suppose that $s$ is irrational and
$\Delta_s$ has no square tiles.
Then $s=\sqrt 2/2$.
\end{lemma}

\startproof
We know that $R^n(s)>1/2$ for every even $n$
and $R^n(s) \in (1/4,1/2)$ for every odd $n$.
If $R^n(s) \in (1/3,1/2)$ for some odd $n$,
then $R^{n+1}(s)<1/2$, because $R$ maps
$(1/3,1/2)$ onto $(0,1/2)$.  Hence
$R^n(s) \in (1/4,1/3)$ for all odd $n$.

For $x \in (1/4,1/3)$ we have
the formula
\begin{equation}
R(x)=\frac{1}{2x}-1
\end{equation}
Hence, if $x>1/2$ and $R(x)\in (1/4,1/3)$ we have
\begin{equation}
R^2(x)=\frac{1-2x}{-2+x}=Mx, \hskip 30 pt
M=\left[\matrix{-\sqrt 2 & 1/\sqrt 2 \cr \sqrt 2 & -\sqrt 2}\right].
\end{equation}
Here $M$ acts as a linear fractional transformation.
Iterating, we have
\begin{equation}
R^{2n}(s)=M^n(s).
\end{equation}
The map $M$, acting as a linear fractional transformation,
fixes $\pm \sqrt 2/2$ and is expanding in a
neighborhood of $\sqrt 2/2$. We conclude that
there is some $n$ such that $M^n(s)<0$ unless
$s=\sqrt 2/2$.  Since we never have
$N^n(s)<0$ we must have $s=\sqrt 2/2$.
\endproof

\begin{lemma}
If $s$ is rational, the $\Delta_s$
has a square tile.
\end{lemma}

\startproof
The proof is the same as in the irrational
case, except we have to worry about what happens
if $R^n(s) \in \{1/4,1/3,1/2\}$ for some odd $n$.
If $R^n(s)=1/4$ for some odd $n$, then
the same argument as in Lemma \ref{quarter}
applies. The point is that $\Delta_{1/4}$
still has $3$ central square tiles.

If $R^n(s)=1/2$ for some odd $n$, then
$R^{n-1}(s)<1/2$.    But then the
argument in Lemma \ref{half} applies
to the even iterate $R^{n-1}(s)$.

If $R^n(s)=1/3$ for some odd $n$,
then $R^{n+1}(s)=1/2$.  But then the
argument in Lemma \ref{half} applies
to the even iterate $R^{n+1}(s)$.
because $\Delta_{1/2}$ still has a
central square tile.
\endproof

Putting together the last several lemmas, we finish
the proof that $\Delta_s$ has a square tile as
long as $s \not = \sqrt 2/2$.

\subsection{The Case of Squares}

Here we prove Statement 2 of Theorem \ref{cor}.
We compare the map $R$ with the {\it Gauss map\/}
\begin{equation}
\label{gauss}
\gamma(s)=1/s-{\rm floor\/}(1/s)
\end{equation}

\begin{lemma}
\label{tech1}
Let $s=[0,a_1,a_2,a_3,...]$ be the
continued fraction expansion of $s$.
Suppose that $a_1$ is even.
Then $s,R(s)<1/2$ and
$R^2(s)=\gamma^2(s)$.
\end{lemma}

\startproof
Let $x=1/s$. We write 
$$x=a_1+\frac{1}{a_2+1/y}.$$
The continued fraction expansion of $1/y$
is $[0,a_3,a_4,...]$.  Hence
$1/y=\gamma^2(s)$.  Note, in particular,
that $1/y<1$.
We need to show that $R^2(s)=1/y$.

Since $a_1>1$, we have $s<1/2$.  Hence 
$$R(s)=1/(2s)-{\rm floor\/}(1/(2s))=x/2-{\rm floor\/}(x/2).$$
We have
\begin{equation}
\label{tech2}
x/2=(a_1/2)+\frac{1}{2a_2+2/y}.
\end{equation}
Hence
$$R(s)=\frac{1}{2a_2+2/y}<1/2.$$
Since $R(s)<1/2$ we have
$$R^2(s)=
(a_2+1/y)-{\rm floor\/}(a_2+1/y)=1/y,$$
as desired.
\endproof

Let $[a_0,a_1,a_2,a_3,...]$ be the continued
fraction expansion for $s$.  We call $s$
{\it oddly even\/} if $s$ is irrational
and $a_k$ is even for all odd $k$.
When $s \in (0,1)$ we have $a_0=0$.

\begin{lemma}
\label{oddly}
Let $s \in (0,1)$ be irrational. Then
$s$ is oddly even if and only if
$R^k(s) \in (0,1/2)$ for
all $k$.
\end{lemma}

\startproof
Suppose that $s$ is oddly even.
By the previous result,
$s$ and $R(s)$ lie in $(0,1/2)$ and
$R^2(s)$ is again oddly even.
Hence $R^2(s)$ and $R^3(s)$ lie
in $(0,1/2)$ and $R^4(s)$ is oddly
even.  And so on.

Conversely, suppose that $s$ is not
oddly even.  Applying Lemma \ref{tech1}
finitely many times if necessary, we
reduce to the case where the first
term in the continued fraction expansion
of $s$ is odd.  If $s>1/2$ we are done.
Otherwise, $s$ lies in one of the
intervals $(1/3,1/2)$, $(1/5,1/4)$,
$(1/7,1/6)$,..., and $R$ maps each
of these intervals onto $(1/2,1)$.
So $R(s)>1/2$ in this case.
Hence, if $R^k(s) \in (0,1/2)$ for all
$k$, then $s$ is oddly even.
\endproof

Lemma \ref{oddly} combines with Theorem \ref{one} to
prove Statement 2 of Theorem \ref{cor} in case
$s \in (0,1)$.  But the case $s>1$ now follows from
the Inversion Lemma and from the fact that the
map $x \to 1/(2x)$ preserves the set of oddly
even numbers.

\subsection{The Density of Shapes}

Now we prove Statement 3 of Theorem \ref{cor}.
This result follows immediate from Theorem
\ref{one} and from the following lemma.

\begin{lemma}
\label{dens}
For almost all $s \in (0,1)$, the 
orbit $\{R^n(s)\}$ is dense in $(0,1)$.
\end{lemma}

We will prove Lemma \ref{dens} through
a series of smaller lemmas.  We start with
a classic result. See [{\bf BKS\/}] for instance.

\begin{lemma}
Almost every orbit of the Gauss map is dense in
$(0,1)$.
\end{lemma}

This result has a well known geometric consequence.
Let $\Sigma_0$ denote the trice-punctured sphere.
Let $T_1(\Sigma)$ denote the unit tangent
bundle of $\Sigma$.  A geodesic in $\Sigma_0$
has a natural {\it lift\/} to the unit
tangent bundle: One just keeps track of
the points on the geodesic as well as their
unit tangent vectors.
We say that a geodesic on $\Sigma_0$
{\it emanates\/} from a cusp of $\Sigma$
if one end of the geodesic is asymptotic
with a cusp.  The lift of such a geodesic
to $\H^2$ has one endpoint on a parabolic
fixed point of the surface fundamental group.
There is a natural measure on the set of
cusps emanating from one of the cusps of
$\Sigma_0$.

\begin{corollary}
Almost every geodesic ray emanating from a
cusp of $\Sigma_0$ lifts to a dense subset of
$T_1(\Sigma_0)$.
\end{corollary}

\startproof
Let $\alpha$ be a geodesic emanating from
a cusp of $\Sigma_0$.  By symmetry, it suffices
to consider the case when some lift
$\widetilde \alpha$ of $\alpha$ is the vertical
geodesic connecting $\infty$ to some
$r \in (0,1)$.
The position of any given
segment of $\widehat \alpha$ is determined,
to arbitrary precision, by finite portions
of the orbit of $r$ under the Gauss map.
If $\alpha$ is
chosen so that this orbit is dense, then we
can approximate any finite geodesic segment
on $\Sigma_0$, up to an arbitrarily small
error, using a portion of $\alpha$.
\endproof

Now let $\Sigma$ be a finite normal covering surface
of $\Sigma_0$.  That is,
$\Sigma_0=\Sigma/G$, where $G$ is
some finite group acting on $\Sigma$.
There is again a natural
measure on the set of geodesics emanating
from a cusp of $\Sigma$.

\begin{corollary}
Almost every geodesic ray emanating from the
cusp of $\Sigma$ lifts to a dense subset of
$T_1(\Sigma)$.
\end{corollary}

\startproof
Let $\alpha$ be a geodesic emanating from our cusp
of $\Sigma$.  Let $\overline \alpha$ be the
projection of $\alpha$ to $\Sigma_0$.
Almost every choice of $\alpha$
leads to $\overline \alpha$ having a dense lift 
in $T_1(\Sigma_0)$. But then the $G$-orbit
of the lift of $\alpha$ is dense in $T_1(\Sigma)$.  Let
$C$ be the closure of the lift of
$\alpha$ in $T_1(\Sigma)$.
We know that $G(C)=T_1(\Sigma)$.  At the same
time, we know that $\overline \alpha$ 
approximates, with arbitrary precision, any
closed loop in $\Sigma_0$.  From this we
see that in fact $C$ is $G$-invariant.
Hence $C=T_1(\Sigma)$, as desired.
\endproof

Let $T$ be the $(2,4,\infty)$ hyperbolic triangle
generating our group $\Gamma$.  
One of the edges of $T$ is contained in the
geodesic circle $C$ fixed pointwise by the
map $z \to 1/(2\overline z)$. We color
$C$ red.   We color the other two edges
blue.  We then lift this coloring to the
universal covering $\H^2$.  This gives
us a pattern of geodesics that is invariant
under the $(2,4,\infty)$ triangle group.  Among other
colored geodesics in $\H^2$, we have
the red circle, and the blue lines connecting
$\infty$ to half-integers.

\begin{lemma}
\label{ae2}
Almost every geodesic emanating from the
cusp of the $(2,4,\infty)$ triangle
has a billiard trajectory which hits the
red edge in a dense set of points.
\end{lemma}

\startproof
We think of $T$ as the $(2,4,\infty)$ orbifold.
There is a surface $\Sigma$ (a $4$-times punctured sphere)
which covers $T$ in the sense of orbifolds, and which
also covers $\Sigma_0$, the thrice punctured sphere.
Say that a geodesic on $\Sigma$ is good if it emanates
from a cusp of $\Sigma$ and lifts to a dense set
in $T_1(\Sigma)$.  Almost every geodesic on
$T$, emanating from the cusp of $T$, has a
preimage which is a good geodesic.  Hence,
almost every geodesic emanating from the cusp
of $T$ has dense image in the unit tangent
bundle of $T$.  But such a geodesic would intersect
each edge of $T$ in a dense set of points.
\endproof

To prove that $\{R^n(s)\}$ is dense for
almost every choice of $s \in (0,1)$, it
suffices to prove that
$\{R^n(s)\}$ is dense for almost every
choice of $s \in (0,1/2)$.
It is also useful to first consider
the alternate map $R_1: (0,1/2) \to (0,1/2)$,
defined as follows:
\begin{itemize}
\item $R_1(s)=R(s)$ if $R(s)<1/2$.
\item $R_1(s)=1-R(s)$ if $R(s)>1/2$.
\end{itemize}

One can describe $R_1$ like this. Starting
with $s_0$, we first reflect in the red
circle $C$ to produce the point $s_1$.
There is some nearest vertical blue line
which separates $s_1$ from $0$.
We reflect in this blue line to produce
$s_2$.  We now repeat these reflections
in blue lines until we arrive at
a point in $(0,1/2)$, and this point is $R_1(s)$.

Each $s \in (0,1/2)$ corresponds to a geodesic
$g_s$ of $T$, which emanates from the cusp. The recipe
is that the lift of $g_s$ to $\H^2$ is the
geodesic connecting $s$ to $\infty$. From the
description of $R_1$ above, we see that we
can recover the action of $R_1$ by looking
at the billiard path of $g_s$.  The
orbit $\{R^1(s)\}$ is dense provided that
$g_s$ intersects the red edge in a dense set.
By Lemma \ref{ae2}, this happens for almost
all $s \in (0,1/2)$.

Knowing that $\{R_1^n(s)\}$ is dense is not quite the
same as knowing the $\{R^n(s)\}$ is dense.  However,
there is a dcomposition of the red edge of $T$ into
intervals $I_1,J_2,I_3,J_4,...$ such that an intersection
point of $g_s$ with $I_j$ corresponds to a value
$s_k$ where $R_1(s_k)=R(s_k)$ and an intersection point
of $g_s$ with $J_j$ corresponds to a value of
$s_k$ where $R_1(s_k)=1-R(s_k)$.  In terms of
billiards, the $I$ intervals are such that
$g_s$ hits an even number of blue edges after
hitting some $I$-interval.   

The upshot of this interval decomposition is that,
since $g_s$ intersects both the $I$-intervals
and the $J$-intervals densely, the orbit
$\{R^n(s)\}$ is also dense.

\newpage

\section{Covering Results}

\subsection{An Area Estimate}

Recall that $X_s^0$ is the portion of
$X_s$ to the left of the central tiles.
Define
\begin{equation}
\lambda(s)=\frac{{\rm Area\/}(\Delta_s \cap X_s^0)}{{\rm Area\/}(X_s^0)}.
\end{equation}

\begin{lemma}[Area]
The function $\lambda(s)$ is uniformly
bounded away from $0$, for all
$s \in (0,1)$.
\end{lemma}

\startproof
By the Insertion Lemma, it suffices
to take $s \in [1/4,1)$.  The
case $s=1/4$ is trivial, so
we consider $s \in (1/4,1)$. 
When $s$ is bounded away from $1/2$ and $1$,
the size of the
largest square in $X_s^0$ is bounded
away from $0$. (Figure 7.1 below shows a picture
of the case when $s$ is very near $1/4$, the
other potential place to worry about.)

As $s \to 1$, the size of $\tau_s$ tends to $0$.
However, the inductive argument 
given to prove Lemma \ref{main1}
shows that the union of squares
isometric to $\tau$ forms the following 
pattern. There is
a bottom row of $2k+1$ squares, then a row of
$2k-1$ squares, then a row of
$2k-3$ squares, and so on, all the way down
to a single top square.
When $s>1$, the number $k$ is such that $T^k(s) \in (5/4,3/2)$.
When $s<1$, the number $k$ is such that $T^k(s) \in (3/4,5/6)$.

Technically, the argument we gave above
establishes the existence of the 
left half of the
picture (and the central column.)
Combining the Inversion Lemma and
the Bilateral Lemma, we see that the reflection
$\rho$ in the vertical line through the top vertex
of $X_s^0$ maps the right half of $\Delta_s \cap X_s^0$
into the left half.  This lets us deduce that
the right half of $X_0^s$ has the same pattern
of squares as the left half.

As $s \to 1$, the bottom row is
more and more nearly filled up with squares
isometric to $\tau$. Hence, the ``triangular pile''
of squares fills more and more of 
$X_s^0$ as $s \to 1$.  Indeed, $\lambda(s) \to 1$
as $s \to 1$.

When $s \to 1/2$, we can again apply the Inversion Lemma
to instead consider the case $s \to 1$, which we have
already treated.
\endproof

\noindent
{\bf Remark:\/}
The Area Lemma  is really
what is responsible for our
result that $\widehat \Lambda_s$ always has
measure $0$.

\subsection{A Length Estimate}

Say that a {\it special edge\/} is an edge of
$X_s^0$ which is contained either in the
bottom edge or the left edge of $X_s$.
Let $b_s$ and $\ell_s$ be the bottom and left edges of $X_s^0$
respectively.  We will first focus on $\ell_s$.  When the
dependence on the parameter is clear, we will set
$\ell=\ell_s$ and $b=b_s$.
Given a special edge $E$ of $X_s^0$, define
\begin{equation}
\lambda(E,s)=\frac{{\rm length\/}(E \cap \Delta_s)}{{\rm length\/}(E)}
\end{equation}
This equation needs some interpretation because $\Delta_s$,
being the union of open tiles, is technically disjoint from
$E$.  What we are talking about here is the portion of $E$
contained in edges of tiles of $\Delta_s$.

\begin{center}
\resizebox{!}{1.2in}{\includegraphics{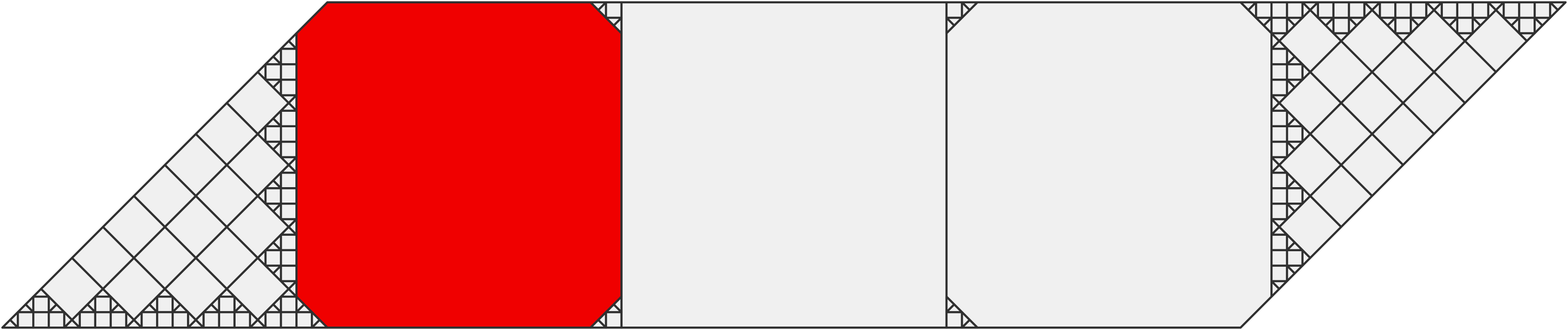}}
\newline
{\bf Figure 7.1\/} $\Delta_s$ for $s=21/80$.
\end{center}

The area estimate given above has an analogue for
lengths, but the result is more subtle.  Figure
6.1 indicates some of the subtlety.
Figure 7.1 shows the tiling for a parameter very
near $1/4$.  The edge $\ell$ is very nearly
covered by the edges of square tiles, whereas
the right half of $b$ is nearly covered by
the edge of a single tile and the left half
meets the tiles in a jagged \footnote{Technically,
in the rational case, the special edges are
completely covered by tile edges, but Figure 7.1
is supposed to suggest the kind of complicated
phenomena which can arise in the irrational
case.} kind of way.
This ``jagged edge'' phenomenon is what makes
our length estimate more subtle than the area
estimate.

\begin{lemma}[Length]
Let $s \in (0,1)$.  Then
$\lambda(\ell,s)$ and
$\lambda(b,s)$ are uniformly
bounded away from $0$ as long
as $1/(2s)$ mod $1$ lies in a
compact subset of $(0,1)$.
\end{lemma}

We will prove the Length Lemma through a series of smaller
lemmas.

\begin{lemma}
There is a uniform lower bound on
$\lambda(\ell,s)$ for $s \in [1/4,2/3]$.
\end{lemma}

\startproof
Our analysis for the case of area gives us a uniform
lower bound for $(b_s,s)$ when $s \in [3/4,3/2]$.  By the
Inversion Lemma, we get the same result for $(\ell_s,s)$ 
when $s \in
[1/3,2/3]$. So, we just have to worry
about $s \in [1/4,1/3]$.

When $s \in [1/4,1/3]$ we have $t=R(s)>1/2$.
Let $O_t$ be the trivial octagonal tile of
$\Delta_t$.
By Statement 3 of the Main Theorem, 
$\phi_s$ extends to $O_t$ and
$\phi_s(O_t)$ has an edge in $\ell_s$.
As long as $t$ is bounded away from $1$,
which corresponds to $s$ being bounded
away from $1/4$, there is a lower bound
to the length of $\phi_s(O_t) \cap \lambda_s$.
Again, we interpret $O_t$ as a closed tile
for this statement.

Consider what happens when $s \to 1/4$ from
above. In this case $t \to 1$ from below.
We check that $Z_s$ nearly covers all
of $\ell_s$.  All that is missing is a small
portion of the very top of $\ell_s$,
coming from $\phi_s(O_t)$. See Figure 7.1.
So, by the Main Theorem, and our analysis
in the Area Lemma, $\lambda(\ell_s,s) \to 1$
as $s \to 1/4$ from above.
\endproof

\begin{lemma}
$\lambda(\ell,s)$ is uniformly bounded
away from $0$ on any compact subset
of $[2/3,1)$.
\end{lemma}

\startproof
When $s=2/3$, there is a large square diamond
whose edge lies in $\ell_s$.  This tile is
stable, as one can see from the arithmetic
graph, so $\lambda(\ell,s)$ uniformly
bounded away from $0$ in a neighborhood of $2/3$.
As long as $t$ is bounded away from $0$,
the tile $\tau_t$ of $\Delta_t$ sharing a vertical
edge with the trivial tile has an edge in 
$\beta_s$ whose length is bounded away from $0$.
So, by the Main Theorem, we get a uniform lower
bound on $\lambda(\ell_s,s)$ for
$s \in (2/3,1)$ as long as $s$ stays away from $1$.
\endproof

Combining the the results above, we see that
$\lambda(\ell,s)$ is uniformly bounded
away from $0$ when $s \in (0,1)$ is
uniformly bounded away from $1$.
This result certainly implies the
result for $\ell_s$ claimed in
the Length Lemma.   Now we turn
our attention to the edge $b_s$.

\begin{lemma}
There is a uniform lower bound on
$\lambda(b,s)$ when $s$ lies 
in a compact subset of $(1/2,1)$.
\end{lemma}

\startproof
We have a uniform lower bound on
$\lambda(\ell,s)$ for $s$ in
a compact subset of $(1/2,1)$.
Now we apply the Inversion Lemma, which
maps $s$ to $1/2s$ and switches the roles of
$\ell$ and $b$.
\endproof

\begin{lemma}
There is a uniform lower bound on
$\lambda(b,s)$ when $s$ lies in 
a compact subset of $(1/4,1/2)$.
\end{lemma}

\startproof
Let $t=R(s)<1$.
The left edge of $Y_t$ is
all of $\ell_t$.  When $s \in (1/4,1/2)$,
the scale factor of $\phi_s$, from the
Main Theorem, is bounded uniformly away
from $0$.  Hence, the bottom edge of
$Z_s$ is uniformly large.
 As long as
$t$ is bounded away from $1$, the
trivial tile in $\Delta_t$ has edges
which take up a uniformly large
fraction of the top and bottom edges
of $X_t$.  By the Main Theorem, the
corresponding tile in $\Delta_s$ has
a uniformly large edge in $b_s$.
\endproof

The Length Lemma for $b_s$ now follows
from the Insertion Lemma and the preceding
two results.

\subsection{The Covering Lemma}

Let $s$ be some parameter, and let $t$ be some
other parameter.
We say that a {\it similar copy\/}
of $X_t^0$ is a set of the form
$T(X_t^0)$ where $T$ is a similarity.
We say that 
$A=T(X_t^0)$ {\it fits nicely\/} over $X_s$
if 
\begin{equation}
T^{-1}(\widehat \Lambda_s \cap A) \subset \widehat \Lambda_t.
\end{equation}
We don't require that
$A$ is actually a subset of $X_s^0$. 

Say that a {\it nice covering\/} of
$X_s^0$ is a covering of the form
\begin{equation}
\label{cov0}
X_s^0 \subset A_1 \cup ... \cup A_m \cup B_1 \cup ... \cup B_n,
\end{equation}
where each $A_i$ is a similar copy of $X_t^0$ which
fits nicely over $X_s^0$ and each $B_j$ is a periodic
tile of $\Delta_s$. 

\begin{lemma}[Covering]
Let $t=R(s)$.
The set $X_s^0$ has a nice covering by
similar copies of $X_t^0$.
All the similarities associated to the pieces
in the cover have the form $I \circ \phi_s$, where
$I$ is an isometry and $\phi_s$ is as in the
Main Theorem.  
\end{lemma}

We will prove this result through a series
of smaller lemmas.

\begin{lemma}
\label{c0}
The Covering Lemma holds for $s \in (1/4,1/3)$.
\end{lemma}

\startproof
Let $A_s$ and $B_s$ be the hexagon and
triangle from the Bilateral Lemma.
For $s \in (1/4,1/3)$ we have $R(s) \in (1/2,1)$.
For these values, a calculation shows that
$B_s \subset Z_s$.  The largest tile in $Z_s$ is
an octagon, $\tau$, the image of the central tile
in $\Delta_t$ under $\phi_s$.  
Figure 7.4 below shows the typical example
of $s=3/10$.

\begin{center}
\resizebox{!}{2.4in}{\includegraphics{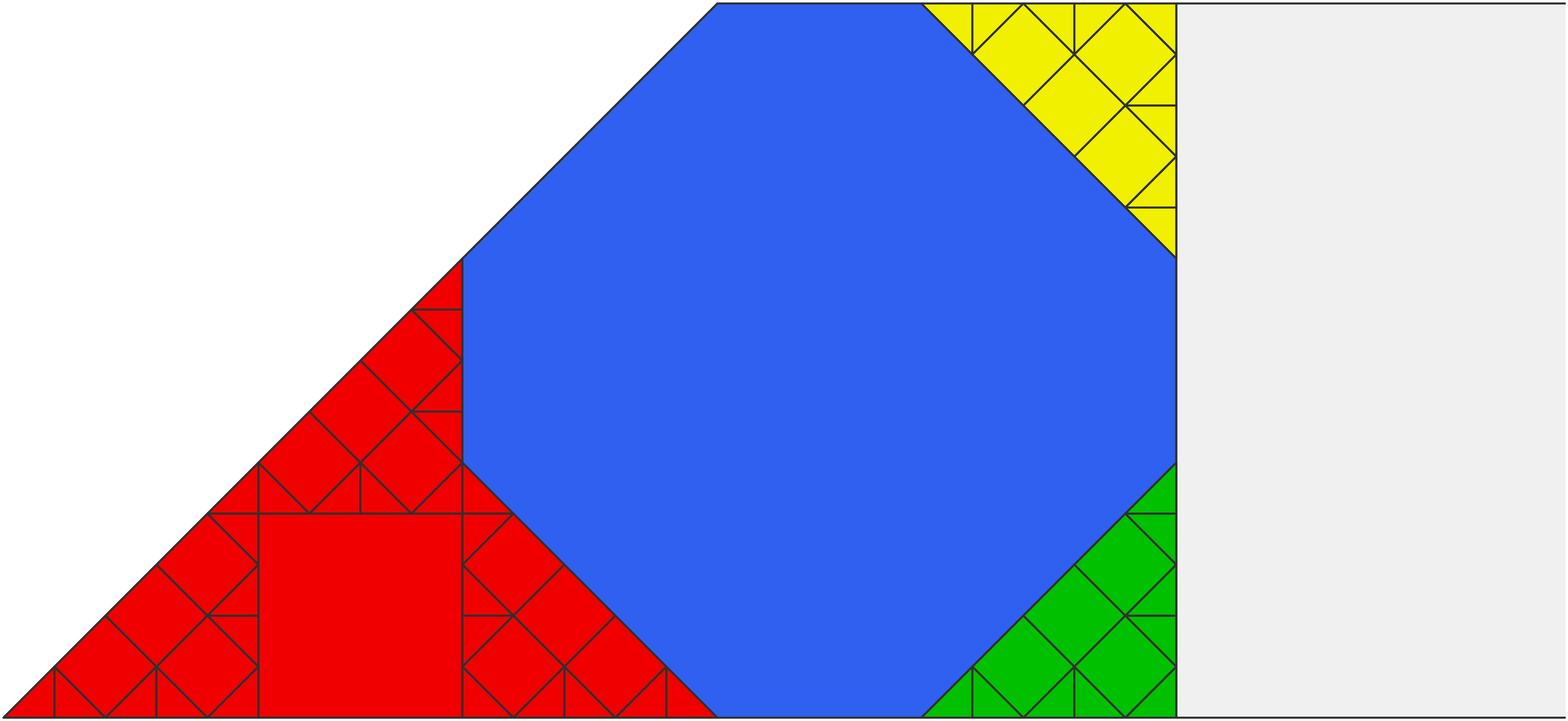}}
\newline
{\bf Figure 7.2\/} $X_s$ (red) and $\tau$ (blue)
and $X_s^0-Z_s-\tau$ (yellow, green)
\end{center}

Let $\rho_1$ be the reflection in the line $D_s$ from
the Bilateral Lemma.
The tiles in the collection 
$\rho_1(\Delta_s \cap Z_s) \cap X_s$
are again tiles of $\Delta_s$ by the Bilateral Lemma, and
$\rho_1(Z_s) \cap X_s$ is precisely the
region of $X_s^0$ lying above $\tau$.  This is the yellow
region in Figure 7.2.  Hence
$\rho_2(Z_s)$ fits nicely over $\Delta_s$. 

Let $\rho_2$ be reflection in the line $H_s$ from the
Bilateral Lemma.  By the Bilateral Lemma,
$\rho_2 \circ \rho_1(Z_s)$ fits nicely over $\Delta_s$. The
intersection $\rho_2 \circ \rho_1(Z_s) \cap X_s$ is
precisely the region of $X_s^0$ lying below and to the
right of $\tau$.  This is the green region in
Figure 7.2.

To finish the proof, we note that $Z_s=\phi_s(X_t^0)$.
The point is that there is just one central tile
in $\Delta_t$.
So, all the sets we are using have the form
$I \circ \rho_s(X_t^0)$, where $I$ is an isometry.
\endproof

\begin{corollary}
The Covering Lemma holds for $s \in (3/2,2)$.
\end{corollary}

\startproof
For $s$ in this range, we have $R(s)=s-1$.
The result now follows from the Inversion Lemma.
\endproof

\begin{lemma}
\label{c1}
The Covering Lemma holds for $s \in (1/3/2/5)$.
\end{lemma}

\startproof
Figure 7.5 shows a typical case, for the parameter
$s=19/50$.

\begin{center}
\resizebox{!}{2.2in}{\includegraphics{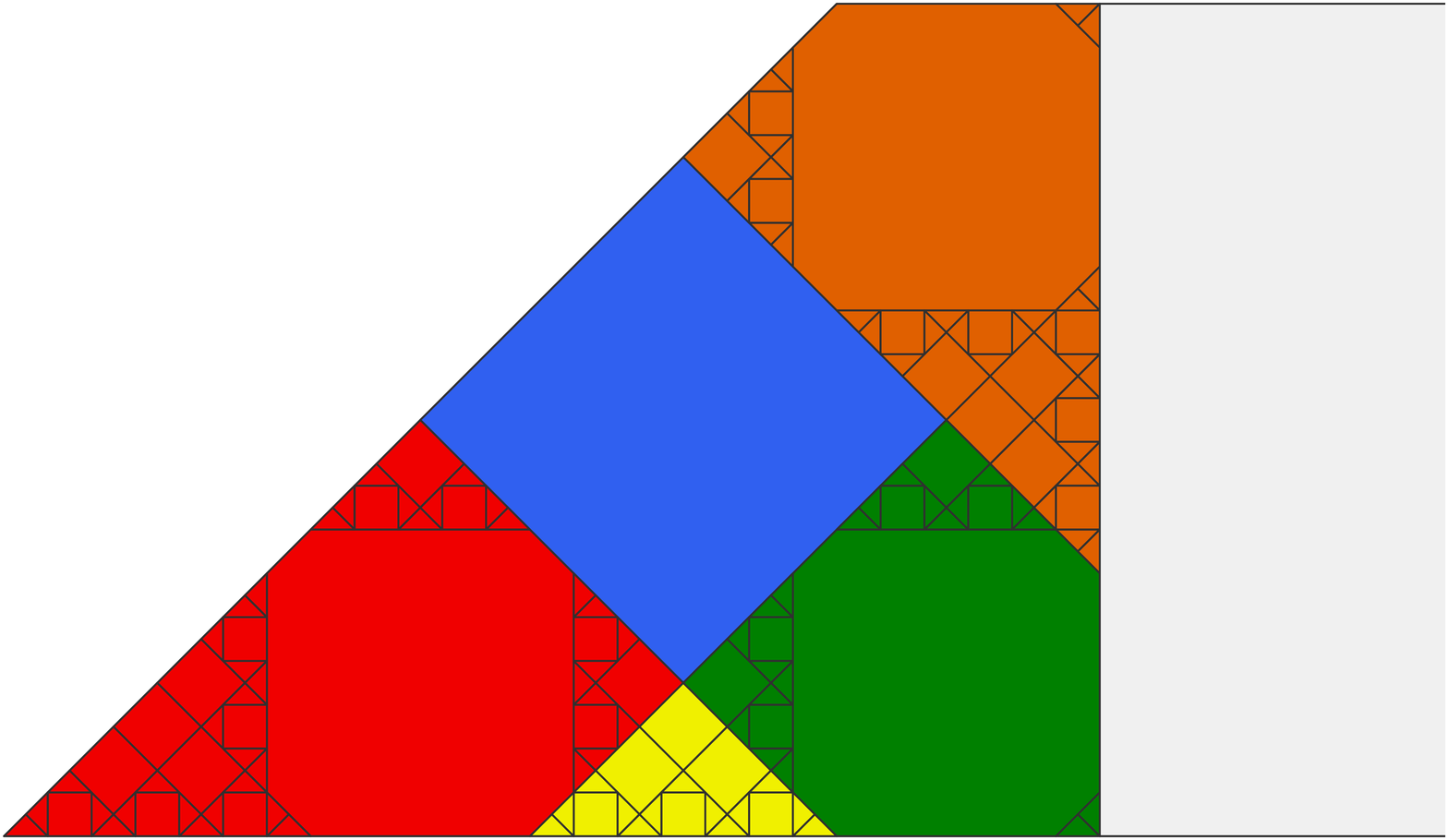}}
\newline
{\bf Figure 7.3\/} $Z_s$ (red) and
$A_s$ (green, blue, orange) and $B_s$ (red, yellow)
\end{center}

For $s \in (1/3,2/5)$, the set $Z_s$ covers the left
half of the triangle $B_s$ from the Bilateral Lemma,
and also $Z_s \subset B_s$.  Let $\rho_1$ be reflection
in the vertical line of symmetry of $B_s$.  By the
Bilateral Lemma, the set $\rho_1(Z_s)$ fits nicely
over $\Delta_s$.  The two pieces
$Z_s$ and $\rho_s(Z_s)$ cover $B_s$.

Let $\rho_2$ be reflection in the line $D_s$ from the
Bilateral Lemma.  The pieces
$\rho_2(Z_s)$ and $\rho_2 \circ \rho_1(Z_s)$ fit nicely
over $\Delta_s$.  This is the orange region in Figure 7.3.

Aside from the blue tile $\tau$, the only region not
yet covered is the green region, and this is covered
by $\rho_3 \circ \rho_2(Z_s)$, where $\rho_2$ is
reflection in the $x$-axis, $H$.  This last-mentioned piece
also fits nicely over $\Delta_s$.

To finish the proof, we note again that $Z_s=\phi_s(X_t^0)$.
The point, again, is that there is just one central tile
in $\Delta_t$.
So, all the sets we are using have the form
$I \circ \rho_s(X_t^0)$, where $I$ is an isometry.
\endproof

\begin{corollary}
The Covering Lemma holds for $s \in (5/4,3/2)$.
Moreover, none of the similar copies
of $X_t^0$ crosses the bottom edge of $X_s^0$.
\end{corollary}

\startproof
This follows from the Inversion Lemma, and from
the fact that the covering we constructed
in Lemma \ref{c1} is such that no similar
copy crosses the left edge of $X_s^0$.
\endproof

\begin{lemma}
\label{induct0}
The Covering Lemma holds for $s \in (1,3/2)$.
\end{lemma}

\startproof
We give the same kind of inductive
proof we used in the proof of
Lemma \ref{main1}.
We already know the result holds for
$s \in (5/4,3/2)$. 
 Given
$s \in (1,5/4)$ we let
$u=T(s)$, where $T$ is the map in Equation
\ref{modular}.  We will show that the
truth of the lemma for $u$ implies
the truth of the lemma for $s$.  In
all cases, we will use the fact that
no tile in the covering crosses the
bottom edge of the relevant region.

Let $t=R(s)$ and $v=R(u)$.
Suppose that $X_u^0$ has
a nice covering by similar copies of
$X_v^0$, with the additional property that
no piece crosses the bottom edge of
$X_u^0$.  Using the map
$\omega_s: Y_u \to W_s$ from the
proof of Lemma \ref{main1}, we get
a nice covering of $W_s$ by finitely
many tiles of $\Delta_s$ and finitely
many similar copies of $X_t^0$.
We also have the piece $Z_s$.

\begin{center}
\resizebox{!}{1.8in}{\includegraphics{Pix/mod1.ps}}
\newline
{\bf Figure 7.4:\/} 
Half of $W_s$ for $s=22/19$.
\end{center}

For convenience, we repeat Figure 5.1 here.
The union $Z_s \cup W_s$ covers
the entire left half of $X_s^0$,
namely the half which lies to the
left of the vertical line through
the top vertex of $X_s^0$. Moreover,
as in the proof of Lemma \ref{main1},
we note that reflection in this vertical
line maps the left right half of the
tiling into the left half.  Thus, we
can extend our covering to all of
$X_s^0$ using this reflection. Some
of the pieces slop over into the central
tile of $\Delta_s$, but this does not
bother us.

As in the proof of Lemma \ref{main1},
this lemma now follows from induction
on the number $k$ needed so that
$T^k(s) \in (5/4,3/2)$.
\endproof

\begin{corollary}
The Covering Lemma holds for $s \in (0,1/2)$.
\end{corollary}

\startproof
Combining the previous results, we see
that the Covering Lemma holds for
all $s \in (1,2)$.  By the Inversion
Lemma, the Covering Lemma holds for
all $s \in (1/4,1/2)$.  But then,
by repeated applications of the
Insertion Lemma, the Covering Lemma
holds for all $s \in (0,1/2)$.
\endproof

\noindent
{\bf Remark:\/}
We have ignored the parameter $1/4$.
The Covering Lemma is vacuous here,
because $R(1/4)=0$.

\begin{lemma}
The Covering Lemma Holds for $s \in (1/2,3/4)$.
\end{lemma}

\startproof
This case is similar to the case treated
in Lemmas \ref{c0}.
Figure 7.5 shows the rather typical
case of $s=7/10$.

\begin{center}
\resizebox{!}{1.8in}{\includegraphics{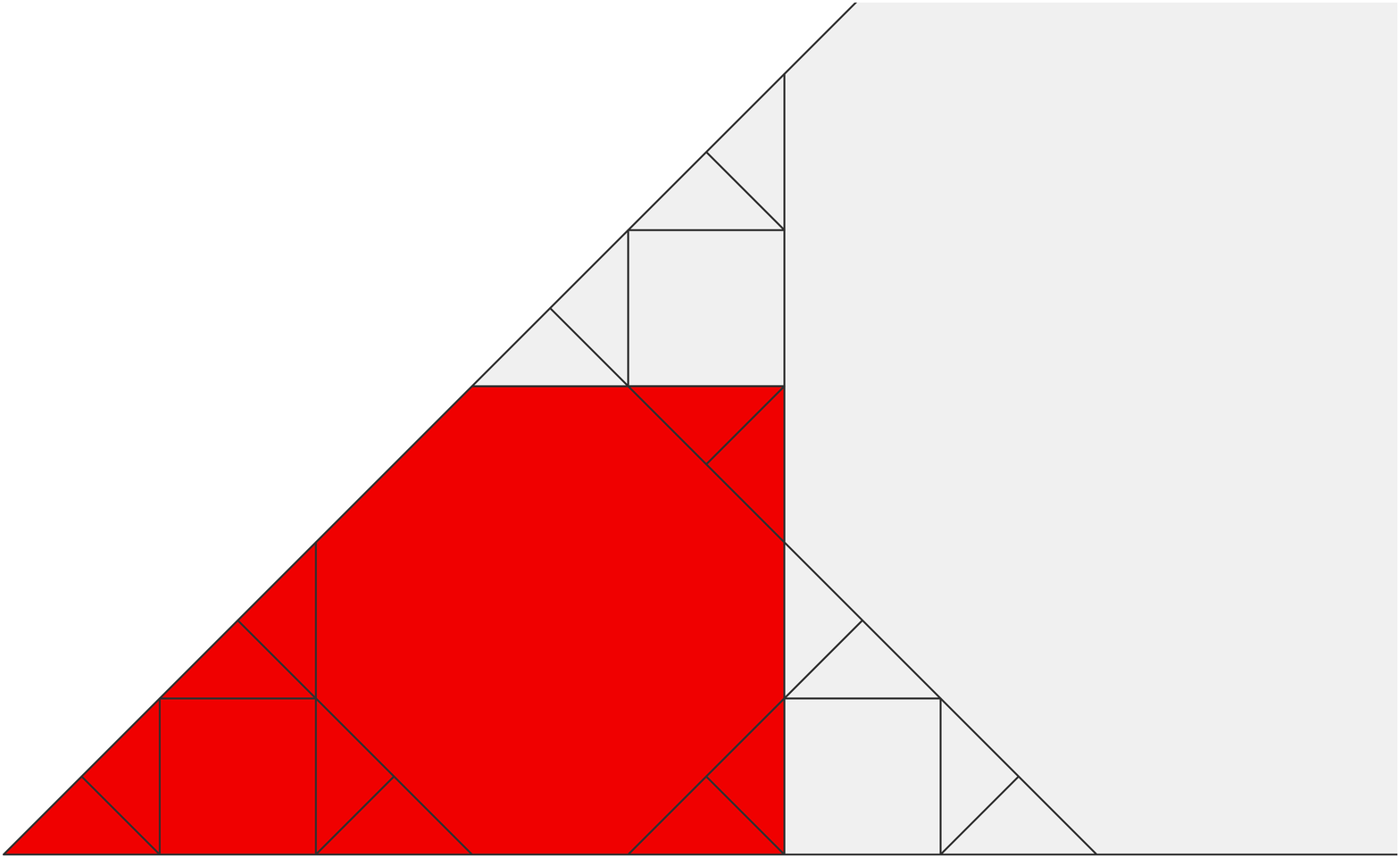}}
\newline
{\bf Figure 7.5:\/} 
$\Delta)s$ (white, red) and $Z_s$ (red) for $s=7/10$.
\end{center}

Using reflections in the lines $V_s$ and $D_s$ from the
Bilateral Lemma, we produce two isometric copies
of $Z_s$ which nicely fit over $\Delta_s$.
\endproof

\begin{lemma}
The Covering Lemma Holds for $s \in (3/4,5/6)$.
\end{lemma}

\startproof
The proof is similar to the case for
$s \in (1/2,3/4)$.  Figure 7.6 shows a
typical example.
\endproof

\begin{center}
\resizebox{!}{1.8in}{\includegraphics{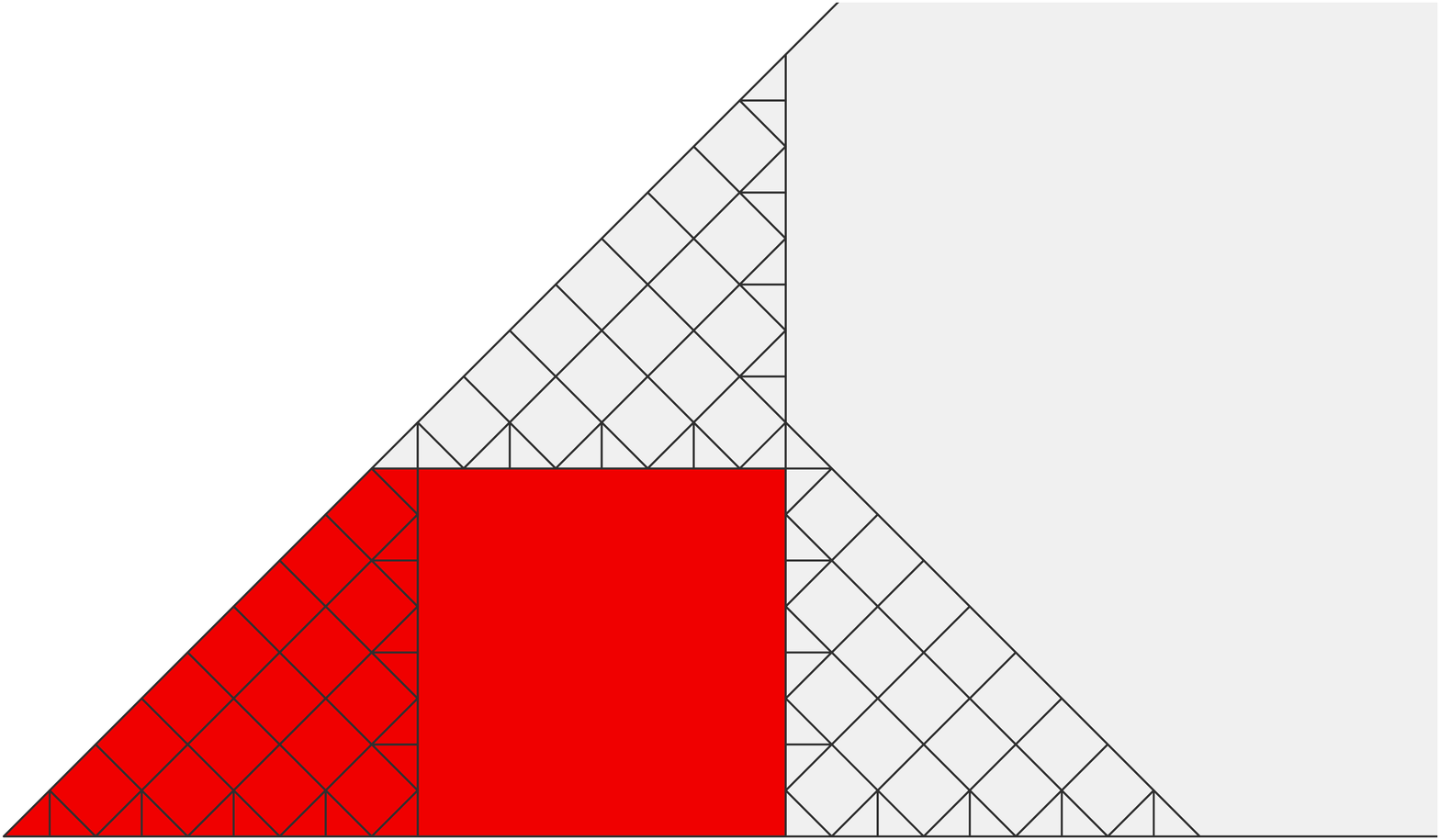}}
\newline
{\bf Figure 7.6:\/} 
$\Delta)s$ (white, red) and $Z_s$ (red) for $s=13/17$.
\end{center}

\begin{lemma}
The Covering Lemma Holds for $s \in (5/6,1)$.
\end{lemma}

\startproof
The proof here relates to the proof
of Lemma \ref{main2} in the same way that
the proof of Lemma \ref{induct0} relates
to the proof of Lemma \ref{main1}.
We omit the details.
\endproof

\newpage
\section{Properties of the Limit Set}

\subsection{Area Zero}

Here we show that $\widehat \Lambda_s$ has
zero area. It suffices to prove this for the left
half of $\widehat \Lambda_s$, which we denote by $S_s$.

\begin{lemma}
\label{cov1}
Let $s$ be irrational.
For any $\epsilon>0$, the
set $X_s^0$ has a nice covering by
similar copies of $X_t^0$, all
having diameter at most $\epsilon$.
\end{lemma}

\startproof
The nice fitting property is hereditary.  If
a similar copy $T_u(X_u^0)$ fits nicely over
$X_t^0$ and a similar copy
$T_t(X_t^0)$ fits nicely over $X_s^0$, then the
similar copy $T_t \circ T_u(X_u)$ fits nicely over
$X_s^0$. The boundary-adapted property is also
hereditary.
We let $s_n=R^n(s)$ and note that, by induction,
$X_s^0$ has a nice and boundary adapted covering
by similar copies of $X_t^0$ where $t=s_n$.
The scale factor of the similarities involved is the
product of the scales factors of $\phi_{s_i}$,
taken over $i=0,...,(n-1)$.  At least
every other time, $s_i<1/2$ and the corresponding
scale factor is less
than $\sqrt 2/2$, so the total scale factor
is less than $2^{-n/4}$. We can make this
less than $\epsilon/D$ by taking $n$ large enough.
Here $D<2$ is the maximum diameter of any set
$X_u^0$, taken over all $u \in (0,1)$.
\endproof

Say that a {\it near disk\/} is a compact set $D$
which is contained in a disk of radius $10r$ and
contains a disk of radius $r$.  For every $s$,
the set $X_s^0$ is a near disk.   Let
$\mu$ denote the $2$-dimensional Lebesgue measure.
The following is an immediate consequence of the
Lebesgue Density Theorem.

\begin{lemma}
\label{density}
Suppose that $S \subset \R^2$ is a bounded
measurable set of positive Lebesgue measure.
Then almost every point $p \in S$ has the
following property.  If $\{D_n\}$ is a
sequence of near disks containing $p$,
having diameter shrinking to $0$, then
$\mu(D_n \cap S)/\mu(S) \to 1,$
as $n \to \infty$.
\end{lemma}

Suppose now that $S_s$ has positive
measure.  We can find a point $p \in S_s$
satisfying the conclusion of Lemma \ref{density}.
By Corollary \ref{cov1} we can find a sequence
$\{D_n\}$ of near disks such that
$p \in D_n$ for all $n$, and
the diameter of $D_n$ tends to $0$, and
$D_n$ is a similar copy of $X_u^0$ for some $u=u_n$.
The Area Lemma now tells us that
$\mu(S_s \cap D_n)<\lambda (1-\epsilon)(D_n)$
for some universal constant $\epsilon>0$.
This contradicts Lemma \ref{density}.  
Hence $\Lambda_s$ has zero area.

\subsection{Characterization of the Limit Set}

Now that we know $\widehat \Lambda_s$ has
area zero, we can give a nicer characterization
of the limit set.

\begin{lemma}
\label{newdef}
A point belongs to $\widehat \Lambda_s$ if and only
if every open neighborhood of the point contains
infinitely many periodic tiles of $\Delta_s$.
\end{lemma}

\startproof
Certainly, $\widehat \Lambda_s$ contains all
such points. For the converse, suppose that
$p \in \widehat \Lambda_s$.
There exists a sequence of points
$q_n \to p$ with the following property
$f^1(q_n),...,f^n(q_n)$ are all defined
and distinct.  Since $\widehat \Lambda_s$
has measure $0$, and the set of points
with undefined orbits has measure $0$,
the set of periodic points has full measure.
So, we can take a new sequence
$\{q_n'\}$ of periodic points converging to $p$,
and we can make $|q_n-q_n'|$ as small as we
like. Making these distances sufficiently
small, we guarantee that $q_n'$ has
periodic at least $n$.  But then every
neighborhood of $p$ intersects infinitely
many periodic tiles.
\endproof

\subsection{Projections of the Limit Set}

Let $s \in (0,1)$ be irrational
and let $\pi$ denote projection
onto the $x$-axis.
We will show that $\pi(\widehat \Lambda_s)$
contains a line segment.  Following this,
we will deal with projections onto lines
parallel to the other $8$th roots
of unity.
We suppress the parameter $s$.
Since $\widehat \Lambda$ is a closed set,
it suffices to prove that $\pi(\Lambda)$
contains a dense subset of a line segment.

There are countably many vertical lines
which contain vertices of tiles in $\Delta$.
We ignore such lines.
Let $L$ be a vertical line which contains both
a point on the bottom edge $b$ of $X$ and a
point on the left edge $\ell$ of $X$.
Assume that $L$ does
not contain a tile vertex.
We will prove that $L$ contains a point of
$\widehat \Lambda$.  From this it follows
that $\pi(\ell) \subset \pi(\widehat \Lambda)$,
which proves what we want.

Suppose $L$ does not
intersect $\widehat \Lambda$.  Then
$L$ only intersects finitely many tiles,
$\tau_1,...,\tau_n$.  We order these
tiles according to when $L$ enters them
as we move upwards along $L$.
Note that $L$ must leave $\tau_i$ and enter
$\tau_{i+1}$ at the same point.  Otherwise,
the segment of $L$ lying between these
two tiles would be the accumulation point of
infinitely many tiles.  For the same reason,
the bottom edge of $\tau_1$ must lie in $b$
and the top edge
of $\tau_n$ must lie in $\ell$.

Since $L$ leaves $\tau_i$ and enters $\tau_{i+1}$
at the same point, $\tau_i$ and $\tau_{i+1}$ must
share an edge, and this shared edge contains
a point of $L$. In particular, if
$L$ leaves $\tau_i$ through a horizontal
edge, then $L$ enters $\tau_{i+1}$ through
a horizontal edge.

Since $L$ is a vertical line, and the tiles of
$\Delta$ are semi-regular octagons with sides
parallel to the $8$th roots of unity, we get
the following result: $L$ enters $\tau_i$ through
a horizontal edge if and only if $L$ leaves
$\tau_i$ through a horizontal edge.

We know that $L$ enters $\tau_1$ through a
horizontal edge. Using the properties above,
we see that $L$ leaves $\tau_n$ through a
horizontal edge. But the top edge of $\tau_n$,
which is contained in $\ell$,
is not horizontal. This is a contradiction.
Hence $L$ does intersect $\widehat \Lambda$.
\newline
\newline
{\bf Remark:\/} 
When $s$ is rational, our argument breaks down because
of the existence of triangular tiles. \newline

Now we deal with projection into a line parallel
to a different $8$th root of unity.  Let
$\omega$ be an $8$th root of unity.
Unless $\omega=\pm i$, we can find lines
perpendicular to $\omega$ which intersect
both a horizontal and a diagonal edge of
$X$.  Once we have such lines, the argument
we gave for $\omega=\pm 1$ works in this
new context.

We just have to worry about the case $\omega=\pm i$.
In this case, we observe that some edge of the
trivial tile in $\Delta$ is vertical.  Thus,
there are horizontal lines connecting $\ell$ to
this vertical edge.  Now we run the same
argument again, interchanging the roles
played by the horizontal and vertical directions.

\subsection{Finite Unions of Lines}

Let $s$ be irrational.
Let $S_s$ denote the left half of $\widehat \Lambda_s$.
We will assume that there exists an irrational
parameter $s$ such that $S_s$ is contained in a
finite union of lines, and we will derive a
contradiction.   

Say that an {\it essential cover\/} of $S_s$ is
a finite union of lines which covers all but
finitely many points of $S_s$.  Let $c(s)$ denote
the cardinality of an essential cover of
$S_s$ having the fewest number of lines.
We call $s$ a {\it minimal failure\/} if
$c(s)$ achieves the minimum possible value of the
function \footnote{We define $c(t)=\infty$
is $S_t$ has no essential cover.} $c: (0,1) \to \N$.
We call a cover realizing $c(s)$ a {\it minimal
essential cover\/}.  For the rest of
our proof we assume that $s$ is a minimal failure.

Referring to the Main Theorem, we can pull
back essential covers by $\phi_s$.
Suppose that $L$ is a line 
intersecting $Z_s^0$ in a
segment. Then we define
$\phi_s^{-1}(L)$ to be the line extending the
segment $\phi_s^{-1}(L \cap Z_s^0)$. 
If $L$ does not intersect $Z_s^0$ in a
line segment we define
$\phi_s^{-1}(L)$ to be the empty set.

Let $t=R(s)$.   By the Main Theorem,
\begin{equation}
\phi_s(S_t) \subset S_s \cap Z_s^0.
\end{equation}
for this reason, $t$ is also
a minimal failure, and the pullback of a
minimal essential
cover of $S_s$ is a
minimal essential cover of $S_t$.

\begin{lemma}
Every line of a minimal essential cover of $S_s$.
contains $v_s$, the bottom left vertex of $X_s$.
\end{lemma}

\startproof
 Let $d(s,L)$ denote the
distance from $v_s$ to a line $L$ of an
essential minimal cover.  Let $d(s)$
denote the maximum of $d(s,L)$, taken over all lines
of $L$ which appear in some essential minimal
cover of $S_s$.
Suppose that $d(s)>0$.
Since $\phi_s$ is a contraction, we
have $d(s)<d(t)$.  Also,
setting $u=R(t)$, we have
\begin{equation}
\label{cx}
d(s)<d(u)/\sqrt 2.
\end{equation}
This is true because at least one of the
maps $\phi_s$ or $\phi_t$ contracts
distances by at least a factor of
$\sqrt 2$.

Equation \ref{cx} is not possible for all
choices of $s$. Note that $X_s$ has diameter
at most $1+\sqrt 2$.  So, $d(s)$ is uniformly
bounded.  Hence, we can choose $s$ so as to
maximize $d(s)$, amongst all minimal
counterexamples, up to a factor of (say)
$9/10$. But then $d(u)$ exceeds the
maximum.  This is a contradiction. The
only way out is that $d(s)=0$.
\endproof

\begin{lemma}
Let $A_s$ denote the hexagon from the
Bilateral Lemma.
Only finitely
many points of $S_s$ lie in the $A_s-\partial X_s$.
\end{lemma}

\startproof
Let $S_s'$ denote the subset of $S_s$ contained
in the interior of $A_s$.
Let $U_s$ denote a minimal essential cover. Let
$\rho$ denote reflection in the $x$-axis, as
in the Bilateral Lemma.
$\rho(S_s')=S_s'$.  But then
$S_s' \subset U_s \cap \rho(U_s).$
This is just a finite set of points, because
all lines in $U_s$ contain the vertex $v_s$
and none of the lines in $\rho(U_s)$ contains $v_s$.
Hence $S_s'$ is a finite set.

Consider points of $S_s$ lying
on the edge of $A_s$ that does not lie in
$\partial X_s$. If this edge contained infinitely
many points of $S_s$ then $U_s$ would have
infinitely many lines. The point is that this
edge of $A_s$ does not contain the vertex $v_s$.
\endproof

\begin{lemma}
Let $B_s$ denote the triangle
 from the
Bilateral Lemma.
Only finitely
many points of $S_s$ lie in the $B_s-\partial X_s$.
\end{lemma}

\startproof
Same argument as above.
\endproof

Now we know that the interior of
$X_s$ contains only finitely many
points of $S_s$.

Let $\ell_+$ (respectively $\ell_-$)
 denote the subset of the
interior of $\ell$ lying above
(respectively below) $(-1,0)$,
a common vertex of $A_s$ and $B_s$.  Let
$\rho$ be reflection in the $x$-axis.
By the Bilateral Lemma,
\begin{equation}
\rho(S_s \cap \ell_+) \subset S_s \cap {\rm interior\/}(X_s).
\end{equation}  Hence $S_s \cap \ell_+$ is
finite.
A similar argument shows that
$S_s \cap \ell_-$ is finite.
Hence $S_s \cap \ell$ is finite.

But then $S_s$ is contained in the top
and bottom (horizontal) edges of $X_s$,
except for finitely many points.
But then the projection of $S_s$ onto the
$y$-axis does not contain a segment. This
contradicts Statement 2 of Theorem \ref{two}.
Hence a minimal failure cannot exist.
This proves
Statement 3 of Theorem \ref{two}.

\subsection{Intersection with the Bad Set}
\label{intersectX}

Here we will show that $\widehat \Lambda_s-\Lambda_s$
has length $0$ for almost every parameter $s$.
What we will actually prove is that
$\widehat \Lambda_s-\Lambda_s$ has length zero
provided that the sequence
$$\bigg\{\frac{1}{2R^n(s)}\hskip 10 pt {\rm mod\/} \ 1\bigg\}$$ has an
accumulation point in $(0,1)$.
Almost every choice of $s$ satisfies this criterion.
Suppose now that $s \in (0,1)$ is a parameter
which satisfies the condition above.

\begin{lemma}
\label{cov3}
Let $s$ be irrational.
For any $\epsilon>0$, some
open neighborhood of the special
segments of $X^0$ has a nice covering
by open subsets of similar copies of
$X_t^0$, all
having diameter at most $\epsilon$.
Moreover, $t$ can be chosen so that
$1/(2t)$ mod $1$ is uniformly
bounded away from $0$ and $1$.
\end{lemma}

\startproof
We simply restrict the covers produced in the
Covering Lemma to neighborhoods of the
special segments.  This gives us everything
in the lemma except the last statement.
The last statement comes from the hypotheses
on $s$.
\endproof

\begin{lemma}
\label{boundary}
$\widehat \Lambda_s \cap \partial X_s$ has
length $0$.
\end{lemma}

\startproof
Referring to the Length Lemma from the
previous chapter, it suffices to show that
$\Lambda_s$ intersects the special segments in sets of
measure $0$.
The proof of this lemma is exactly
like what we did for $S_s$, except that we
use the Lebesgue Density Theorem for length
instead of for area, and we use
Lemma \ref{cov3} in place of
Corollary \ref{cov1}.
\endproof

Let $B_s \subset X_s$ denote the set of points where
some iterate of $f_s$ is undefined.
Let $B_s^+$ (respectively $B_s^-$) denote the set of points where
some positive (respectively negative) 
iterate of $f$ is undefined.
We will show that $\widehat \Lambda \cap B_s^+$ has
length $0$.  The case of $\widehat \Lambda_s \cap B_s^-$
is quite similar.

\begin{lemma}
Let $s>1/2$.
If $\widehat \Lambda_s \cap B_s^+$ has positive
length, then $\Lambda_s \cap \partial X_s$ has
positive length.
\end{lemma}

\startproof
We suppress the parameter $s$.
Suppose that $\widehat \Lambda \cap B^+$ has
positive length.

Recall that $X'=F_1 \cup F_2$ and
$f': X' \to X'$ is such that $f=(f')^2$.
It is easier to work with $f'$ because
the discontinuity set is simpler.
Let $S'_n \subset X'$ denote the set of
points where one of the first $n$ iterates
of $f'$ is not defined.

We have
\begin{equation}
B_+ \subset \bigcup_{n=1}^{\infty}S_n'
\end{equation}
The set on the right is countable.  So,
if $\widehat \Lambda \cap B_+$ has positive
length, there is some $n$ such that
$\widehat \Lambda \cup S_n'$ has positive length.
We also have the equation
\begin{equation}
f'(S_n'-S_1') \subset S'_{n-1}.
\end{equation}
Using this equation, together with the fact that $f'$
is a piecewise isometry, we see by induction
that $\widehat \Lambda \cap S_1'$ has positive
length.

By symmetry, the $\pi/2$ rotation carries
$F_1$ to $F_2$ and $\widehat \Lambda \cap F_1$
to $\widehat \Lambda \cap F_2$ and
$S_1' \cap F_1$ to $S_1' \cap F_2$.  Hence,
if $\widehat \Lambda$ has positive length
intersection with
$S_1'$, then $\widehat \Lambda$ also has
positive length intersection with $S_1' \cap F_1$.
(Recall that $X=F_1$.)

\begin{center}
\resizebox{!}{1.5in}{\includegraphics{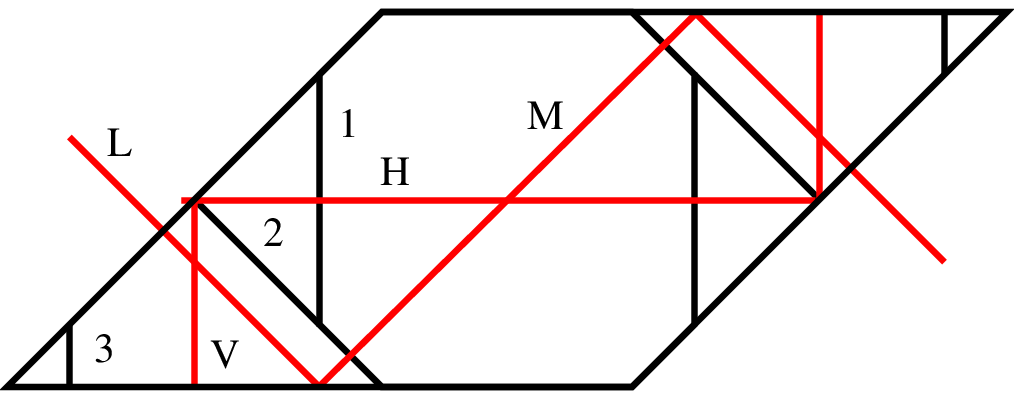}}
\newline
{\bf Figure 8.1:\/} The 
set $S_1' \cap X$ in black, and symmetry lines in red.
\end{center}

The red lines labelled $H,V,L$ in Figure 8.1 
are lines of symmetry from
Bilateral Lemma I and Bilateral Lemma II.
The fact that line $M$ is a line of symmetry
follows from the Inversion Lemma and the fact
that $L_t$ is a line of symmetry for $A_t$
when $t=1/2s$.

Let $R_L$ denote reflection in line $L$, etc.
We now have $3$ cases.
\begin{enumerate}
\item Suppose $\widehat \Lambda$ intersects the
segment labelled $1$ in a set of positive
length.  We apply $R_L$.
Lemma \ref{newdef} and the Bilateral Lemma II
combine to say that $\widehat \Lambda \cap \partial X$
has positive length.
\item Suppose $\widehat \Lambda$ intersects the
line labelled $2$ in a set of positive length.
On at least one side of line $2$ there are
tiles of $\Delta_s$ which accumulate on
line $2$ in a set of positive length.
Depending on which side has this property 
we either apply $R_H$ or $R_V$.
Lemma \ref{newdef} and Bilateral Lemma I
combine to say that $\widehat \Lambda \cap \partial X$
has positive length.
\item  Suppose $\widehat \Lambda$ intersects the
line labelled $3$ in a set of positive length.
This is the same argument as in Case 2, except
we use one of the reflections
$R_L \circ R_V$
or $R_M \circ R_V$. 
\end{enumerate}
If $\widehat \Lambda$ intersects $S_1'$ on the right
hand side of the picture in a set of positive length,
then we use the rotational symmetry to reduce to the
cases above.
\endproof

Combining the previous result with
Lemma \ref{boundary}, we see that
$\widehat \Lambda_s \cap B_s^+$ has
zero length when $s>1/2$. 
Now consider the case $t<1/2$.
Suppose that 
$\widehat \Lambda_t \cap B_t^+$ has positive
length.
Let $s=1-t$ so that $t=R(s)$.  By the Main Theorem,
$\widehat \Lambda_s \cap B_s^+$ has positive length.
But then, by the previous result,
$\widehat \Lambda_s \cap \partial X_s$ has
positive length.  Lemma \ref{boundary} rules
this out.

\subsection{Hyperbolic Symmetry}

In this section, we prove Theorem \ref{three}.
Let $i=\sqrt{-1}$, as usual.
Let $\Gamma'$ denote the group of maps
of $\C \cup \infty$ generated by the following maps.
\begin{enumerate}
\item $z \to \overline z$.
\item $z \to -z$.
\item $z \to z-1$.
\item $z \to 1/(2z)$.
\end{enumerate}
Let $\Gamma \subset \Gamma$ be the index
$2$ subgroup which preserves $\H^2$.

\begin{lemma}
If $s$ and $t$ lie in the same $\Gamma'$ orbit,
then $(X_s,f_s)$ and $(X_t,f_t)$ are locally
equivalent.
\end{lemma}

\startproof
Local equivalence is an equivalence relation, so we
just have to check this result on the generators of
$\Gamma$.  Note that local equivalence is an empty
condition for rational parameters.  So, we only
work with irrational parameters.

Complex conjugation fixes $\R$ pointwise. So, for
this generator there is nothing to prove.

If $t=-s$ then the two systems are identical, by definition.

If $t=1/2s$, then the two systems are conjugate, by the
Inversion Lemma.

The one nontrivial case is when $t=s-1$.  There are several
cases to consider.  By symmetry, it suffices to consider
the case when $s>0$.  There are several cases to consider.
When $s>2$, the result follows from the Insertion Lemma.

Suppose that $s \in (1,2)$.  By the Covering Lemma,
$X_s^0$ is nicely covered by similar copies of
$X_t^0$.  Moreover, in our construction, at least one
of these similar copies, namely $Z_s$, lies entirely
in $X_s^0$.  The local equivalence follows immediately
in this case.  To make this work, we need to throw out
points on the finitely many lines extending the sides
of the similar copies of $X_t^0$.

Suppose that $s \in (1/2,1)$ and $t'=s-1$.  Then $t'<0$ and
we can switch to the parameter $t=-t'=1-s$.  Now we
apply the Covering Lemma again, just as in the previous case.

Finally, suppose that $s \in (0,1/2)$.  Again, we 
consider $t=1-s \in (1/2,1)$.  Switching the roles of
$s$ and $t$ we reduce to the previous case.
\endproof

It just remains to recognize $\Gamma$.
\begin{itemize}
\item $\Gamma$ contains the element
$\rho_1(z)=-\overline z$. This is the reflection
in the geodesic joining $0$ to $\infty$.
\item $\Gamma$ contains the element
$\rho_2(z)=\overline{1-z}$. This is the reflection
in the geodesic ray joining $1/2+i/2$ to $\infty$.
\item $\Gamma$ contains the element
$\rho_3(z)= 1/(2\overline z)$.  
This is a reflection in the geodesic
joining $i/\sqrt 2$ to $1/2 + i/2$.
\end{itemize}
In short, $\Gamma$ contains
$\rho_j$ for $j=1,2,3$, and these elements generate
the $(2,4,\infty$ triangle group generated by
reflections in the sides of the triangle with
vertices $i/\sqrt 2$, and $1/2+i/2$, and
$\infty$.  We omit the
proof that these elements generate $\Gamma$, since
what we have done already gives a proof of
Theorem \ref{three}.

\newpage

\section{The Eight Calculations}
\label{calc}

\subsection{The Fiber Bundle Picture}

\label{fiber}

Our strategy is to reduce all our calculatons
to statements about convex lattice polyhedra.
In this section, we explain the main idea
behing this reduction.
For the most part, we only need to make
our calculations for parameters
$a \in [1/4,2]$, and so this interval will
make a frequent appearance in the discussion.

In our standard normalization,
$X_s$ be the parallelogram with vertices
\begin{equation}
\label{standard}
(\epsilon_1 + \epsilon_2 s,\epsilon_2 s), \hskip 30 pt
\epsilon_1,\epsilon_2 \in \{-1,1\}.
\end{equation}
We define
\begin{equation}
{\cal X\/}=\{(x,y,s)|\ (s,y) \in X_s)\}
\subset \R^2 \times [1/4,2]
\end{equation}
The space ${\cal X\/}$ is both a convex
lattice polyhedron and a fiber bundle over
$[1/4,2]$ such that
the fiber above $s$ is the parallelogram $X_s$.
See \S \ref{somepoly} for the list of vertices of
$\cal X$.
The maps $f_s: X_s \to X_s$ piece together
to give a fiber-preserving
map $F: {\cal X\/} \to {\cal X\/}$.

\begin{lemma}
\label{pa}
The map $F$ is a piecewise affine map.
\end{lemma}

\startproof
Consider some point
$p \in X_s$.  There is some vector
$$V_s=(A + Bs,C+Ds)$$
Here $A,B,C,D$ are integers.
such that $f_s(p)=p+V_s$.  If we perturb both $V$ and
$s$, the integers $A,B,C,D$ do not change.
So, in a neighborhood of $p$ in $\cal X$, the map $F$
has the form
$$(F(x,y,s)=(x+Bs,y+Ds,s)+(A,C,0).$$
This is a locally affine map of $\R^3$.
\endproof

Let ${\cal X\/}(S)$ be the subset of
our fiber bundle lying over a set
$S \subset [1/4,2]$.  Given the nature of
the map, we find it useful to split
our fiber bundle into $3$ pieces, namely
\begin{equation}
{\cal X\/}=
{\cal X\/}[1/4,1/2] \cup
{\cal X\/}[1/2,1] \cup
{\cal X\/}[1,2].
\end{equation}

A {\it maximal domain\/} of ${\cal X\/}(S)$ is a
maximal subset on which $F$ is entirely defined
and continuous. The map $F$ acts as an affine map
on each maximal subset.  In \S \ref{verify}
we explain how we verify the following
experimentally discovered facts.

\begin{itemize}

\item  ${\cal X\/}[1/4,1/2]$ is partitioned into $19$ maximal
domains, each of which is a convex rational polytope.
The vertices in the partition of
 ${\cal X\/}$ are of the form
$(a/q,b/q,p/q)$ where
$a,b,p,q \in \Z$ and $$(p,q) \in \{(1,4),(2,7),(3,10),(1,3),(1,2)\}.$$
These polyhedra are permuted by the map
$\iota_1(x,y,s)=(-x,-y,s)$.

\item ${\cal X\/}[1/2,1]$ is partitioned into $13$ maximal
domains, each of which is a convex rational polytope.
 The vertices in the partition of
 ${\cal X\/}[1/2,1]$ are of the form
$(a/q,b/q,p/q)$ where
$a,b,p,q \in \Z$ and $$(p,q) \in \{(1,2),(2,3),(3,4),(1,1)\}.$$
These polyhedra are permuted by $\iota_1$.

\item The $19$ polyhedra in the partition of
${\cal X\/}[1,2]$ are all images of the polyhedra in the
partition of ${\cal X\/}[1/4,1/2]$ under the map
$\iota_2 \circ F$, where
$$\iota_2(x,y,s) \to ((x+y)/2s,(x-y)/2s,1/2s).$$
These polyhedra are permuted by $\iota_1$.

\end{itemize}

In \S \ref{move} we will explain the action of
the map $F$ on each polyhedron.  We scale
all the polyhedra by
a factor of $420=10 \times 7 \times 4$ so that we
can make all our calculations using integer arithmetic,
One can think of this rescaling as
a way of clearing all denominators in advance
of the calculations.
  
We will list the $420$-scaled
polyhedra in \S \ref{partition}.
For now, call them 
$$\alpha_0,...,\alpha_{18},\beta_0,...,\beta_{12},
\gamma_0,...,\gamma_{18}.$$  We have labeled so that
list these so that 
\begin{equation}
\label{contain1}
\alpha_i \subset 
{\cal X\/}[105,210], \hskip 10 pt
\beta_i \subset
{\cal X\/}[210,420], \hskip 10 pt
\gamma_i \subset {\cal X\/}[420,840].
\end{equation}
${\cal X\/}[105,210]$ is our name
for the $420$-scaled version of
${\cal X\/}[1/4,1/2]$, etc.
When we want to discuss these polyhedra
all at once, we call them $P_0,...,P_{50}$.

\subsection{Computational Methods}

Here we describe the main features of
our calculations. As we mentioned above,
all our calculations boil down to calculations
involving convex lattice polytopes.
\newline
\newline
{\bf Operations on Vectors:\/}
The only operations we perform on
vectors are vector addition
and subtraction, the dot and cross product,
scaling a vector by an integer, and
dividing a vector by $d \in \Z$ provided 
$d$ divides all the coordinates.
These operations in turn only use
plus, minus, and times, and
integer division.
\newline
\newline
{\bf Avoiding Computer Error:\/}
We represent itnegers as {\it longs\/}, a
$64$ bit integer data type.  One
can represent any integer strictly  between
$-2^{-63}$ and $2^{63}$. (The extra bit
gives the sign of the number.)  There
are two possible sources of error:
overflow error and division errors.

We subject all our calculations to an
overflow checker, to make sure that the
computer never attempts a basic operation
(plus, minus, times) in with either the
inputs or the output is out of range.
To give an example, if we want to
take the cross product $V_1 \times V_2$,
we first check that all entries in
$V_1$ and $V_2$ are less than $2^{30}$
in size.  This guarantees that all
intermediate answers, as well as
the final answer, will be in the
legal range for longs. 

We also check our division operations.
Before we compute $n/d$, we make
sure that $n \equiv 0$ mod $d$.
The java operation $n\%d$ does this.
Once we know that $n\%d=0$, we know
that the computer correctly
computes the integer $n/d$.
\newline
\newline
{\bf No Collinearities:\/}
Given a polyhedron $P$, let
$P_1,...,P_n$ denote the vertices.
We first check that
\begin{equation}
\label{normal}
(P_k-P_i) \times (P_j-P_i) \not = 0,
\hskip 30 pt \forall\ i<j<k \in \{1,...,n\}.
\end{equation}
This guarantees that no three points
in our vertex list of $P$ are collinear.
We found the polyhedra of interest to us in
an experimental way, and initially they had
many such collinearities.  We
detected collinearities by the failure
of Equation \ref{normal}, and then
removed all the redundant points.
\newline
\newline
{\bf Face Lists:\/}
For each of our polyhedra $P$, we find
and then store the list of faces of the
polyhedron.  To do this, we consider
each subset $S=\{S_1,...,S_m\}$ having at least $3$ members.
We check for three things.
\begin{enumerate}
\item $S$ lies in a single plane.  We
compute a normal $N=(S_2-S_1) \times (S_3-S_1)$
and then check that $N \cdot S_i$ is independent of $i$.
Assuming this holds, let $D=N \cdot S_i$.
\item $S$ lies on $\partial P$.  To check this,
we compute the normal $N$ as above, and then
check that either $\max N \cdot P_i \leq D$
or $\min N \cdot P_i \geq D$.
\item We check that $S$ is maximal with respect
to sets satisfying the first two properties.
\end{enumerate}
\noindent
{\bf Improved Normals:\/}
We noticed computationally that all of
the normals to all of the polyhedron
faces can be scaled so that they
have the following form:
At least two of the three
coordinates lie in $\{-1,0,1\}$
and the third coordinate lies
in $\{-8,...,8\}$.  When we
use the normals in practice, we
make this scaling.  This is one
more safeguard against overflow
error.
\newline
\newline
\noindent
{\bf No Face Redundancies:\/}
Once we have the face list,
we check that each vertex of each
polyhedron lies in exactly $3$ faces.
In particular, all the vertices of
our polyhedra are genuine vertices.
\newline
\newline
{\bf Containment Algorithm:\/}
Suppose we want to check if a vector
$V$ lies in a polyhedron $P$.  For
each face $S$ of $P$, we
let $N$ be the (scaled) normal to $S$.
we set $D=N \cdot S_0$, and 
then we verify the following.
\begin{itemize}
\item If $\max N \cdot P_i \leq D$ then
$V \cdot N \leq D$.
\item If $\min N \cdot P_i \geq D$ then
$V \cdot N \geq D$.
\end{itemize}
If this always holds then $V$ lies on the
same side of $S$ as does $P$, for all
faces $S$.  In this case, we know
that $V \in P$. If we want to check that
$V \in {\rm interior\/}(P)$, we make the
same tests, except that we require strict
inequalities.
\newline
\newline
{\bf Disjointness Algorithm:\/}
Let $\Z_{10}=\{-10,...,10\}$.
To prove that two polyhedra
$P$ and $Q$ have disjoint interiors, we produce
(after doing a search)
an integer vecto $W \in \Z_{10}^3$ such that
\begin{equation}
\max W \cdot P_i \leq \min W \cdot Q_j.
\end{equation}
\newline
\newline
{\bf Volumes:\/} For many of the lattice polytopes
$P$ we consider, we compute $6\ {\rm volume\/}(P) \in \Z$.
To compute this volume, we decompose $P$ into prisms
by choosing a vertex of $v \in P$ and then computing
\begin{equation}
\sum_{f \in P} {\rm 6 volume\/}([v,f]).
\end{equation}
The sum is taken over all faces $f$ of $P$ and
$[v,f]$ denotes the prism obtained by taking
the convex hull of $v \cup f$.  We compute the
volumes of these prisms by taking various triple
products of the vectors involved.

\subsection{Verifying the Partition}
\label{verify}

Now we explain how we verify that the
polyhedra we mention in \S \ref{fiber}
(and list in \S \ref{partition}) really
are correct.
Let $\cal{RX}$ denote the polytope obtained from
$\cal X$ by rotating $90$ degrees.
We have $F=(F')^2$, where
$F'$ is the original PET which swaps 
$\cal X$ and $\cal{RX}$.
\newline
\newline
{\bf Pairwise Disjointness:\/}
Using the Disjointness Algorithm, we check
that $P_i$ and $P_j$ have disjoint interiors for
all $i \not = j \in \{0,...,50\}$.   We check
the same thing for
$F(P_i)$ and $F(P_j)$.
\newline
\newline
{\bf Containment:\/}
Using the Containment Algorithm, we check that
\begin{itemize}
\item $P_i \subset {\cal X\/}$ for $i=0,...,50$.
\item $F(P_i) \subset {\cal X\/}$ for $i=0,...,50$.
\item $F'(P_i) \subset {\cal RX\/}$ for $i=0,...,50$.
\end{itemize}
We also see, by inspection, that $F$ has a different
action on $P_i$ and $P_j$ whenever $P_i$ and
$P_j$ share a (non-horizontal) face.
These checks show that each $P_i$ is a maximal
domain for the action of $F$ 
\newline
\newline
{\bf Filling:\/}
It remains to check 
we check that ${\cal X\/}$ is partitioned
into $P_0,...,P_{50}$.  We check that
\begin{equation}
\label{volcheck}
\sum_{i=0}^{50} {\rm volume\/}(P_i)=
{\rm volume\/}({\cal X\/}).
\end{equation}
The same equation shows that
${\cal X\/}$ is also partitioned
into $F(P_0),...,F(P_{50})$.

\subsection{Calculation 1}

  Let $H=F^{-1}$,
and let $\cal H$ be the partition of
${\cal X\/}[1/4,1]$ by the
polyhedra $F(P_0),...,F(P_{31})$.
Then $\cal H$ is the partition by
maximal domains for $F^{-1}$.
We rename the members of 
$\cal H$ as $H_0,...,H_{31}$.

We have the partition
\begin{equation}
{\cal X\/}[1/4,1]=
{\cal A\/}[1/4,1] \cup
{\cal B\/}[1/4,1] \cup
{\cal C\/}[1/4,1]
\end{equation}
Here $\cal A$ is such that the fiber of
$\cal A$ over $s$ is the hexagon $A_s$.
The polyhedron $\cal B$ has the same
definition relative to the triangle $\beta_s$.
The polyhedron $\cal C$ is obtained from
$\cal B$ by reflecting in the line $x=y=0$.

The map $\mu$ from Calculation 1 acts on each
of $\cal A$ and $\cal B$ and $\cal C$ as a
reflection.  We verify that each polyhedron
$\alpha_i$ and $\beta_i$
is a subset of one of these $3$
big pieces.  Thus, $\mu$ acts on each polyhedron as
a reflection. The new partition
\begin{equation}
{\cal X\/}[1/4,1]=\bigcup_{i=0}^{18}\mu(\alpha_i) \cup
\bigcup_{i=0}^{12}\mu(\beta_i
\end{equation}
is the partition for the map 
\begin{equation}
G=\mu \circ f \circ \mu^{-1}.
\end{equation}
We call this partition $\cal G$,
and we rename its members $G_0,...,G_{31}$.

So, in summary $\cal G$ is the partition for
$G$ and $\cal H$ is the partition for $H$.
Next, we find a list of $48$ pairs $i,j$
so that
$${\rm interior\/}(G_i) \cap {\rm interior\/}(H_j) \not = \emptyset$$
only if $(i,j)$ lies on our list.  More
precisely, we use the Separation Algorithm
to show that all other pairs have disjoint interiors.

Finally, we consider the grid 
\begin{equation}
\Gamma=\{(20i,20k,105+10k)\ i=-42,...,42,\ j=-21,...,21,\ k=0,...,31\}.
\end{equation}
We check the identity
$G=H$ on each point of $\Gamma$ and we also check that at least
one point of $\Gamma$ is contained in each intersection
$G_i \cap H_j$ for each of our $48$ pairs.  This
suffices to establish the identity on all of
${\cal X\/}[1/4,1]$.

\subsection{Calculation 2}

Calculation 2 follows the same scheme as Calculation 1.
Here we just explain the differences in the calculation.
\begin{itemize}
\item We set $H=F$ and $G=\nu F^{-1} \nu$.
\item $\cal H$ is the partition
consisting of $\alpha_0,...,\alpha_{18}$.
\item ${\cal X\/}[1/4,1/2]$ is partitioned into
$5$ smaller polyhedra, coming from 
$P_s$, $Q_s$, the central tile,
$\iota(P_s)$ and $\iota(Q_s)$.  the map
$\nu$ acts as a reflection on each piece.
For $i=0,...,18$, the polyhedron $F(\alpha_i)$ is contained
in one of the $5$ pieces, so that $\nu$ acts
isometrically on $F(\alpha_i)$.
\item $\cal G$ be the partition of
${\cal X\/}[1/4,1/2]$ by the polyhedra
$$\nu \circ F(\alpha_0),...,\nu \circ F(\alpha_{18}).$$
\item We find a list of $27$ pairs
$(i,j)$ such that $G_i$ and $H_j$ do not
have disjoint interiors.
\end{itemize}
The rest of the calculation is the same.

\subsection{Calculation 3}

Let $s \in [5/4,2]$ and let $t=s-1 \in [1/4,1]$.
We want to show that $\phi_s$ conjugates
$f_t|Y_t$ to $f_s|Z_s$ and that
every orbit of $f_s$ intersects $Z_s$,
except the following orbits.
\begin{itemize}
\item Those in the trivial tile $(\alpha_0 \cup \beta_0)_s$ of $\Delta_s$.
\item Those in the set 
$$\tau_s=\phi_s\Big((\alpha_0 \cup \beta_0)_t\Big).$$
\end{itemize}
Once we are done, we will know that $\tau_s$ is
in fact a tile of $\Delta_s$, and that $\tau_s$
has period $2$.

For this section we set ${\cal X\/}={\cal X\/}[5/4,2]$.
Let $\cal Y$ denote the subset of $\cal X$ whose fiber
over $s$ is the set $Y_s$.  Define
$\cal Z$ in a similar way.
The maps $\phi_s: Y_t \to Z_s$ piece together to
give an isometry $\phi: {\cal Y\/} \to {\cal Z\/}$.
The map is given by
\begin{equation}
\phi(x,y,z)=(x \pm 1,y \pm 1,z-1)
\end{equation}
Whether we add or subtract $1$ to the first two
coordinates depends on whether the point
$(x,y,z)$ lies in the left half of $\cal Y$ or
in the right half.

For what we describe next, we always refer to
open polyhedra, and our equalities are meant
to hold up to sets of codimension $1$, namely
the boundaries of our polyhedra.

We have
\begin{equation}
{\cal Y\/}=\alpha_1 \cup ... \cup \alpha_{18} \cup
\beta_1 \cup ... \cup \beta_{12}.
\end{equation}

For each $i=1,...,18$ we
check computationally 
that there is some $k=k_i$ with the
following three properties.
\begin{enumerate}
\item The first $k_i+1$ iterates of $F^{-1}$ are defined on
$\phi \circ F(\alpha_i)$.  This amounts to checking that 
\begin{equation}
P_{ij}=F^{-i} \circ \phi \circ F(\alpha_i)
\end{equation}
is contained in some $\alpha_a$ or $\beta_b$ for suitable
indices $a$ and $b$, and for all $i=0,...,k_i$.
\item  
\begin{equation}
\label{fr}
P_{ij} \cap {\cal Z\/}=\emptyset, \hskip 30 pt
j=1,...,k_i.
\end{equation}
Equation \ref{fr} shows that
$$F|P_{i0}=F^{k_i}$$
That is, on $P_{i0}$, the map $F$ returns to
${\cal Z\/}$ as $F^{k_i}$.
To establish Equation \ref{fr},
we use the Separation Algorithm so show
that $$P_{ij}\cap \phi(\alpha_a)=\emptyset, \hskip 30 pt
P_{ij} \cap \phi(\beta_b)=\emptyset$$ 
for all
$a=1,...,18$ and $b=1,...,12$, and all relevant
indices $i$ and $j$.  This suffices because
$\cal Z$ is partitioned into the polyhedra
$$\cal Z=\phi(\alpha_1)\cup...\cup\phi(\alpha_{18})\cup\phi(\beta_1)\cup...\cup\phi(\beta_{12}).$$
\item $P_{i,k_i+1}=\phi(\alpha_i)$. 
\end{enumerate}
We make all the same calculations for $\beta_1,...,\beta_{12}$,
finding an integer $\ell_i$ which works for
$\beta_i$.  We define $Q_{ij}$ with respect to $\beta_i$
just as we defined $P_{ij}$ with respect to $\alpha_i$.

Our calculations above show 
that $\phi$ conjugates $F|{\cal Y\/}$ to $F|{\cal Z\/}$.
Also, by construction, the boundary of
${\cal Z\/}$ is contained in the union of the boundaries
of the polyhedra $\phi(\alpha_i) \cup \phi(\beta_j)$.  Hence,
$Z_s$ is a clean set for all $s \in [5/4,2]$.

We still want to see that all orbits except those of
period $1$ and $2$ actually intersect $\cal Z$.
We check the following.
\begin{enumerate}
\item $F$ is entirely defined on $\phi(\alpha_0)$ and has order $2$.
\item Both $\phi(\alpha_0)$ and $F \circ \phi(\alpha_0)$ are disjoint
from $\cal Z$.  We use the same trick with the Separation
Algorithm to do this.
\end{enumerate}

We claim that the open polyhedra in the following
union are pairwise disjoint.
\begin{equation}
\label{grand}
\bigcup_{i=1}^{18} \bigcup_{j=0}^{k_i}P_{ij} 
\cup \bigcup_{i=1}^{12} \bigcup_{i=0}^{\ell_i} Q_{ij}
\cup \bigcup_{j=0}^1 F^j \circ \phi(\alpha_0) 
\cup \bigcup_{j=0}^1 F^j \circ \phi(\beta_0) 
\cup \alpha_0 \cup \beta_0.
\end{equation}
Suppose, for instance, that
$P_{ab}$ and $P_{cd}$ were not disjoint.
Then $P_{a,b+e}$ and $P_{c,d+f}$ would not
be disjoint for $e>0$ and $f>0$ such that
$b+e=k_a+1$ and $c+f=k_b+1$.  But we know
that these last polyhedra are disjoint because
they respectively equal the disjoint
polyhedra $\phi(\alpha_a)$ and $\phi(\alpha_c)$.
Similar arguments work for the other cases.

Similar to Equation \ref{volcheck},
we compute the sum of the volumes
of the polyhedra in Equation \ref{grand} and see
that it coincides with the volume of 
${\cal X\/}$.  Thus,
${\cal X\/}$ is partitioned into
the polyhedra in Equation \ref{grand}.
This fact implies the all orbits except those of
period $1$ and $2$ actually intersect $\cal Z$.

Finally, we see by process of elimination that
$\tau_s$ really is a tile of $\Delta_s$.  All
other points not in the interior of
$\tau_s$ either have undefined orbits, or
lie in the trivial tile, or have orbits
which intersect $\cal Z$.  Thus $f_s$
cannot be defined on any point of the
boundary of $\tau_s$.  Since $f_s$ is
defined, and has period $2$, on the
interior of $\tau_s$, we see that
$\tau_s$ is a tile of $\Delta_s$ having
period $2$.

\subsection{Calculation 4}

Calculation 4 follows the same scheme as Calculation 3.
Here we will describe the differences between the two
calculations.

\begin{itemize}
\item We consider the behavior of polyhedra on the
interval $s \in [1/2,3/4]$ rather than on $[5/4,1]$.
Here $t=1-s \in [1/4,1/2]$.
\item The map $\phi$ is not an isometry here, but
rather a volume preserving affine map.
The formula is
\begin{equation}
\phi(x,y,z)=(x \pm (1-2z),y \pm (1-2z),1-z).
\end{equation}
The choice of plus or minus again depends on whether
$(x,y,z)$ lies in the left of the right half of
$\cal Y$.
\item $\cal Y$ is partitioned into the
tiles $\alpha_1,...,\alpha_{18}$.  The $B$-tiles are
not needed here.
\item The tiles $\tau$ and $\iota(\tau)$
already belong to $\cal Z$. The work in
Calculation 3 shows that $\tau_s$ and
$\iota_s$ are indeed
period $2$ tiles of $\Delta_s$.
This time, $\tau$ and $\iota(\tau)$ are
amongst the images of $\alpha_1,...,\alpha_{18}$
under $\phi$.
\item Using the notation from the previous section,
the partition in Equation \ref{grand} becomes
\begin{equation}
\label{grand2}
\bigcup_{i=1}^{18} \bigcup_{j=0}^{k_i}P_{ij} 
\cup \alpha_0
\end{equation}
\end{itemize}
The rest of the calculation is the same.

\subsection{Calculation 5}

Calculation 5 follows the same scheme as
Calculation 3, except that we don't need to keep
track of the volumes.  Let $T$ and
$\omega$ be as in Calculation 5.
Let $s \in (1,4/3]$ and let $t=T(s) \in (1,2]$.

We define $\cal W$ and $\cal Y$ as the
global versions of $W_s$ and $Y_u$, as
in Calculation 3. We are interested 
in ${\cal Y\/}[1,2]$ and ${\cal W\/}[1,4/3]$.
Similar to Calculation 3, we have a global map
$\omega: {\cal Y\/} \to {\cal W\/}$.
We have the formula
\begin{equation}
\label{omega}
\omega(x,y,z)=\big(\omega_s(x,y),s\big), \hskip 30 pt
s=T^{-1}(z).
\end{equation}
We want to see that $\omega$ conjugates
$F|{\cal Y\/}$ to $F|{\cal W\/}$.

We have
\begin{equation}
{\cal Y\/}=\gamma_1 \cup ... \cup \gamma_{18}.
\end{equation}
By the same methods used in Calculation 3, we check,
for each $i=1,...,18$,
that there is some $k=k_i$ with the
following three properties.
\begin{enumerate}
\item The first $k_i+1$ iterates of $F^{-1}$ are defined on
$\omega \circ F(\gamma_i)$.  Define
\begin{equation}
P_{ij}=F^{-i} \circ \omega \circ F(\gamma_i)
\end{equation}
\item  
\begin{equation}
\label{frX}
P_{ij} \cap {\cal W\/}=\emptyset, \hskip 30 pt
j=1,...,k_i.
\end{equation} 
\item $P_{i,k_i+1}=\omega(\alpha_i)$. 
\end{enumerate}
These facts imply
that $\omega$ conjugates $F|{\cal Y\/}$ to $F|{\cal W\/}$.

Finally, the set $Z_s$ is clean for each $s$ for the following
reasons.
\begin{itemize}
\item The top edge of $Z_s^0$ and the 
bottom edge of $\iota(Z_s^0)$ are contained in the union
of slices of the sets $\omega \circ F(\partial \gamma_i)$.  
\item The vertical edges of $Z_s$ are contained in the
set $\partial \tau_s \cup \iota(\partial \tau_s)$.
\item The remaining edges of $Z_s$ lie in the $\partial X_s$.
\end{itemize}

\subsection{Calculation 6}

As we mentioned in \S 5, Calculation 6 practically
amounts to inspecting the partition.  For Statement 1,
we let $\tau$ be the polyhedron which restricts to
$\tau_s$ for $s \in [1,3/2]$.  We list this polyhedron
in \S \ref{somepoly}.  We check that
$F$ is entirely defined on (the interior of) $\tau$
and that $F^2(\tau)=\tau$.

For each polyhedron $P$, let $P_s$ denote the
intersection of $P$ with the horizontal plane
of height $s$.

For Statement 2, let ${\cal Z\/}$ be the
polyhedron which restricts to $Z_s^0$ for 
$s \in [1,5/4]$.   We compute that
\begin{equation}
{\cal Z\/} \subset F(\gamma_{13}).
\end{equation}
We also try a single point $(x,y,s) \in F(\gamma_{13})$
and check that $f_s^{-1}(p)=p+\delta_s$.  Since
$F(\gamma_{13})$ is a domain of continuity for
$F^{-1}$, the same result holds for all points
in $F(\gamma_{13})$, including all the points in
$\cal Z$.  This proves Statement 2.

For Statement 3, let
$\cal K$ be the union of two polyhedra which intersect
the fiber $X_s$ in 
$X_s-Z_s-W_s$, for $s=(1,5/4]$. We see by inspection that
\begin{equation}
\label{inspect1}
{\cal K\/}=F(\gamma_2) \cup F(\gamma_8) \cup F(\gamma_{11}) \cup F(\gamma_{17}); \hskip 30 pt
\gamma_j \subset {\cal Z\/}, \hskip 20 pt j=2,8,11,17.
\end{equation}
This proves Statement 3.
\newline
\newline
{\bf Remark:\/}
We could have made an explicit computation to 
establish Equation \ref{inspect1}, but this
is something that is obvious from a glance
at just $2$ planar pictures.
We just have to check Equation
\ref{inspect1} at the parameters
$s=1$ and $s=3/2$ because every polyhedron $P$
in sight, when restricted to the fibers above
$[1,5/4]$, is the convex hull of
$P_1 \cup P_{5/4}$.

\subsection{Calculation 7}

Calculation 7 follows the same scheme as
Calculation 5.  Here are the differences.
\begin{itemize}
\item Here we are interested in
${\cal Y\/}[1/2,1]$ and
${\cal W\/}[3/4,1]$.
\item Here we use the formula from
\ref{modular2} to define tha map $\omega$
in Equation \ref{omega}.
\item Here we have
\begin{equation}
{\cal Y\/}=\beta_1 \cup ... \cup \beta_{12}.
\end{equation}
\end{itemize}
The rest of the calculation is the same.

\subsection{Calculation 8}

Calculation 8 works essentially the same was as
Calculation 6. 
  For Statement 1,
we let $\tau$ be the polyhedron which restricts to
$\tau_s$ for $s \in [3/4,1]$.  We list this polyhedron
in \S \ref{somepoly}.  We check that
$F$ is entirely defined on (the interior of) $\tau$
and that $F^2(\tau)=\tau$.

For Statement 2, we let ${\cal Z\/}^*$ be the
polyhedron which intersects $X_s$ in 
$(Z_s^0)^*$ for $s \in [3/4,1)$.  We
compute that
\begin{equation}
{\cal Z\/}^* \subset \beta_7
\end{equation}
and we finish the proof of Statement 2 just as in
Calcultion 6.

For Statement 3, we define $\cal K$ as
in Calculation 6 and we see by inspection that
\begin{equation}
{\cal K\/}=\beta_2 \cup \beta_6 \cup \beta_8 \cup \beta_{12},
\hskip 30 pt
F(\beta_j) \subset {\cal Z\/}, \hskip 20 pt
j=2,6,8,12.
\end{equation}
This proves Statement 3.

\newpage
\section{The Raw Data}

\subsection{Auxilliary Polyhedra}
\label{somepoly}

Here we list the coordnates of the
auxilliary polyhedra which arise in our
calculations.  All our polyhedra are scaled
by a factor of $420$.

Here is ${\cal X\/}[1/4,2]$.

{\footnotesize
\begin{eqnarray}
\nonumber
\left[\matrix{-525\cr -105\cr 105}\right]\left[\matrix{315\cr -105\cr 105}\right]\left[\matrix{-315\cr 105\cr 105}\right]\left[\matrix{525\cr 105\cr 105}\right]\left[\matrix{-1260\cr -840\cr 840}\right]\left[\matrix{-420\cr -840\cr 840}\right]\left[\matrix{420\cr 840\cr 840}\right]\left[\matrix{1260\cr 840\cr 840}\right]
\end{eqnarray}
\/}

Here is ${\cal A\/}[1/4,1]$.  This set intersects the fiber $X_s$ in
the hexagon $A_s$, for $s \in [1/4,1]$.

{\footnotesize
\begin{eqnarray}
\nonumber
\left[\matrix{315\cr 105\cr 105}\right]\left[\matrix{-315\cr 105\cr 105}\right]\left[\matrix{-315\cr -105\cr 105}\right]\left[\matrix{315\cr -105\cr 105}\right]\left[\matrix{420\cr 0\cr 105}\right]\left[\matrix{-420\cr 0\cr 105}\right]\left[\matrix{0\cr 420\cr 420}\right]\left[\matrix{420\cr 0\cr 420}\right]\left[\matrix{-420\cr 0\cr 420}\right]\left[\matrix{0\cr -420\cr 420}\right]
\end{eqnarray}
\/}

Here is ${\cal B\/}[1/4,1]$.  This set intersects the fiber $X_s$ in
the triangle $\beta_s$, for $s \in [1/4,1]$.

{\footnotesize
\begin{eqnarray}
\nonumber
\left[\matrix{-315\cr -105\cr 105}\right]\left[\matrix{-420\cr 0\cr 105}\right]\left[\matrix{-525\cr -105\cr 105}\right]\left[\matrix{-840\cr -420\cr 420}\right]\left[\matrix{-420\cr 0\cr 420}\right]\left[\matrix{0\cr -420\cr 420}\right]
\end{eqnarray}
\/}

Here is ${\cal P\/}[1/4,1/2]$.  This set intersects the fiber $X_s$ in
the pentagon $P_s$, for $s \in [1/4,1/2]$.

{\footnotesize
\begin{eqnarray}
\nonumber
\left[\matrix{-315\cr -105\cr 105}\right]\left[\matrix{-105\cr -105\cr 105}\right]\left[\matrix{-105\cr 105\cr 105}\right]\left[\matrix{-315\cr 105\cr 105}\right]\left[\matrix{-630\cr -210\cr 210}\right]\left[\matrix{-210\cr -210\cr 210}\right]\left[\matrix{-210\cr 210\cr 210}\right]
\end{eqnarray}
\/}

Here is ${\cal Q\/}[1/4,1/2]$.  This set intersects the fiber $X_s$ in
the triangle $Q_s$, for $s \in [1/4,1/2]$.

{\footnotesize
\begin{eqnarray}
\nonumber
\left[\matrix{-525\cr -105\cr 105}\right]\left[\matrix{-315\cr -105\cr 105}\right]\left[\matrix{-315\cr 105\cr 105}\right]\left[\matrix{-630\cr -210\cr 210}\right]
\end{eqnarray}
\/}

Here is the period $2$ tile $\tau[1,5/4]$ from Calculation 6.
{\footnotesize
\begin{eqnarray}
\nonumber
\left[\matrix{-525\cr -315\cr 315}\right]\left[\matrix{-315\cr -315\cr 315}\right]\left[\matrix{-315\cr -105\cr 315}\right]\left[\matrix{-525\cr -105\cr 315}\right]\left[\matrix{-420\cr -420\cr 420}\right]
\end{eqnarray}
\/}

Here is the period $2$ tile $\tau[3/4,1]$ from Calculation 8.
{\footnotesize
\begin{eqnarray}
\nonumber
\left[\matrix{-420\cr -420\cr 420}\right]\left[\matrix{-210\cr -630\cr 630}\right]\left[\matrix{-210\cr -210\cr 630}\right]\left[\matrix{-630\cr -210\cr 630}\right]\left[\matrix{-630\cr -630\cr 630}\right]
\end{eqnarray}
\/}

Here is the domain ${\cal Z\/}[1,5/4]$ from Calculation 6.
{\footnotesize
\begin{eqnarray}
\nonumber
\left[\matrix{-420\cr -420\cr 420}\right]\left[\matrix{-840\cr -420\cr 420}\right]\left[\matrix{-945\cr -525\cr 525}\right]\left[\matrix{-525\cr -525\cr 525}\right]\left[\matrix{-525\cr -315\cr 525}\right]\left[\matrix{-735\cr -315\cr 525}\right]
\end{eqnarray}
\/}

Here is the domain ${\cal Z\/}[3/4,1]$ from Calculation 8.
{\footnotesize
\begin{eqnarray}
\nonumber
\left[\matrix{-420\cr -420\cr 420}\right]\left[\matrix{-840\cr -420\cr 420}\right]\left[\matrix{-945\cr -525\cr 525}\right]\left[\matrix{-525\cr -525\cr 525}\right]\left[\matrix{-525\cr -315\cr 525}\right]\left[\matrix{-735\cr -315\cr 525}\right]
\end{eqnarray}
\/}

\subsection{The Polyhedra in the Partition}
\label{partition}

Define
\begin{equation}
\iota_1(x,y,s)=(-x,-y,s), \hskip 30 pt
\iota_2(x,y,s)=\bigg(\frac{x+y}{2s},\frac{x-y}{2s},\frac{1}{2s}\bigg).
\end{equation}

The partition of ${\cal X\/}[1/4,1/2]$ consists of the $19$
polyhedra
\begin{equation}
\label{Alist}
\alpha_0,\alpha_1,...,\alpha_9,\iota_1(\alpha_1),...,\iota_1(\alpha_9).
\end{equation}

The partition of ${\cal X\/}[1/2,1]$ consists of the $13$
polyhedra
\begin{equation}
\beta_0,\beta_1,...,\beta_6,\iota_1(\beta_1),...,\iota_1(\beta_6).
\end{equation}

The partition of  ${\cal X\/}[1,2]$
consists of the $19$ polyhedra
\begin{equation}
\iota_2 \circ F(\alpha_i), \hskip 30 pt i=0,...,18.
\end{equation}
Here $\alpha_0,...,\alpha_{18}$ are the polyhedra from Equation \ref{Alist}.

As we mentioned in the last chapter, we scale each polyhedron
by a factor of $420$ so that all the entries are integers.  Here
are the polyhedra.

{\footnotesize
\begin{eqnarray}
\nonumber
A0=\left[\matrix{105\cr 105\cr 105}\right]\left[\matrix{-105\cr 105\cr 105}\right]\left[\matrix{-105\cr -105\cr 105}\right]\left[\matrix{105\cr -105\cr 105}\right]\left[\matrix{210\cr 210\cr 210}\right]\left[\matrix{-210\cr 210\cr 210}\right]\left[\matrix{-210\cr -210\cr 210}\right]\left[\matrix{210\cr -210\cr 210}\right]
\end{eqnarray}

\begin{eqnarray}
\nonumber
A1=\left[\matrix{420\cr 0\cr 105}\right]\left[\matrix{525\cr 105\cr 105}\right]\left[\matrix{315\cr 105\cr 105}\right]\left[\matrix{280\cr 140\cr 140}\right]
\end{eqnarray}

\begin{eqnarray}
\nonumber
A2=\left[\matrix{105\cr 105\cr 105}\right]\left[\matrix{280\cr 140\cr 140}\right]\left[\matrix{140\cr 140\cr 140}\right]\left[\matrix{210\cr -210\cr 210}\right]
\end{eqnarray}

\begin{eqnarray}
\nonumber
A3=\left[\matrix{280\cr 140\cr 140}\right]\left[\matrix{140\cr 140\cr 140}\right]\left[\matrix{210\cr 210\cr 210}\right]\left[\matrix{210\cr -210\cr 210}\right]\left[\matrix{420\cr 0\cr 210}\right]
\end{eqnarray}

\begin{eqnarray}
\nonumber
A4=\left[\matrix{525\cr 105\cr 105}\right]\left[\matrix{420\cr 0\cr 140}\right]\left[\matrix{420\cr 0\cr 210}\right]\left[\matrix{630\cr 210\cr 210}\right]\left[\matrix{210\cr 210\cr 210}\right]
\end{eqnarray}

\begin{eqnarray}
\nonumber
A5=\left[\matrix{315\cr -105\cr 105}\right]\left[\matrix{420\cr 0\cr 105}\right]\left[\matrix{315\cr 105\cr 105}\right]\left[\matrix{420\cr 0\cr 120}\right]\left[\matrix{378\cr -42\cr 126}\right]
\end{eqnarray}

\begin{eqnarray}
\nonumber
A6=\left[\matrix{420\cr 0\cr 120}\right]\left[\matrix{462\cr 42\cr 126}\right]\left[\matrix{420\cr 0\cr 140}\right]\left[\matrix{420\cr 140\cr 140}\right]\left[\matrix{280\cr 140\cr 140}\right]\left[\matrix{210\cr 210\cr 210}\right]
\end{eqnarray}

\begin{eqnarray}
\nonumber
A7=\left[\matrix{315\cr 105\cr 105}\right]\left[\matrix{420\cr 0\cr 120}\right]\left[\matrix{378\cr -42\cr 126}\right]\left[\matrix{420\cr 0\cr 140}\right]
\end{eqnarray}

\begin{eqnarray}
\nonumber
A8=\left[\matrix{315\cr -105\cr 105}\right]\left[\matrix{315\cr 105\cr 105}\right]\left[\matrix{105\cr 105\cr 105}\right]\left[\matrix{105\cr -105\cr 105}\right]\left[\matrix{420\cr 0\cr 140}\right]\left[\matrix{280\cr 140\cr 140}\right]\left[\matrix{210\cr -210\cr 210}\right]\left[\matrix{420\cr 0\cr 210}\right]
\end{eqnarray}

\begin{eqnarray}
\nonumber
A9=\left[\matrix{420\cr 0\cr 105}\right]\left[\matrix{525\cr 105\cr 105}\right]\left[\matrix{420\cr 0\cr 120}\right]\left[\matrix{462\cr 42\cr 126}\right]\left[\matrix{420\cr 140\cr 140}\right]\left[\matrix{280\cr 140\cr 140}\right]
\end{eqnarray}

\begin{eqnarray}
\nonumber
B0=\left[\matrix{105\cr 105\cr 105}\right]\left[\matrix{-105\cr 105\cr 105}\right]\left[\matrix{-105\cr -105\cr 105}\right]\left[\matrix{105\cr -105\cr 105}\right]\left[\matrix{210\cr 210\cr 210}\right]\left[\matrix{-210\cr 210\cr 210}\right]\left[\matrix{-210\cr -210\cr 210}\right]\left[\matrix{210\cr -210\cr 210}\right]
\end{eqnarray}

\begin{eqnarray}
\nonumber
B1=\left[\matrix{420\cr 0\cr 105}\right]\left[\matrix{525\cr 105\cr 105}\right]\left[\matrix{315\cr 105\cr 105}\right]\left[\matrix{280\cr 140\cr 140}\right]
\end{eqnarray}

\begin{eqnarray}
\nonumber
B2=\left[\matrix{105\cr 105\cr 105}\right]\left[\matrix{280\cr 140\cr 140}\right]\left[\matrix{140\cr 140\cr 140}\right]\left[\matrix{210\cr -210\cr 210}\right]
\end{eqnarray}

\begin{eqnarray}
\nonumber
B3=\left[\matrix{280\cr 140\cr 140}\right]\left[\matrix{140\cr 140\cr 140}\right]\left[\matrix{210\cr 210\cr 210}\right]\left[\matrix{210\cr -210\cr 210}\right]\left[\matrix{420\cr 0\cr 210}\right]
\end{eqnarray}

\begin{eqnarray}
\nonumber
B4=\left[\matrix{525\cr 105\cr 105}\right]\left[\matrix{420\cr 0\cr 140}\right]\left[\matrix{420\cr 0\cr 210}\right]\left[\matrix{630\cr 210\cr 210}\right]\left[\matrix{210\cr 210\cr 210}\right]
\end{eqnarray}

\begin{eqnarray}
\nonumber
B5=\left[\matrix{315\cr -105\cr 105}\right]\left[\matrix{420\cr 0\cr 105}\right]\left[\matrix{315\cr 105\cr 105}\right]\left[\matrix{420\cr 0\cr 120}\right]\left[\matrix{378\cr -42\cr 126}\right]
\end{eqnarray}

\begin{eqnarray}
\nonumber
B6=\left[\matrix{420\cr 0\cr 120}\right]\left[\matrix{462\cr 42\cr 126}\right]\left[\matrix{420\cr 0\cr 140}\right]\left[\matrix{420\cr 140\cr 140}\right]\left[\matrix{280\cr 140\cr 140}\right]\left[\matrix{210\cr 210\cr 210}\right]
\end{eqnarray}
\/}

\subsection{The Action of the Map}
\label{move}

In this section we explain the action of the map on each of
the polyhedra listed above.  To each polyhedron
we associate a $4$-tuple of integers.  The list
$V=(u_1,v_1,u_2,v_2)$ tells us that
\begin{equation}
\label{affine}
F_V\left[\matrix{x\cr y \cr s}\right]=
\left[\matrix{1&0&2v_1-2v_2\cr 0&1&2v_1+2v_2\cr 0&0&1}\right]
\left[\matrix{x\cr y \cr s}\right]+
\left[\matrix{-2u_1\cr 2u_2\cr 0}\right].
\end{equation}

\noindent
{\bf Remark;\/}
Equation \ref{affine} gives the equation for the
action on the unscaled polyhedra.  When we acts on
the scaled polyhedra listed above, we need to scale
the translation part of the map by $420$.  That is,
$-2u_1$ and $2u_2$ need to be replaced by
$-840u_1$ and $840u_2$.
\newline

The polyhedra $\alpha_0$ and $\beta_0$ correspond to the
trivial tiles. The vectors associated to these
are $a_0=b_0=(0,0,0,0)$.  Below we will list
the vectors $a_1,...,a_9,b_1,...,b_6$.
We have the relations
\begin{equation}
a_{9+i}=-a_i, \hskip 30 pt b_{6+i}=-b_i.
\end{equation}
The vector $c_i$ associated to $\iota_2 \circ F(\alpha_i)$ is
given by the following rule:
\begin{equation}
a_i=(u_1,v_1,u_2,v_2) \hskip 10 pt
\Longrightarrow \hskip 10 pt
c_i=(-v_2,-u_2,-v_1,-u_1).
\end{equation}

Recall that the map $F$ is really the composition
$(F')^2$, where $F'$ maps the bundle
${\cal X\/}[1/4,2]$ to the polyhedron obtained by
rotating ${\cal X\/}[1/4,2]$ by $90$ degrees about the $z$-axis.
To get the action of $F'$ we simply replace
each vector $V=a_1,a_2,...$ by $V'$, where
\begin{equation}
V=(u_1,v_1,u_2,v_2) \hskip 10 pt
\Longrightarrow \hskip 10 pt
V'=(u_1,v_1,0,0).
\end{equation}

Here are the vectors.
$$
a_1=(1,2,0,-2), \hskip 20 pt
a_2=(0,-1,-1,-2) \hskip 20 pt
a_3=(0,-1,-1,-1).
$$
$$
a_4=(1,1,0,-1), \hskip 20 pt
a_5=(0,-2,-1,-2) \hskip 20 pt
a_6=(1,2,1,1).
$$
$$
a_7=(0,-2,-1,-1), \hskip 20 pt
a_8=(0,-1,0,1) \hskip 20 pt
a_9=(1,2,1,2).
$$
$$
b_1=(1,0,-1,-1), \hskip 20 pt
b_2=(1,1,0,-1) \hskip 20 pt
b_3=(1,1,1,0).
$$
$$
b_4=(0,-1,-1,0), \hskip 20 pt
b_5=(0,-1,-1,-1) \hskip 20 pt
b_6=(1,1,1,1).
$$

\newpage
\section{References}

[{\bf AG\/}] A. Goetz and G. Poggiaspalla, {\it Rotations by $\pi/7$\/}, Nonlinearity {\bf 17\/}
(2004) no. 5 1787-1802
\newline
\newline
[{\bf AKT\/}] R. Adler, B. Kitchens, and C. Tresser,
{\it Dynamics of non-ergodic piecewise affine maps of the torus\/},
Ergodic Theory Dyn. Syst {\bf 21\/} (2001) no. 4 959-999
\newline
\newline
[{\bf BKS\/}] T. Bedford, M. Keane, and C. Series, eds.,
{\it Ergodic Theory, Symbolic Dynamics, and Hyperbolic Spaces\/}, Oxford University Press, Oxford (1991).
\newline
\newline
[{\bf B\/}] J, Buzzi, {\it Piecewise isometrries have zero topological entropy (English summary)\/} Ergodic Theory and Dynamical Systems {\bf 21\/}  (2001) no. 5 pp 1371-1377
\newline
\newline
[{\bf GH1\/}] E Gutkin and N. Haydn, {\it Topological entropy of generalized polygon exchanges\/}, Bull. Amer. Math. Soc., {\bf 32\/} (1995) no. 1., pp 50-56
\newline
\newline
[{\bf GH2\/}] E Gutkin and N. Haydn, {\it Topological entropy polygon exchange transformations and polygonal billiards\/}, Ergodic Theory and Dynamical Systems {\bf 17\/} (1997) no. 4., pp 849-867
\newline
\newline
[{\bf H\/}] H. Haller, {\it Rectangle Exchange Transformations\/}, Monatsh Math. {\bf 91\/}
(1985) 215-232
\newline
\newline
[{\bf Hoo\/}] W. Patrick Hooper, {\it Renormalization of Polygon Exchage Maps arising from Corner Percolation\/} Invent. Math. 2012.
\newline
\newline
[{\bf LKV\/}] J. H. Lowenstein, K. L. Koupsov, F. Vivaldi, {\it Recursive Tiling and Geometry of piecewise rotations by $\pi/7$\/}, nonlinearity {\bf 17\/} (2004) no. 2.
[{\bf Low\/}] J. H. Lowenstein, {\it Aperiodic orbits of piecewise rational
rotations of convex polygons with recursive tiling\/}, Dyn. Syst. {\bf 22\/}
(2007) no. 1 25-63
\newline
\newline
[{\bf R\/}] G. Rauzy, {\it Exchanges d'intervalles et transformations induites\/},
Acta. Arith. {\bf 34\/} 315-328 (1979)
\newline
\newline
[{\bf S0\/}] R.E. Schwartz {\it The Octagonal PET II: Topology of the Limit Sets\/},
preprint  (2012)
\newline
\newline
[{\bf S1\/}] R.E. Schwartz {\it Outer Billiards, Quarter Turn Compositions, and Polytope Exchange 
Transformations\/}, preprint (2011)
\newline
\newline
[{\bf S2\/}] R. E. Schwartz, {\it Outer Billiards on Kites\/},
Annals of Math Studies {\bf 171\/}, Princeton University Press (2009)
\newline
\newline
[{\bf S3\/}] R. E. Schwartz, {\it Outer Billiards on the Penrose Kite:
Compactification and Renormalization\/}, Journal of Modern Dynamics, 2012.
\newline
\newline
[{\bf T\/}] S. Tabachnikov, {\it Billiards\/}, Soci\'{e}t\'{e} Math\'{e}matique de France, 
``Panoramas et Syntheses'' 1, 1995
\newline
\newline
[{\bf VL\/}] F. Vivaldi and J. H. Lowenstein, {|it Arithmetical properties of a family
of irrational piecewise rotations\/}, {\it Nonlinearity\/} {\bf 19\/}:1069--1097 (2007).
\newline
\newline
[{\bf Y\/}] J.-C. Yoccoz, {\small {\it Continued Fraction Algorithms for Interval Exchange Maps: An Introduction\/}\/}, Frontiers in Number Theory, Physics, and Geometry Vol 1, P. Cartier, B. Julia, P. Moussa, P. Vanhove (editors) Springer-Verlag 4030437 (2006)
\newline
\newline
[{\bf Z\/}] A. Zorich, {\it Flat Surfaces\/}, Frontiers in Number Theory, Physics, and Geometry Vol 1, P. Cartier, B. Julia, P. Moussa, P. Vanhove (editors) Springer-Verlag 4030437 (2006)

\end{document}